# Infinite-dimensional Log-Determinant divergences II: Alpha-Beta divergences

Hà Quang Minh

*Istituto Italiano di Tecnologia, Via Morego 30, Genova 16163, ITALY*

**Abstract**

This work presents a parametrized family of divergences, namely Alpha-Beta Log-Determinant (Log-Det) divergences, between positive definite unitized trace class operators on a Hilbert space. This is a generalization of the Alpha-Beta Log-Determinant divergences between symmetric, positive definite matrices to the infinite-dimensional setting. The family of Alpha-Beta Log-Det divergences is highly general and contains many divergences as special cases, including the recently formulated infinite-dimensional affine-invariant Riemannian distance and the infinite-dimensional Alpha Log-Det divergences between positive definite unitized trace class operators. In particular, it includes a parametrized family of metrics between positive definite trace class operators, with the affine-invariant Riemannian distance and the square root of the symmetric Stein divergence being special cases. For the Alpha-Beta Log-Det divergences between covariance operators on a Reproducing Kernel Hilbert Space (RKHS), we obtain closed form formulas via the corresponding Gram matrices.

*Email address:* minh.haquang@iit.it (Hà Quang Minh)

# 1. Introduction

Symmetric Positive Definite (SPD) matrices play an important role in many areas of mathematics, statistics, machine learning, optimization, computer vision, and related fields, see e.g. [1, 2, 3, 4, 5, 6, 7, 8, 9, 10, 11]. The set $\text{Sym}^{++}(n)$ of $n \times n$ SPD matrices is an open convex cone and can also be equipped with a Riemannian manifold structure. Among the most studied Riemannian metrics on $\text{Sym}^{++}(n)$ are the classical affine-invariant metric [1, 2, 3, 5, 12] and the more recent Log-Euclidean metric [4, 9, 13]. The convex cone structure of $\text{Sym}^{++}(n)$, on the other hand, gives rise to distance-like functions such as the Alpha Log-Determinant divergences [14], which have been shown to be special cases of the Alpha-Beta Log-Determinant divergences [15]. These divergences are fast to compute and have been shown to work well in various applications [7, 16, 8]. The present work aims to generalize the Alpha-Beta Log-Determinant divergences to the infinite-dimensional setting.

**Finite-dimensional Alpha-Beta Log-Determinant divergences**. We recall that for $A, B \in \text{Sym}^{++}(n)$, the Alpha-Beta Log-Determinant (Log-Det) divergence between $A$ and $B$ is a parametrized family of divergences defined by (see [15])

$$D^{(\alpha,\beta)}(A,B) = \frac{1}{\alpha\beta} \log \det \left[ \frac{\alpha(AB^{-1})^\beta + \beta(AB^{-1})^{-\alpha}}{\alpha + \beta} \right], \quad (1)$$
$$\alpha \neq 0, \beta \neq 0, \alpha + \beta \neq 0.$$

*Remark* 1. To keep our presentation compact, in the following we consider the case $\alpha > 0$, $\beta > 0$, as well as the limiting cases $\alpha = 0, \beta = 0$. Since $D^{(\alpha,\beta)}(A,B) = D^{(-\alpha,-\beta)}(B,A)$, the case $\alpha < 0, \beta < 0$ is essentially identical to the previous case. We do not consider the cases $\alpha, \beta$ have opposite signs, since in those cases the well-definedness and finiteness of $D_r^{(\alpha,\beta)}(A,B)$ depends on the spectrum of $AB^{-1}$ (see Theorem 2 in [15]), that is it is not a valid divergence on all of $\text{Sym}^{++}(n)$.

The parametrized family of divergences defined by Eq.(1) is highly general and admits as special cases many metrics and distance-like functions on $\text{Sym}^{++}(n)$, including in particular the following:

1. The affine-invariant Riemannian distance [3], corresponding to the limiting case $D^{(0,0)}(A,B)$, with

$$D^{(0,0)}(A,B) = \frac{1}{2}d_{\text{aiE}}^2(A,B) = \frac{1}{2}||\log(B^{-1/2}AB^{-1/2})||_F^2, \qquad (2)$$

where $\log(A)$ denotes the principal logarithm of the matrix $A$ and $||\ ||_F$ denotes the Frobenius norm.

2. The Alpha Log-Determinant divergences [14], corresponding to $D^{(\alpha,1-\alpha)}(A,B)$, $0 < \alpha < 1$, with

$$D^{(\alpha,1-\alpha)}(A,B) = \frac{1}{\alpha(1-\alpha)}\log\left[\frac{\det[\alpha A + (1-\alpha)B]}{\det(A)^\alpha \det(B)^{1-\alpha}}\right]. \qquad (3)$$

A special case of this divergence is the symmetric Stein divergence (also called the Jensen-Bregman LogDet divergence), corresponding to $D^{(1/2,1/2)}(A,B)$, whose square root is a metric on $\text{Sym}^{++}(n)$ [16], with

$$D^{(1/2,1/2)}(A,B) = 4d_{\text{stein}}^2(A,B) = 4\log\frac{\det(\frac{A+B}{2})}{\sqrt{\det(A)\det(B)}}. \qquad (4)$$

3. The limiting cases $\beta = 0$ and $\alpha = 0$ correspond to, respectively,

$$D^{(\alpha,0)}(A,B) = \frac{1}{\alpha^2}\left\{\text{tr}((A^{-1}B)^\alpha - I) - \alpha\log\det(A^{-1}B)\right\}, \qquad (5)$$

$$D^{(0,\beta)}(A,B) = \frac{1}{\beta^2}\left\{\text{tr}((B^{-1}A)^\beta - I) - \beta\log\det(B^{-1}A)\right\}, \qquad (6)$$

with $D^{(1,0)}(A,B) = \text{tr}(A^{-1}B - I) - \log\det(A^{-1}B)$ and $D^{(0,1)}(A,B) = \text{tr}(B^{-1}A - I) - \log\det(B^{-1}A)$.

**Contributions of this work**. The current work is a continuation and generalization of the author's recent work [17]. In [17], we generalized the Alpha Log-Det divergences between SPD matrices [14] to the infinite-dimensional Alpha Log-Determinant divergences between positive definite unitized trace class operators in a Hilbert space. In the current work, we present a formulation for the Alpha-Beta Log-Det divergences between positive definite unitized trace class operators, generalizing the Alpha-Beta divergences between SPD matrices as defined by Eq.(1). As in the finite-dimensional setting, the formulation we present here is general and admits as special cases many metrics and distance-like functions between positive definite unitized trace class operators, including in particular the following: the infinite-dimensional affine-invariant Riemannian distance [18]; the infinite-dimensional Alpha Log-Det divergences [17], a special case of which is the infinite-dimensional symmetric Stein divergence. For the divergences between reproducing kernel Hilbert spaces (RKHS) covariance operators, we obtain closed form formulas for the Alpha-Beta Log-Det divergences via the corresponding Gram matrices.

**Organization**. We provide a summary of the main results of the paper in Section 2, including our definition of the infinite-dimensional Alpha-Beta Log-Det divergences. The key concepts involved are described in Section 3. The motivations and derivations leading to our definition of the Alpha-Beta Log-Det divergences are presented in Section 4. We then show in Section 5 that both the affine-invariant Riemannian distance and the Alpha Log-Det divergences are special cases of the Alpha-Beta Log-Det divergences. All mathematical proofs are presented in Appendix A.

## 2. Summary of main results

We present a summary of our main results in this section, with the detailed technical descriptions provided in subsequent sections. Throughout the paper, let $\mathcal{H}$ denote a separable Hilbert space, with $\dim(\mathcal{H}) = \infty$, unless explicitly stated otherwise. Let $\mathcal{L}(\mathcal{H})$ be the Banach space of bounded linear operators on $\mathcal{H}$ and $\mathrm{Sym}(\mathcal{H}) \subset \mathcal{L}(\mathcal{H})$ be the subspace of self-adjoint, bounded operators on $\mathcal{H}$. For $A \in \mathcal{L}(\mathcal{H})$, we write $A > 0$ to denote that $A$ is a self-adjoint positive definite operator. Let $\mathrm{Tr}(\mathcal{H})$ denote the Banach algebra of trace class operators on $\mathcal{H}$. The set of *positive definite unitized trace class operators* on $\mathcal{H}$ is then defined to be

$$\mathrm{PTr}(\mathcal{H}) = \{A + \gamma I > 0 \ : \ A = A^*, A \in \mathrm{Tr}(\mathcal{H}), \gamma \in \mathbb{R}\}. \tag{7}$$

The main purpose of the current work is the generalization of the Alpha-Beta Log-Det divergence between SPD matrices, as defined in Eq. (1), to that between positive definite unitized trace class operators in $\mathrm{PTr}(\mathcal{H})$. The following is our definition of the Alpha-Beta (Log-Det) divergences in the infinite-dimensional setting.

**Definition 1** (**Alpha-Beta Log-Determinant Divergences**). *Assume that* $\dim(\mathcal{H}) = \infty$. *Let* $\alpha > 0$, $\beta > 0$ *be fixed. Let* $r \in \mathbb{R}$, $r \neq 0$ *be fixed. For* $(A + \gamma I), (B + \mu I) \in \mathrm{PTr}(\mathcal{H})$, *the* $(\alpha, \beta)$-*Log-Det divergence* $D_r^{(\alpha,\beta)}[(A + \gamma I), (B + \mu I)]$ *is defined to be*

$$D_r^{(\alpha,\beta)}[(A + \gamma I), (B + \mu I)]$$
$$= \frac{1}{\alpha\beta} \log \left[ \left(\frac{\gamma}{\mu}\right)^{r(\delta - \frac{\alpha}{\alpha+\beta})} \det{}_\mathrm{X} \left( \frac{\alpha(\Lambda + \frac{\gamma}{\mu}I)^{r(1-\delta)} + \beta(\Lambda + \frac{\gamma}{\mu}I)^{-r\delta}}{\alpha + \beta} \right) \right], \quad (8)$$

*where* $\Lambda + \frac{\gamma}{\mu}I = (B + \mu I)^{-1/2}(A + \gamma I)(B + \mu I)^{-1/2}$, $\delta = \frac{\alpha\gamma^r}{\alpha\gamma^r + \beta\mu^r}$. *Equivalently,*

$$D_r^{(\alpha,\beta)}[(A + \gamma I), (B + \mu I)]$$
$$= \frac{1}{\alpha\beta} \log \left[ \left(\frac{\gamma}{\mu}\right)^{r(\delta - \frac{\alpha}{\alpha+\beta})} \det{}_\mathrm{X} \left( \frac{\alpha(Z + \frac{\gamma}{\mu}I)^{r(1-\delta)} + \beta(Z + \frac{\gamma}{\mu}I)^{-r\delta}}{\alpha + \beta} \right) \right], \quad (9)$$

*where* $Z + \frac{\gamma}{\mu}I = (A + \gamma I)(B + \mu I)^{-1}$.

*Remark* 2. In Definition 1, $\det_\mathrm{X}$ denotes the *extended Fredholm determinant* defined in [17] (see Section 3 below). For $\gamma = 1$, we have $\det_\mathrm{X}(A + \gamma I) = \det(A + I)$, with $\det$ on the right hand side being the Fredholm determinant. For $\dim(\mathcal{H}) < \infty$, $\det_\mathrm{X}(A + \gamma I) = \det(A + \gamma I)$, with $\det$ on the right hand side being the standard matrix determinant.

The quantity $D_r^{(\alpha,\beta)}[(A + \gamma I), (B + \mu I)]$ where $\alpha > 0, \beta > 0$, as stated in Definition 1, can be extended to the cases $\alpha > 0, \beta = 0$ and $\alpha = 0, \beta > 0$, $\forall r \in \mathbb{R}, r \neq 0$, via limiting arguments. The following is our definition in these cases.

**Definition 2** (**Limiting cases - I**). *Assume that* $\dim(\mathcal{H}) = \infty$. *Let* $\alpha > 0, \beta > 0$, $r \neq 0$ *be fixed. For* $(A + \gamma I), (B + \mu I) \in \mathrm{PTr}(\mathcal{H})$, *the Log-Det divergence* $D_r^{(\alpha,0)}[(A + \gamma I), (B + \mu I)]$ *is defined to be*

$$D_r^{(\alpha,0)}[(A + \gamma I), (B + \mu I)] = \frac{r}{\alpha^2} \left[ \left(\frac{\mu}{\gamma}\right)^r - 1 \right] \log \frac{\mu}{\gamma} \quad (10)$$
$$+ \frac{1}{\alpha^2} \mathrm{tr}_X([(A + \gamma I)^{-1}(B + \mu I)]^r - I)$$
$$- \frac{1}{\alpha^2} \left(\frac{\mu}{\gamma}\right)^r \log \det{}_\mathrm{X}[(A + \gamma I)^{-1}(B + \mu I)]^r.$$

Similarly, $D_r^{(0,\beta)}[(A+\gamma I), (B+\mu I)]$ is defined to be

$$D_r^{(0,\beta)}[(A+\gamma I), (B+\mu I)] = \frac{r}{\beta^2}\left[\left(\frac{\gamma}{\mu}\right)^r - 1\right]\log\frac{\gamma}{\mu} \quad (11)$$
$$+ \frac{1}{\beta^2}\text{tr}_X([(B+\mu I)^{-1}(A+\gamma I)]^r - I)$$
$$- \frac{1}{\beta^2}\left(\frac{\gamma}{\mu}\right)^r \log\det{}_X[(B+\mu I)^{-1}(A+\gamma I)]^r.$$

The following result confirms that the quantity $D_r^{(\alpha,\beta)}$, as defined in Definitions 1 and 2, is in fact a divergence on $\text{PTr}(\mathcal{H})$.

**Theorem 1 (Positivity).** *Assume the hypothesis stated in Definitions 1 and 2. Then*

$$D_r^{(\alpha,\beta)}[(A+\gamma I), (B+\mu I)] \geq 0 \quad (12)$$
$$D_r^{(\alpha,\beta)}[(A+\gamma I), (B+\mu I)] = 0 \iff A = B, \ \gamma = \mu. \quad (13)$$

**Theorem 2 (Special cases - I).** *The following are some of the most important special cases of Definitions 1 and 2.*

1. *The infinite-dimensional affine-invariant Riemannian distance $d_{\text{aiHS}}[(A+\gamma I), (B+\mu I)]$ [18], which corresponds to the limiting case $\lim_{\alpha\to 0} D_r^{(\alpha,\alpha)}[(A+\gamma I), (B+\mu I)]$, where $r = r(\alpha)$ is smooth, with $r(0) = 0$, $r'(0) \neq 0$, and $r(\alpha) \neq 0$ for $\alpha \neq 0$. The limit is given by*

$$\lim_{\alpha\to 0} D_r^{(\alpha,\alpha)}[(A+\gamma I), (B+\mu I)] = \frac{[r'(0)]^2}{8}d_{\text{aiHS}}^2[(A+\gamma I), (B+\mu I)]. \quad (14)$$

   *In particular, for $r = 2\alpha$,*

$$\lim_{\alpha\to 0} D_{2\alpha}^{(\alpha,\alpha)}[(A+\gamma I), (B+\mu I)] = \frac{1}{2}d_{\text{aiHS}}^2[(A+\gamma I), (B+\mu I)]. \quad (15)$$

   *This is the content of Theorem 9.*

2. *The infinite-dimensional Alpha Log-Determinant divergences $d_{\text{logdet}}^\alpha[(A+\gamma I), (B+\mu I)]$ [17], with*

$$D_{\pm 1}^{(\alpha,1-\alpha)}[(A+\gamma I), (B+\mu I)] = d_{\text{logdet}}^{\pm(1-2\alpha)}[(A+\gamma I), (B+\mu I)], \quad (16)$$
$$0 \leq \alpha \leq 1.$$

   *This is the content of Theorem 10.*

Since the limit $\lim_{\alpha \to 0} D_r^{(\alpha,\alpha)}[(A+\gamma I),(B+\mu I)]$ in the first part of Theorem 2 is unique, up to the multiplicative factor $[r'(0)]^2/8$, we define the quantity $D_0^{(0,0)}[(A+\gamma I),(B+\mu I)]$ as follows.

**Definition 3** (**Limiting cases - II**). *For $(A+\gamma I), (B+\mu I) \in \mathrm{PTr}(\mathcal{H})$, the Log-Det divergence $D_0^{(0,0)}[(A+\gamma I),(B+\mu I)]$ is defined to be*

$$D_0^{(0,0)}[(A+\gamma I),(B+\mu I)] = \lim_{\alpha \to 0} D_{2\alpha}^{(\alpha,\alpha)}[(A+\gamma I),(B+\mu I)]$$
$$= \frac{1}{2} d_{\mathrm{aiHS}}^2[(A+\gamma I),(B+\mu I)]. \tag{17}$$

Since $d_{\mathrm{aiHS}}[(A+\gamma I),(B+\mu I)]$ is a metric on $\mathrm{PTr}(\mathcal{H})$, $D_0^{(0,0)}[(A+\gamma I),(B+\mu I)]$ is automatically a symmetric divergence on $\mathrm{PTr}(\mathcal{H})$. In fact, it is a member of the parametrized family $D_{2\alpha}^{(\alpha,\alpha)}[(A+\gamma I),(B+\mu I)]$, $\alpha \geq 0$, of symmetric divergences on $\mathrm{PTr}(\mathcal{H})$, as stated in the following result.

**Theorem 3** (**Special cases - II**). *The parametrized family $D_{2\alpha}^{(\alpha,\alpha)}[(A+\gamma I),(B+\mu I)]$, $\alpha \geq 0$, is a family of symmetric divergences on $\mathrm{PTr}(\mathcal{H})$, with $\alpha = 0$ corresponding to the infinite-dimensional affine-invariant Riemannian distance above and $\alpha = 1/2$ corresponding to the infinite-dimensional symmetric Stein divergence, which is given by $\frac{1}{4} d_{\mathrm{logdet}}^0[(A+\gamma I),(B+\mu I)]$.*

**Finite-dimensional case**. For $\gamma = \mu$, we have $\delta = \frac{\alpha}{\alpha+\beta}$, so that Eq. (9) becomes

$$D_r^{(\alpha,\beta)}[(A+\gamma I),(B+\gamma I)]$$
$$= \frac{1}{\alpha\beta} \log \det\nolimits_X \left( \frac{\alpha[(A+\gamma I)(B+\gamma I)^{-1}]^{\frac{r\beta}{\alpha+\beta}} + \beta[(A+\gamma I)(B+\gamma I)^{-1}]^{-\frac{r\alpha}{\alpha+\beta}}}{\alpha+\beta} \right). \tag{18}$$

In the finite-dimensional case, where $A$ and $B$ are two $n \times n$ SPD matrices, setting $\gamma = 0$ and recalling that $\det_X = \det$ for finite matrices, we obtain

$$D_r^{(\alpha,\beta)}(A,B) = \frac{1}{\alpha\beta} \log \det \left( \frac{\alpha(AB^{-1})^{\frac{r\beta}{\alpha+\beta}} + \beta(AB^{-1})^{-\frac{r\alpha}{\alpha+\beta}}}{\alpha+\beta} \right). \tag{19}$$

In particular, by setting $r = \alpha + \beta$, we recover Eq. (1). For $\gamma = \mu$, Eq. (10) becomes

$$D_r^{(\alpha,0)}[(A+\gamma I),(B+\gamma I)] \tag{20}$$
$$= \frac{1}{\alpha^2} \left\{ \mathrm{tr}_X([(A+\gamma I)^{-1}(B+\gamma I)]^r - I) - \log \det\nolimits_X[(A+\gamma I)^{-1}(B+\gamma I)]^r \right\},$$

which reduces to Eq. (5) when $A, B \in \text{Sym}^{++}(n)$, $\gamma = 0$, and $r = \alpha$. Similarly, Eq. (11) becomes

$$D_r^{(0,\beta)}[(A+\gamma I), (B+\gamma I)] \quad (21)$$
$$= \frac{1}{\beta^2}\left\{\text{tr}_X([(B+\gamma I)^{-1}(A+\gamma I)]^r - I) - \log\det_X[(B+\gamma I)^{-1}(A+\gamma I)]^r\right\},$$

which reduces to Eq. (6) when $A, B \in \text{Sym}^{++}(n)$, $\gamma = 0$, and $r = \beta$.

*Remark* 3. As in the cases of the Log-Hilbert-Schmidt distance [19], the infinite-dimensional affine-invariant Riemannian distance [18, 20], and the infinite-dimensional Alpha Log-Det divergences [17], we show below that in general, the infinite-dimensional formulation is *not* obtainable as the limit of the finite-dimensional version as the dimension approaches infinity.

*Remark* 4. Except for the case $r = \alpha + \beta$, the quantity $r$ in $D_r^{(\alpha,\beta)}$ that we introduce here, to the best of our knowledge, has no equivalence in the existing literature in the finite-dimensional setting.

*Remark* 5. Throughout the paper, we employ the following notations. Using the identity $(B+\mu I)^{-1} = \frac{1}{\mu}I - \frac{B}{\mu}(B+\mu I)^{-1}$, we write the operator $(B+\mu I)^{-1/2}(A+\gamma I)(B+\mu I)^{-1/2}$ as

$$(B+\mu I)^{-1/2}(A+\gamma I)(B+\mu I)^{-1/2} = \Lambda + \frac{\gamma}{\mu}I \in \text{PTr}(\mathcal{H}), \quad (22)$$

where $\Lambda = (B+\mu I)^{-1/2}A(B+\mu I)^{-1/2} - \frac{\gamma}{\mu}B(B+\mu I)^{-1} \in \text{Tr}(\mathcal{H})$. This notation is employed in Eq. (8). Similarly, in Eq. (9), we write

$$(A+\gamma I)(B+\mu I)^{-1} = \frac{\gamma}{\mu}I + A(B+\mu I)^{-1} - \frac{\gamma}{\mu}B(B+\mu I)^{-1} = Z + \frac{\gamma}{\mu}I, \quad (23)$$

where $Z = A(B+\mu I)^{-1} - \frac{\gamma}{\mu}B(B+\mu I)^{-1} \in \text{Tr}(\mathcal{H})$.

**Metric properties**. Consider now a special case, where $\alpha = \beta$ and $r = \alpha + \beta$. For simplicity, we consider operators $(A+\gamma I)$ and $(B+\mu I)$ with $\gamma = \mu$. For $\gamma > 0, \gamma \in \mathbb{R}$ fixed, we define the following subset of $\text{PTr}(\mathcal{H})$

$$\text{PTr}(\mathcal{H})(\gamma) = \{A + \gamma I > 0 : A^* = A, A \in \text{Tr}(\mathcal{H})\}. \quad (24)$$

*Remark* 6. Throughout the paper, we assume, unless stated otherwise, that $\dim(\mathcal{H}) = \infty$, and the condition $A + \gamma I > 0$ automatically implies that $\gamma > 0$. When $\dim(\mathcal{H}) < \infty$, we can set $\gamma = 0$.

**Theorem 4** (**Metric property**). *Let $\gamma > 0, \gamma \in \mathbb{R}$ be fixed. The square root function $\sqrt{D_{2\alpha}^{(\alpha,\alpha)}[(A+\gamma I),(B+\gamma I)]}$ is a metric on $\mathrm{PTr}(\mathcal{H})(\gamma)$ for all $\alpha \geq 0$.*

We thus have a family of metrics between positive definite operators of the form $(A+\gamma I) \in \mathrm{PTr}(\mathcal{H})(\gamma)$, parametrized by the parameter $\alpha \geq 0$. In particular, with $\alpha = 0$ in Theorem 4, we obtain the affine-invariant Riemannian distance, and with $\alpha = \frac{1}{2}$ we obtain the following metric, which is the square root of the infinite-dimensional Stein divergence

$$\sqrt{D_1^{(1/2,1/2)}[(A+\gamma I),(B+\gamma I)]} = 2\sqrt{\log\left[\frac{\det_X\left[\frac{(A+\gamma I)+(B+\gamma I)}{2}\right]}{\det_X(A+\gamma I)^{1/2}\det_X(B+\gamma I)^{1/2}}\right]}. \tag{25}$$

The corresponding finite-dimensional result [15], where $A, B \in \mathrm{Sym}^{++}(n)$, is recovered by setting $\gamma = 0$ in Theorem 4. In particular, with $\alpha = 1/2$ and $A, B \in \mathrm{Sym}^{++}(n)$, we obtain the corresponding result of [16].

*Remark* 7. The analysis of $\sqrt{D_{2\alpha}^{(\alpha,\alpha)}[(A+\gamma I),(B+\mu I)]}$, where $\gamma \neq \mu$, is technically more involved and will be presented in a separate work.

## 3. Positive definite unitized trace class operators

To generalize the Alpha-Beta Log-Determinant divergences from the finite to infinite-dimensional setting, we need to employ the following concepts

- Positive definite operators $\mathbb{P}(\mathcal{H})$.
- Extended (or unitized) trace class operators $\mathrm{Tr}_X(\mathcal{H})$.
- Positive definite unitized trace class operators $\mathrm{PTr}(\mathcal{H})$.
- Extended Fredholm determinant $\det_X$ on $\mathrm{Tr}_X(\mathcal{H})$.

- Exponential, logarithm, and power functions for operators in $\mathrm{PTr}(\mathcal{H})$ and their products.

We discuss in detail below the logarithm and power functions of products of operators in $\mathrm{PTr}(\mathcal{H})$. Other concepts are briefly reviewed and we refer to [17] for the detailed motivations leading to the definitions of these concepts. Throughout the following, we assume that $\dim(\mathcal{H}) = \infty$, unless stated explicitly otherwise.

**Positive definite operators**. We recall that an operator $\mathcal{A} \in \mathcal{L}(\mathcal{H})$ is said to be positive definite if there exists a constant $M_A > 0$ such that

$$\langle x, Ax \rangle \geq M_A ||x||^2 \quad \forall x \in \mathcal{H}.$$

This is equivalent to saying that $A$ is both strictly positive and invertible. We denote by $\mathbb{P}(\mathcal{H})$ the set of all positive definite operators on $\mathcal{H}$.

**Extended trace class operators**. Let $\mathrm{Tr}(\mathcal{H})$ denote the set of trace class operators on $\mathcal{H}$, the set of extended (or unitized) trace class operators on $\mathcal{H}$ is defined to be

$$\mathrm{Tr}_X(\mathcal{H}) = \{A + \gamma I \; : \; A \in \mathrm{Tr}(\mathcal{H}), \gamma \in \mathbb{R}\}.$$

Equipped with the extended trace class norm

$$||A + \gamma I||_{\mathrm{tr}_X} = ||A||_{\mathrm{tr}} + |\gamma| = \mathrm{tr}|A| + |\gamma|,$$

$\mathrm{Tr}_X(\mathcal{H})$ becomes a Banach algebra. For $(A + \gamma I) \in \mathrm{Tr}_X(\mathcal{H})$, its *extended trace* is defined to be

$$\mathrm{tr}_X(A + \gamma I) = \mathrm{tr}(A) + \gamma.$$

Thus by this definition $\mathrm{tr}_X(I) = 1$, in contrast to usual trace definition, according to which $\mathrm{tr}(I) = \infty$.

**Extended Fredholm determinant**. For $(A + \gamma I) \in \mathrm{Tr}_X(\mathcal{H}), \gamma \neq 0$, its extended Fredholm determinant is defined to be

$$\det{}_X(A + \gamma I) = \frac{1}{\gamma}\det\left(\frac{A}{\gamma} + I\right),$$

where the determinant on the right hand side is the Fredholm determinant. For $\gamma = 1$, we recover the Fredholm determinant. In the case $\dim(\mathcal{H}) < \infty$, we define $\det_X(A + \gamma I) = \det(A + \gamma I)$, the standard matrix determinant.

**Positive definite unitized trace class operators**. Having defined both positive definite operators and extended trace class operators, the set of positive definite unitized trace class operators $\text{PTr}(\mathcal{H}) \subset \text{Tr}_X(\mathcal{H})$ is then defined to be the intersection

$$\text{PTr}(\mathcal{H}) = \text{Sym}(\mathcal{H}) \cap \mathbb{P}(\mathcal{H}) = \{A + \gamma I > 0 \ : A^* = A, \ A \in \text{Tr}(\mathcal{H}) \ \gamma \in \mathbb{R}\}.$$

**Exponential, logarithm, and power functions**. Consider the exponential function $\exp : \mathcal{L}(\mathcal{H}) \to \mathcal{L}(\mathcal{H})$ defined by

$$\exp(A) = \sum_{j=0}^{\infty} \frac{A^j}{j!}.$$

The following result shows that $\exp$ maps $\text{Tr}_X(\mathcal{H})$ to $\text{Tr}_X(\mathcal{H})$.

**Lemma 1.** *Let* $(A + \gamma I) \in \text{Tr}_X(\mathcal{H})$. *Then* $\exp(A + \gamma I) \in \text{Tr}_X(\mathcal{H})$.

Consider next the inverse function $\log = \exp^{-1} : \mathcal{L}(\mathcal{H}) \to \mathcal{L}(\mathcal{H})$. For any $(A + \gamma I) \in \text{PTr}(\mathcal{H})$, $\log(A + \gamma I)$ is always well-defined as follows. Let $\{\lambda_k\}_{k=1}^{\infty}$ be the eigenvalues of $A$ with corresponding orthonormal eigenvectors $\{\phi_k\}_{k=1}^{\infty}$. Then

$$A = \sum_{k=1}^{\infty} \lambda_k \phi_k \otimes \phi_k, \quad \log(A + \gamma I) = \sum_{k=1}^{\infty} \log(\lambda_k + \gamma) \phi_k \otimes \phi_k, \qquad (26)$$

where $\phi_k \otimes \phi_k : \mathcal{H} \to \mathcal{H}$ is a rank-one operator defined by $(\phi_k \otimes \phi_k) w = \langle \phi_k, w \rangle \phi_k$ $\forall w \in \mathcal{H}$. Moreover, $\log(A + \gamma I) \in \text{Sym}(\mathcal{H}) \cap \text{Tr}_X(\mathcal{H})$ and assumes the form

$$\log(A + \gamma I) = A_1 + \gamma_1 I, \quad A_1 \in \text{Sym}(\mathcal{H}) \cap \text{Tr}(\mathcal{H}), \gamma_1 \in \mathbb{R}.$$

By Proposition 6 in [17], for any $\alpha \in \mathbb{R}$, the power function $(A + \gamma I)^\alpha$ is then well-defined via the expression

$$(A + \gamma I)^\alpha = \exp[\alpha \log(A + \gamma I)] \in \text{PTr}(\mathcal{H}).$$

For the purposes of the current work, we need to go beyond the set $\text{PTr}(\mathcal{H})$. Specifically, for two operators $(A + \gamma I), (B + \mu I) \in \text{PTr}(\mathcal{H})$, we show that

$$\log[(A + \gamma I)(B + \mu I)^{-1}], \quad [(A + \gamma I)(B + \mu I)^{-1}]^\alpha, \alpha \in \mathbb{R} \qquad (27)$$

are all well-defined and are elements of $\text{Tr}_X(\mathcal{H})$, even though they are no longer necessarily self-adjoint.

First, let $B \in \mathcal{L}(\mathcal{H})$ be any invertible operator, then for any $A \in \mathcal{L}(\mathcal{H})$, we have

$$\exp(BAB^{-1}) = \sum_{j=0}^{\infty} \frac{(BAB^{-1})^j}{j!} = B\left(\sum_{j=0}^{\infty} \frac{A^j}{j!}\right) B^{-1} = B\exp(A)B^{-1}.$$

Thus for $(A + \gamma I) \in \mathrm{PTr}(\mathcal{H})$, the logarithm of $B(A + \gamma I)B^{-1} = BAB^{-1} + \gamma I \in \mathrm{Tr}_X(\mathcal{H})$ is also well-defined and is given by

$$\log[B(A+\gamma I)B^{-1}] = B\log(A+\gamma I)B^{-1}$$
$$= B(A_1 + \gamma_1 I)B^{-1} = BA_1B^{-1} + \gamma_1 I \in \mathrm{Tr}_X(\mathcal{H}). \quad (28)$$

Using Eq. (28), we obtain the following results.

**Proposition 1.** *Let $(A+\gamma I), (B+\mu I) \in \mathrm{PTr}(\mathcal{H})$. Let $\Lambda + \frac{\gamma}{\mu}I = (B+\mu I)^{-1/2}(A+\gamma I)(B+\mu I)^{-1/2}$. Then*

1. *The logarithm $\log[(A+\gamma I)(B+\mu I)^{-1}] \in \mathrm{Tr}_X(\mathcal{H})$ is well-defined and is given by*

$$\log[(A+\gamma I)(B+\mu I)^{-1}] = (B+\mu I)^{1/2} \log\left(\Lambda + \frac{\gamma}{\mu}I\right)(B+\mu I)^{-1/2}.$$
(29)

2. *For any $\alpha \in \mathbb{R}$, the power function $[(A + \gamma I)(B + \mu I)^{-1}]^{\alpha} \in \mathrm{Tr}_X(\mathcal{H})$ is well-defined and is given by*

$$[(A+\gamma I)(B+\mu I)^{-1}]^{\alpha} = (B+\mu I)^{1/2}\left(\Lambda + \frac{\gamma}{\mu}I\right)^{\alpha}(B+\mu I)^{-1/2}. \quad (30)$$

3. *For any $p, q \in \mathbb{R}$, any $\alpha, \beta \in \mathbb{R}$ such that $\alpha + \beta \neq 0$,*

$$\mathrm{det}_X\left[\frac{\alpha[(A+\gamma I)(B+\mu I)^{-1}]^p + \beta[(A+\gamma I)(B+\mu I)^{-1}]^q}{\alpha + \beta}\right]$$
$$= \mathrm{det}_X\left[\frac{\alpha(\Lambda + \frac{\gamma}{\mu}I)^p + \beta(\Lambda + \frac{\gamma}{\mu}I)^q}{\alpha + \beta}\right]. \quad (31)$$

## 4. Infinite-Dimensional Alpha-Beta Log-Determinant divergences

We now show the motivations and derivations leading to Definition 1. We recall that in the case $\dim(\mathcal{H}) < \infty$, the Log-Det divergences were motivated by Ky Fan's inequality [21] on the log-concavity of the determinant, which states that for $A, B \in \text{Sym}^{++}(n)$, $\det(\alpha A + (1-\alpha)B) \geq \det(A)^\alpha \det(B)^{1-\alpha}$, $0 \leq \alpha \leq 1$, with equality if and only if $A = B$ ($0 < \alpha < 1$). This inequality has recently been generalized to the infinite-dimensional setting for the extended Fredholm determinant (Theorem 1 in [17]). The following is a further generalization of Theorem 1 in [17].

**Theorem 5.** *Let $0 \leq \alpha \leq 1$. For $(A + \gamma I), (B + \mu I) \in \text{PTr}(\mathcal{H})$, for any $p, q \in \mathbb{R}$,*

$$\det_X[\alpha(A+\gamma I)^p + (1-\alpha)(B+\mu I)^q]$$
$$\geq \left(\frac{\gamma^p}{\mu^q}\right)^{\alpha-\delta} \det_X(A+\gamma I)^{p\delta} \det_X(B+\mu I)^{q(1-\delta)}, \quad (32)$$

*where $\delta = \frac{\alpha\gamma^p}{\alpha\gamma^p + (1-\alpha)\mu^q}$, $1-\delta = \frac{(1-\alpha)\mu^q}{\alpha\gamma^p + (1-\alpha)\mu^q}$. For $0 < \alpha < 1$, equality happens if and only if*

$$\left(\frac{A}{\gamma}+I\right)^p = \left(\frac{B}{\mu}+I\right)^q \text{ and } \gamma^p = \mu^q \iff (A+\gamma I)^p = (B+\mu I)^q. \quad (33)$$

*In particular, for $\gamma = \mu \neq 1$, equality happens if and only if simultaneously*

$$p = q \text{ and } A = B. \quad (34)$$

In particular, for $p = q = 1$, we recover Theorem 1 in [17]. From Theorem 5, we immediately have the following result.

**Corollary 1.** *Let $\alpha > 0$, $\beta > 0$. For $(A+\gamma I), (B+\mu I) \in \text{PTr}(\mathcal{H})$, for any $p, q \in \mathbb{R}$,*

$$\det_X\left[\frac{\alpha(A+\gamma I)^p + \beta(B+\mu I)^q}{\alpha+\beta}\right]$$
$$\geq \left(\frac{\gamma^p}{\mu^q}\right)^{\frac{\alpha}{\alpha+\beta}-\delta} \det_X(A+\gamma I)^{p\delta} \det_X(B+\mu I)^{q(1-\delta)}, \quad (35)$$

*where $\delta = \frac{\alpha\gamma^p}{\alpha\gamma^p+\beta\mu^q}$, $1-\delta = \frac{\beta\mu^q}{\alpha\gamma^p+\beta\mu^q}$. Equality happens if and only if $(A+\gamma I)^p = (B+\mu I)^q$. For $\gamma = \mu \neq 1$, equality happens if and only if simultaneously $p = q$ and $A = B$.*

Motivated by Theorem 5 and Corollary 1, we first define the following quantity.

**Definition 4.** *Let $\alpha > 0$, $\beta > 0$ be fixed. For $(A + \gamma I), (B + \mu I) \in \mathrm{PTr}(\mathcal{H})$, for $p, q \in \mathbb{R}$, define*

$$D_{(p,q)}^{(\alpha,\beta)}[(A + \gamma I), (B + \mu I)]$$
$$= \frac{1}{\alpha\beta} \log \left[ \left(\frac{\gamma}{\mu}\right)^{(p+q)(\delta - \frac{\alpha}{\alpha+\beta})} \det{}_\mathrm{X} \left( \frac{\alpha(\Lambda + \frac{\gamma}{\mu}I)^p + \beta(\Lambda + \frac{\gamma}{\mu}I)^{-q}}{\alpha + \beta} \right) \right], \quad (36)$$

*where $\Lambda + \frac{\gamma}{\mu}I = (B + \mu I)^{-1/2}(A + \gamma I)(B + \mu I)^{-1/2}$, $\delta = \frac{\alpha(\frac{\gamma}{\mu})^{p+q}}{\alpha(\frac{\gamma}{\mu})^{p+q}+\beta}$.*

The following theorem gives sufficient conditions for $p, q \in \mathbb{R}$, with $\alpha > 0, \beta > 0$ being fixed, so that for a given pair of operators $(A + \gamma I), (B + \mu I) \in \mathrm{PTr}(\mathcal{H})$, the quantity $D_{(p,q)}^{(\alpha,\beta)}[(A + \gamma I), (B + \mu I)]$ in Definition 4 is nonnegative, with equality if and only if $A = B$ and $\gamma = \mu$.

**Theorem 6.** *Let $\alpha > 0$, $\beta > 0$ be fixed. For $(A + \gamma I), (B + \mu I) \in \mathrm{PTr}(\mathcal{H})$, assume that $p, q \in \mathbb{R}$ satisfy the following conditions*

$$p + q \neq 0, \quad (37)$$

$$\alpha p \left(\frac{\gamma}{\mu}\right)^{p+q} = \beta q. \quad (38)$$

*Then the quantity $D_{(p,q)}^{(\alpha,\beta)}[(A + \gamma I), (B + \mu I)]$ satisfies*

$$D_{(p,q)}^{(\alpha,\beta)}[(A + \gamma I), (B + \mu I)] \geq 0, \quad (39)$$

$$D_{(p,q)}^{(\alpha,\beta)}[(A + \gamma I), (B + \mu I)] = 0 \iff A = B, \gamma = \mu. \quad (40)$$

Subsequently, we assume that conditions (37) and (38) are satisfied. We see that $p$ and $q$ are not uniquely determined by (38). One way to enforce the uniqueness of $p$ and $q$ is by fixing the sum $p + q$. This is the approach we adopt in this work, which leads to Definition 1.

**Theorem 7.** *Under the hypothesis of Theorem 6, assume further that $p + q = r$, $r \in \mathbb{R}, r \neq 0$, $r$ fixed. Under this condition, in Definition 4, we have*

$$\delta = \frac{\alpha(\frac{\gamma}{\mu})^r}{\alpha(\frac{\gamma}{\mu})^r + \beta}, \quad p = r(1 - \delta) = \frac{\beta r}{\alpha(\frac{\gamma}{\mu})^r + \beta}, \quad q = r\delta = \frac{\alpha r(\frac{\gamma}{\mu})^r}{\alpha(\frac{\gamma}{\mu})^r + \beta}. \quad (41)$$

*Plugging the expressions for p and q in Eq. (41) into Definition 4, we obtain Definition 1. Furthermore, the two formulas given in Eqs. (8) and (9) in Definition 1 are equivalent.*

We now show how $D^{(\alpha,\beta)}_{(p,q)}[(A+\gamma I),(B+\mu I)]$ can be expressed concretely in terms of the Fredholm determinant.

**Theorem 8.** *Let $\alpha > 0$, $\beta > 0$ be fixed. For $(A+\gamma I), (B+\mu I) \in \mathrm{PTr}(\mathcal{H})$, assume that $p, q \in \mathbb{R}$ satisfy conditions (37) and (38) in Theorem 6. Then*

$$D^{(\alpha,\beta)}_{(p,q)}[(A+\gamma I),(B+\mu I)] = \frac{(p+q)(\delta - \frac{\alpha}{\alpha+\beta})}{\alpha\beta}\left(\log\frac{\gamma}{\mu}\right) \qquad (42)$$
$$+ \frac{1}{\alpha\beta}\log\left(\frac{\alpha(\frac{\gamma}{\mu})^p + \beta(\frac{\gamma}{\mu})^{-q}}{\alpha+\beta}\right) + \frac{1}{\alpha\beta}\log\det\left[\frac{\alpha(\Lambda+\frac{\gamma}{\mu}I)^p + \beta(\Lambda+\frac{\gamma}{\mu}I)^{-q}}{\alpha(\frac{\gamma}{\mu})^p + \beta(\frac{\gamma}{\mu})^{-q}}\right].$$

## 5. Special cases of the Alpha-Beta Log-Determinant divergences

We now describe several important special cases of Definition 1, including the infinite-dimensional affine-invariant Riemannian distance, the infinite-dimensional Alpha Log-Det divergences [17], and the infinite-dimensional Beta Log-Det divergences.

*5.1. Affine-invariant Riemannian distance*

Let $\mathrm{HS}(\mathcal{H})$ denote the space of Hilbert-Schmidt operators on $\mathcal{H}$, which is defined by

$$\mathrm{HS}(\mathcal{H}) = \{A \in \mathcal{L}(\mathcal{H}) \;:\; ||A||^2_{\mathrm{HS}} = \mathrm{tr}(A^*A) < \infty\},$$

where $||\;||_{\mathrm{HS}}$ is the Hilbert-Schmidt norm. If $A$ is Hilbert-Schmidt, then $A$ is compact and possesses a countable set of eigenvalues $\{\lambda_k\}_{k=1}^{\infty}$. If $A$ is furthermore self-adjoint, then the Hilbert-Schmidt norm of $A$ is given by

$$||A||^2_{\mathrm{HS}} = \sum_{k=1}^{\infty}\lambda_k^2.$$

We recall the infinite-dimensional Hilbert manifold of positive definite unitized Hilbert-Schmidt operators on $\mathcal{H}$, considered in [18]

$$\Sigma(\mathcal{H}) = \{A+\gamma I > 0 \;:\; A = A^*, A \in \mathrm{HS}(\mathcal{H}), \gamma \in \mathbb{R}\}.$$

In the case $\dim(\mathcal{H}) = \infty$, the set $\mathrm{PTr}(\mathcal{H})$ of positive definite unitized trace class operators on $\mathcal{H}$ is a strict subset of $\Sigma(\mathcal{H})$. The manifold $\Sigma(\mathcal{H})$ can be equipped with the following Riemannian metric, as formulated by [18]. For each $P \in \Sigma(\mathcal{H})$, on the tangent space $T_P(\Sigma(\mathcal{H})) \cong \mathcal{H}_\mathbb{R} = \{A + \gamma I \;:\; A = A^*, A \in \mathrm{HS}(\mathcal{H}), \gamma \in \mathbb{R}\}$, we define the following inner product

$$\langle A + \gamma I, B + \mu I \rangle_P = \langle P^{-1/2}(A + \gamma I) P^{-1/2}, P^{-1/2}(B + \mu I) P^{-1/2} \rangle_{\mathrm{eHS}},$$

where $\langle \,,\, \rangle_{\mathrm{eHS}}$ is the extended Hilbert-Schmidt inner product, defined by

$$\langle A + \gamma I, B + \mu I \rangle_{\mathrm{eHS}} = \langle A, B \rangle_{\mathrm{HS}} + \gamma \mu.$$

The Riemannian metric given by $\langle \,,\, \rangle_P$ then makes $\Sigma(\mathcal{H})$ an infinite-dimensional Riemannian manifold. Under this metric, the geodesic distance between $(A+\gamma I), (B+\mu I)$ is given by

$$d_{\mathrm{aiHS}}[(A + \gamma I), (B + \mu I)] = ||\log[(B + \mu I)^{-1/2}(A + \gamma I)(B + \mu I)^{-1/2}]||_{\mathrm{eHS}}. \tag{43}$$

We now show that the affine-invariant distance $d_{\mathrm{aiHS}}[(A + \gamma I), (B + \mu I)]$ is a limiting case of $D_r^{(\alpha,\beta)}[(A + \gamma I), (B + \mu I)]$, as $\alpha \to 0, \beta \to 0$. In this section, we consider $\beta = \alpha$, in which case Definition 1 reduces to the following.

**Definition 5.** *In Definition 1, with $\alpha = \beta$, we have*

$$D_r^{(\alpha,\alpha)}[(A + \gamma I), (B + \mu I)]$$
$$= \frac{1}{\alpha^2} \log \left[ \left(\frac{\gamma}{\mu}\right)^{r(\delta - \frac{1}{2})} \det\nolimits_{\mathrm{X}} \left( \frac{(\Lambda + \frac{\gamma}{\mu} I)^{r(1-\delta)} + (\Lambda + \frac{\gamma}{\mu} I)^{-r\delta}}{2} \right) \right], \tag{44}$$

*where $\delta = \frac{(\frac{\gamma}{\mu})^r}{(\frac{\gamma}{\mu})^r + 1}$, $1 - \delta = \frac{1}{(\frac{\gamma}{\mu})^r + 1}$.*

By Theorem 8, we have the following formula, which expresses $D_r^{(\alpha,\alpha)}[(A + \gamma I), (B + \mu I)]$ concretely in terms of the Fredholm determinant.

$$D_r^{(\alpha,\alpha)}[(A + \gamma I), (B + \mu I)] = \frac{r(\delta - \frac{1}{2})}{\alpha^2} \log\left(\frac{\gamma}{\mu}\right) + \frac{1}{\alpha^2} \log\left(\frac{(\frac{\gamma}{\mu})^p + (\frac{\gamma}{\mu})^{-q}}{2}\right)$$
$$+ \frac{1}{\alpha^2} \log \det \left[ \frac{(\Lambda + \frac{\gamma}{\mu} I)^p + (\Lambda + \frac{\gamma}{\mu} I)^{-q}}{(\frac{\gamma}{\mu})^p + (\frac{\gamma}{\mu})^{-q}} \right], \tag{45}$$

where $\delta = \frac{(\frac{\gamma}{\mu})^r}{(\frac{\gamma}{\mu})^r+1}$, $1-\delta = \frac{1}{(\frac{\gamma}{\mu})^r+1}$, $p = r(1-\delta)$, $q = r\delta$.

The following is the main result in this section.

**Theorem 9 (Affine-Invariant Riemannian Distance).** *Let $(A + \gamma I), (B + \mu I) \in \mathrm{PTr}(\mathcal{H})$. Assume that $r = r(\alpha)$ is smooth, with $r(0) = 0$, $r'(0) \neq 0$, and $r(\alpha) \neq 0$ for $\alpha \neq 0$. Then*

$$\lim_{\alpha \to 0} D_r^{(\alpha,\alpha)}[(A+\gamma I),(B+\mu I)] = \frac{[r'(0)]^2}{8} d_{\mathrm{aiHS}}^2[(A+\gamma I),(B+\mu I)]. \qquad (46)$$

*In particular, for $r = 2\alpha$, we have*

$$\lim_{\alpha \to 0} D_{2\alpha}^{(\alpha,\alpha)}[(A+\gamma I),(B+\mu I)] = \frac{1}{2} d_{\mathrm{aiHS}}^2[(A+\gamma I),(B+\mu I)]. \qquad (47)$$

*Remark* 8. We stress that, as they are currently stated, the limits in Theorem 9 are valid for $(A+\gamma I), (B+\mu I) \in \mathrm{PTr}(\mathcal{H})$, that is $A$ and $B$ must be trace class operators. The generalization of Theorem 9 to the entire Hilbert manifold $\Sigma(\mathcal{H})$, where $A$ and $B$ are Hilbert-Schmidt operators, will be presented in an upcoming work.

*5.2. Infinite-dimensional Alpha Log-Determinant divergences*

We now show that the formulation for the infinite-dimensional Alpha Log-Determinant divergences in [17] is a special case of the present formulation, with $\beta = 1 - \alpha$ and $r = \pm 1$. Let $\dim(\mathcal{H}) = \infty$. We recall that for $-1 < \alpha < 1$, the Log-Det $\alpha$-divergence $d_{\mathrm{logdet}}^\alpha[(A+\gamma I),(B+\mu I)]$ for $(A+\gamma I),(B+\mu I) \in \mathrm{PTr}(\mathcal{H})$ is defined in [17] to be

$$d_{\mathrm{logdet}}^\alpha[(A+\gamma I),(B+\mu I)]$$
$$= \frac{4}{1-\alpha^2} \log \left[ \frac{\det_X\left(\frac{1-\alpha}{2}(A+\gamma I) + \frac{1+\alpha}{2}(B+\mu I)\right)}{\det_X(A+\gamma I)^q \det_X(B+\mu I)^{1-q}} \left(\frac{\gamma}{\mu}\right)^{q-\frac{1-\alpha}{2}} \right], \qquad (48)$$

where $q = \frac{(1-\alpha)\gamma}{(1-\alpha)\gamma+(1+\alpha)\mu}$ and $1-q = \frac{(1+\alpha)\mu}{(1-\alpha)\gamma+(1+\alpha)\mu}$, with the limiting cases $\alpha = \pm 1$ given by

$$d_{\mathrm{logdet}}^1[(A+\gamma I),(B+\mu I)] = \left(\frac{\gamma}{\mu}-1\right)\log\frac{\gamma}{\mu} + \mathrm{tr}_X[(B+\mu I)^{-1}(A+\gamma I) - I]$$
$$- \frac{\gamma}{\mu} \log \det_X[(B+\mu I)^{-1}(A+\gamma I)]. \qquad (49)$$

$$d_{\mathrm{logdet}}^{-1}[(A+\gamma I),(B+\mu I)] = \left(\frac{\mu}{\gamma}-1\right)\log\frac{\mu}{\gamma} + \mathrm{tr}_X\left[(A+\gamma I)^{-1}(B+\mu I) - I\right]$$
$$- \frac{\mu}{\gamma} \log \det_X[(A+\gamma I)^{-1}(B+\mu I)]. \qquad (50)$$

**Definition 6.** *In Definition 1, with $0 < \alpha < 1$ and $\beta = 1 - \alpha$, we have*

$$D_r^{(\alpha,1-\alpha)}[(A+\gamma I),(B+\mu I)] \tag{51}$$
$$= \frac{1}{\alpha(1-\alpha)}\log\left[\left(\frac{\gamma}{\mu}\right)^{r(\delta-\alpha)} \det\nolimits_X\left(\alpha\left(\Lambda+\frac{\gamma}{\mu}I\right)^{r(1-\delta)}+(1-\alpha)\left(\Lambda+\frac{\gamma}{\mu}I\right)^{-r\delta}\right)\right].$$

*where $\delta = \frac{\alpha(\frac{\gamma}{\mu})^r}{\alpha(\frac{\gamma}{\mu})^r+1-\alpha}$, $1-\delta = \frac{1-\alpha}{\alpha(\frac{\gamma}{\mu})^r+1-\alpha}$.*

The following result shows that $D_r^{(\alpha,1-\alpha)}[(A+\gamma I),(B+\mu I)]$ for the cases $r = \pm 1$ are precisely $d_{\text{logdet}}^{1-2\alpha}[(A+\gamma I),(B+\mu I)]$ and $d_{\text{logdet}}^{2\alpha-1}[(A+\gamma I),(B+\mu I)]$, respectively.

**Theorem 10 (Alpha Log-Determinant Divergences).** *Let $0 < \alpha < 1$ be fixed. For $(A+\gamma I), (B+\mu I) \in \text{PTr}(\mathcal{H})$,*

$$\begin{aligned}&D_1^{(\alpha,1-\alpha)}[(A+\gamma I),(B+\mu I)]\\&=\frac{\delta-\alpha}{\alpha(1-\alpha)}\log\frac{\gamma}{\mu}+\frac{1}{\alpha(1-\alpha)}\log\left[\frac{\det_X[\alpha(A+\gamma I)+(1-\alpha)(B+\mu I)]}{\det_X(A+\gamma I)^\delta \det_X(B+\mu I)^{1-\delta}}\right]\\&=d_{\text{logdet}}^{1-2\alpha}[(A+\gamma I),(B+\mu I)],\end{aligned} \tag{52}$$

*where $\delta = \frac{\alpha\gamma}{\alpha\gamma+(1-\alpha)\mu}$. Similarly,*

$$D_{-1}^{(\alpha,1-\alpha)}[(A+\gamma I),(B+\mu I)] = d_{\text{logdet}}^{2\alpha-1}[(A+\gamma I),(B+\mu I)]. \tag{53}$$

*At the endpoints $\alpha = 0$ and $\alpha = 1$,*

$$\lim_{\alpha\to 1}D_1^{(\alpha,1-\alpha)}[(A+\gamma I),(B+\mu I)] = d_{\text{logdet}}^{-1}[(A+\gamma I),(B+\mu I)] \tag{54}$$

$$\lim_{\alpha\to 0}D_1^{(\alpha,1-\alpha)}[(A+\gamma I),(B+\mu I)] = d_{\text{logdet}}^{1}[(A+\gamma I),(B+\mu I)]. \tag{55}$$

In particular, in Theorem 10, for $\gamma = \mu$, we have $\delta = \alpha$, and

$$\begin{aligned}&D_1^{(\alpha,1-\alpha)}[(A+\gamma I),(B+\gamma I)]\\&=\frac{1}{\alpha(1-\alpha)}\log\left[\frac{\det_X[\alpha(A+\gamma I)+(1-\alpha)(B+\gamma I)]}{\det_X(A+\gamma I)^\alpha \det_X(B+\gamma I)^{1-\alpha}}\right].\end{aligned} \tag{56}$$

This is the direct generalization of the finite-dimensional formula given by Eq. (6) in [14].

*Remark* 9 (**Beta Log-Determinant Divergences**). In the finite-dimensional setting in [15], the authors call $D^{1,\beta}(A,B)$ the Beta Log-Determinant divergence between $A, B \in \text{Sym}^{++}(n)$. Similarly, in the case $\dim(\mathcal{H}) = \infty$, let $\beta > 0$ be fixed and let $r \in \mathbb{R}$, $r \neq 0$ be fixed. For $(A + \gamma I), (B + \mu I) \in \text{PTr}(\mathcal{H})$, we then have the corresponding infinite-dimensional Beta Log-Determinant divergence

$$D_r^{(1,\beta)}[(A + \gamma I), (B + \mu I)]$$
$$= \frac{1}{\beta} \log \left[ \left( \frac{\gamma}{\mu} \right)^{r(\delta - \frac{1}{1+\beta})} \det\nolimits_X \left( \frac{(\Lambda + \frac{\gamma}{\mu} I)^{r(1-\delta)} + \beta (\Lambda + \frac{\gamma}{\mu} I)^{-r\delta}}{1 + \beta} \right) \right], \quad (57)$$

where $\Lambda + \frac{\gamma}{\mu} I = (B + \mu I)^{-1/2}(A + \gamma I)(B + \mu I)^{-1/2}$, $\delta = \frac{(\frac{\gamma}{\mu})^r}{(\frac{\gamma}{\mu})^r + \beta}$, $1 - \delta = \frac{\beta}{(\frac{\gamma}{\mu})^r + \beta}$. However, we do not explore this divergence in detail in this work.

*5.3. Other limiting cases*

We consider next two other limiting cases, namely $\beta \to 0$ when $\alpha > 0$ is fixed, and $\alpha \to 0$ when $\beta > 0$ is fixed. In particular, our definitions of $D_r^{(\alpha,0)}[(A + \gamma I), (B + \mu I)]$, $\alpha > 0$, and $D_r^{(0,\beta)}[(A + \gamma I), (B + \mu I)]$, $\beta > 0$, as given in Definition 2, are based on the respective limits in Theorems 11 and 12 below.

**Theorem 11 (Liming case $\alpha > 0, \beta \to 0$).** *Let $\alpha > 0$ be fixed. Assume that $r = r(\beta)$ is smooth, with $r(0) = r(\beta = 0)$. Then*

$$\lim_{\beta \to 0} D_r^{(\alpha,\beta)}[(A + \gamma I), (B + \mu I)] = \frac{r(0)}{\alpha^2} \left[ \left( \frac{\mu}{\gamma} \right)^{r(0)} - 1 \right] \log \frac{\mu}{\gamma} \quad (58)$$
$$+ \frac{1}{\alpha^2} \text{tr}_X([(A + \gamma I)^{-1}(B + \mu I)]^{r(0)} - I)$$
$$- \frac{1}{\alpha^2} \left( \frac{\mu}{\gamma} \right)^{r(0)} \log \det\nolimits_X[(A + \gamma I)^{-1}(B + \mu I)]^{r(0)}.$$

**Theorem 12 (Limit case $\alpha \to 0, \beta > 0$).** *Let $\beta > 0$ be fixed. Assume that $r = r(\alpha)$ is smooth, with $r(0) = r(\alpha = 0)$. Then*

$$\lim_{\alpha \to 0} D_r^{(\alpha,\beta)}[(A + \gamma I), (B + \mu I)] = \frac{r(0)}{\beta^2} \left[ \left( \frac{\gamma}{\mu} \right)^{r(0)} - 1 \right] \log \frac{\gamma}{\mu} \quad (59)$$
$$+ \frac{1}{\beta^2} \text{tr}_X([(B + \mu I)^{-1}(A + \gamma I)]^{r(0)} - I)$$
$$- \frac{1}{\beta^2} \left( \frac{\gamma}{\mu} \right)^{r(0)} \log \det\nolimits_X[(B + \mu I)^{-1}(A + \gamma I)]^{r(0)}.$$

**Special cases.** Let us now describe several special cases of Theorems 11 and 12, including their specialization to the finite-dimensional setting.

(i) For $\gamma = \mu$, we have

$$\lim_{\beta \to 0} D_r^{(\alpha,\beta)}[(A + \gamma I), (B + \gamma I)] = \frac{1}{\alpha^2} \text{tr}_X([(A + \gamma I)^{-1}(B + \gamma I)]^{r(0)} - I)$$
$$- \frac{1}{\alpha^2} \log\det_X[(A + \gamma I)^{-1}(B + \gamma I)]^{r(0)}, \quad (60)$$

$$\lim_{\alpha \to 0} D_r^{(\alpha,\beta)}[(A + \gamma I), (B + \gamma I)] = \frac{1}{\beta^2} \text{tr}_X([(B + \gamma I)^{-1}(A + \gamma I)]^{r(0)} - I)$$
$$- \frac{1}{\beta^2} \log\det_X[(B + \gamma I)^{-1}(A + \gamma I)]^{r(0)}. \quad (61)$$

In particular, for $r = \alpha + \beta$, we have $r(\beta = 0) = \alpha$, $r(\alpha = 0) = \beta$, so that

$$\lim_{\beta \to 0} D_{\alpha+\beta}^{(\alpha,\beta)}[(A + \gamma I), (B + \gamma I)] \quad (62)$$
$$= \frac{1}{\alpha^2} \left\{ \text{tr}_X([(A + \gamma I)^{-1}(B + \gamma I)]^{\alpha} - I) - \alpha \log\det_X[(A + \gamma I)^{-1}(B + \gamma I)] \right\},$$

$$\lim_{\alpha \to 0} D_{\alpha+\beta}^{(\alpha,\beta)}[(A + \gamma I), (B + \gamma I)] \quad (63)$$
$$= \frac{1}{\beta^2} \left\{ \text{tr}_X([(B + \gamma I)^{-1}(A + \gamma I)]^{\beta} - I) - \beta \log\det_X[(B + \gamma I)^{-1}(A + \gamma I)] \right\}.$$

These are the direct generalizations of the corresponding formulas in the finite-dimensional setting. In fact, for $A, B \in \text{Sym}^{++}(n)$, $n \in \mathbb{N}$, by setting $\gamma = 0$, we obtain

$$\lim_{\beta \to 0} D_{\alpha+\beta}^{(\alpha,\beta)}[A, B] = \frac{1}{\alpha^2} \left\{ \text{tr}([(A^{-1}B)^{\alpha} - I) - \alpha \log\det(A^{-1}B) \right\}, \quad (64)$$

$$\lim_{\alpha \to 0} D_{\alpha+\beta}^{(\alpha,\beta)}[A, B] = \frac{1}{\beta^2} \left\{ \text{tr}([(B^{-1}A)^{\beta} - I) - \beta \log\det(B^{-1}A) \right\}. \quad (65)$$

These are precisely the finite-dimensional expressions given by Eqs. (5) and 6, which are Eqs. (23) and (22) in [15], respectively.

(ii) If $r(0) = r(\beta = 0) = 1$, we have for $\alpha > 0$ fixed,

$$\lim_{\beta \to 0} D_r^{(\alpha,\beta)}[(A + \gamma I), (B + \mu I)] = \frac{1}{\alpha^2} \left( \frac{\mu}{\gamma} - 1 \right) \log \frac{\mu}{\gamma}$$
$$+ \frac{1}{\alpha^2} \left\{ \text{tr}_X[(A + \gamma I)^{-1}(B + \mu I) - I] - \frac{\mu}{\gamma} \log\det_X[(A + \gamma I)^{-1}(B + \mu I)] \right\}$$
$$= \frac{1}{\alpha^2} d_{\text{logdet}}^{-1}[(A + \gamma I), (B + \mu I)]. \quad (66)$$

Similarly, if $r(0) = r(\alpha = 0) = 1$, we have for $\beta > 0$ fixed,

$$\lim_{\alpha \to 0} D_r^{(\alpha,\beta)}[(A+\gamma I), (B+\mu I)] = \frac{1}{\beta^2}\left(\frac{\gamma}{\mu} - 1\right)\log\frac{\gamma}{\mu}$$
$$+ \frac{1}{\beta^2}\left\{\operatorname{tr}_X[(B+\mu I)^{-1}(A+\gamma I) - I] - \frac{\gamma}{\mu}\log\det\nolimits_X[(B+\mu I)^{-1}(A+\gamma I)]\right\}$$
$$= \frac{1}{\beta^2} d_{\text{logdet}}^1[(A+\gamma I), (B+\mu I)]. \tag{67}$$

In particular, if $r \equiv 1$ as a constant function, then with $\beta = 1 - \alpha$, we have

$$\lim_{\alpha \to 1} D_1^{(\alpha,1-\alpha)}[(A+\gamma I), (B+\mu I)] = d_{\text{logdet}}^{-1}[(A+\gamma I), (B+\mu I)]$$
$$\lim_{\alpha \to 0} D_1^{(\alpha,1-\alpha)}[(A+\gamma I), (B+\mu I)] = d_{\text{logdet}}^1[(A+\gamma I), (B+\mu I)],$$

which are precisely the limiting cases stated in Eqs. (54) and (55) in Theorem 10.

## 6. Properties of the Alpha-Beta Log-Determinant divergences

The following results establish several important results of $D_r^{(\alpha,\beta)}$ as defined above, which generalize those from both the finite-dimensional setting [14, 15] and the infinite-dimensional Alpha Log-Det divergences [17].

**Theorem 13** (**Dual symmetry**).

$$D_r^{(\beta,\alpha)}[(B+\mu I), (A+\gamma I)] = D_r^{(\alpha,\beta)}[(A+\gamma I), (B+\mu I)]. \tag{68}$$

*In particular, for $\beta = \alpha$, we have*

$$D_r^{(\alpha,\alpha)}[(B+\mu I), (A+\gamma I)] = D_r^{(\alpha,\alpha)}[(A+\gamma I), (B+\mu I)]. \tag{69}$$

**Special case: Dual symmetry of the infinite-dimensional Alpha Log-Det divergences**. By Theorem 10, we have for $0 \leq \alpha \leq 1$,

$$D_1^{(\alpha,1-\alpha)}[(A+\gamma I), (B+\mu I)] = D_1^{(1-\alpha,\alpha)}[(B+\mu I), (A+\gamma I)]$$
$$\iff d_{\text{logdet}}^{1-2\alpha}[(A+\gamma I), (B+\mu I)] = d_{\text{logdet}}^{-(1-2\alpha)}[(B+\mu I), (A+\gamma I)]. \tag{70}$$

This is precisely the dual symmetry of the infinite-dimensional Alpha Log-Det divergences (Theorem 4 in [17]).

**Theorem 14** (**Dual invariance under inversion**).

$$D_r^{(\alpha,\beta)}[(A+\gamma I)^{-1},(B+\mu I)^{-1}] = D_{-r}^{(\alpha,\beta)}[(A+\gamma I),(B+\mu I)] \qquad (71)$$

**Special case: Dual invariance under inversion of the infinite-dimensional Alpha Log-Det divergences.** By Theorem 10, we have

$$D_1^{(\alpha,1-\alpha)}[(A+\gamma I)^{-1},(B+\mu I)^{-1}] = D_{-1}^{(\alpha,1-\alpha)}[(A+\gamma I),(B+\mu I)]$$
$$\iff d_{\text{logdet}}^{1-2\alpha}[(A+\gamma I)^{-1},(B+\mu I)^{-1}] = d_{\text{logdet}}^{-(1-2\alpha)}[(A+\gamma I),(B+\mu I)]. \qquad (72)$$

This is precisely the dual invariance under inversion of the infinite-dimensional Alpha Log-Det divergences (Theorem 5 in [17]).

**Theorem 15** (**Affine invariance**). *For any* $(A+\gamma I),(B+\mu I) \in \text{PTr}(\mathcal{H})$ *and any invertible* $(C+\nu I) \in \text{Tr}_X(\mathcal{H})$, $\nu \neq 0$,

$$D_r^{(\alpha,\beta)}[(C+\nu I)(A+\gamma I)(C+\nu I)^*,(C+\nu I)(B+\mu I)(C+\nu I)^*]$$
$$= D_r^{(\alpha,\beta)}[(A+\gamma I),(B+\mu I)]. \qquad (73)$$

**Theorem 16** (**Invariance under unitary transformations**). *For any* $(A+\gamma I),(B+\mu I) \in \text{PTr}(\mathcal{H})$ *and any* $C \in \mathcal{L}(\mathcal{H})$, *with* $CC^* = C^*C = I$,

$$D_r^{(\alpha,\beta)}[C(A+\gamma I)C^*, C(B+\mu I)C^*] = D_r^{(\alpha,\beta)}[(A+\gamma I),(B+\mu I)]. \qquad (74)$$

**Theorem 17.**

$$D_r^{(\alpha,\beta)}[(A+\gamma I),(B+\mu I)] = D_r^{(\alpha,\beta)}\left[\left(\Lambda+\frac{\gamma}{\mu}I\right), I\right]. \qquad (75)$$

**Theorem 18.** *Let* $\omega \in \mathbb{R}, \omega \neq 0$ *be arbitrary. Then*

$$D_{\omega r}^{(\omega\alpha,\omega\beta)}[(A+\gamma I),(B+\mu I)] = \frac{1}{\omega^2}D_r^{(\alpha,\beta)}\left[\left(\Lambda+\frac{\gamma}{\mu}I\right)^\omega, I\right]. \qquad (76)$$

The following two properties are important for proving that the square root function $\sqrt{D_{2\alpha}^{(\alpha,\alpha)}[(A+\gamma I),(B+\gamma I)]}$, is a metric on $\text{PTr}(\mathcal{H})(\gamma)$. We focus on the case $\alpha > 0$, since for $\alpha = 0$, $\sqrt{D_{2\alpha}^{(\alpha,\alpha)}[(A+\gamma I),(B+\mu I)]} = \frac{1}{\sqrt{2}}d_{\text{aiHS}}[(A+\gamma I),(B+\mu I)]$ is automatically a metric on $\text{PTr}(\mathcal{H})$.

**Theorem 19** (**Convergence in trace norm**). *Let $\alpha > 0$ be fixed. Let $\mathcal{H}$ be a separable Hilbert space. Let $A, B : \mathcal{H} \to \mathcal{H}$ be self-adjoint, trace class operators such that $(I + A) > 0$, $(I + B) > 0$. Let $\{A_n\}_{n \in \mathbb{N}}$, $\{B_n\}_{n \in \mathbb{N}}$ be sequences of self-adjoint, trace-class operators such that $\lim_{n \to \infty} ||A_n - A||_{\text{tr}} = 0$, $\lim_{n \to \infty} ||B_n - B||_{\text{tr}} = 0$. Then*

$$\lim_{n \to \infty} D_{2\alpha}^{(\alpha,\alpha)}[(I + A_n), (I + B_n)] = D_{2\alpha}^{(\alpha,\alpha)}[(I + A), (I + B)]. \quad (77)$$

**Theorem 20** (**Triangle inequality**). *Let $\alpha > 0$ be fixed. Let $\mathcal{H}$ be a separable Hilbert space. Let $\gamma > 0$, $\gamma \in \mathbb{R}$ be fixed. Let $A, B, C : \mathcal{H} \to \mathcal{H}$ be self-adjoint, trace class operators such that $(A + \gamma I) > 0$, $(B + \gamma I) > 0$, $(C + \gamma I) > 0$. Then*

$$\sqrt{D_{2\alpha}^{(\alpha,\alpha)}[(A + \gamma I), (B + \gamma I)]} \leq \sqrt{D_{2\alpha}^{(\alpha,\alpha)}[(A + \gamma I), (C + \gamma I)]}$$
$$+ \sqrt{D_{2\alpha}^{(\alpha,\alpha)}[(C + \gamma I), (B + \gamma I)]}. \quad (78)$$

In particular, for $\alpha = 1/2$ and $\gamma = 1$, we obtain the following triangle inequality.

**Theorem 21** (**Triangle inequality- square root of symmetric Stein divergence**). *Let $\mathcal{H}$ be a separable Hilbert space. Let $A, B, C : \mathcal{H} \to \mathcal{H}$ be self-adjoint trace-class operators with $A + I > 0$, $B + I > 0$, $C + I > 0$. Then*

$$\sqrt{\log \frac{\det(\frac{A+B}{2} + I)}{\sqrt{\det(A + I)\det(B + I)}}} \leq \sqrt{\log \frac{\det(\frac{A+C}{2} + I)}{\sqrt{\det(A + I)\det(C + I)}}}$$
$$+ \sqrt{\log \frac{\det(\frac{C+B}{2} + I)}{\sqrt{\det(C + I)\det(B + I)}}}. \quad (79)$$

**Theorem 22** (**Diagonalization**). *Let $\alpha \geq 0$ be fixed. Let $\mathcal{H}$ be a separable Hilbert space. Let $\gamma > 0$, $\gamma \in \mathbb{R}$, be fixed. Let $A, B : \mathcal{H} \to \mathcal{H}$ be self-adjoint trace class operators, such that $A + \gamma I > 0$, $B + \gamma I > 0$. Let $\text{Eig}(A), \text{Eig}(B) : \ell^2 \to \ell^2$ be diagonal operators with the diagonals consisting of the eigenvalues of $A$ and $B$, respectively, in decreasing order. Then*

$$D_{2\alpha}^{(\alpha,\alpha)}[(\text{Eig}(A) + \gamma I), (\text{Eig}(B) + \gamma I)] \leq D_{2\alpha}^{(\alpha,\alpha)}[(A + \gamma I), (B + \gamma I)]. \quad (80)$$

## 7. Alpha-Beta Log-Det divergences between RKHS covariance operators

Let $\mathcal{X}$ be an arbitrary non-empty set. We now compute the Alpha-Beta Log-Det divergences between covariance operators on an RKHS induced by a positive definite kernel $K$ on $\mathcal{X} \times \mathcal{X}$. In this case, we have explicit formulas for $D_r^{(\alpha,\beta)}$ via the corresponding Gram matrices. We recall that similar formulas exist in the cases of the Log-Hilbert-Schmidt distance [19], the infinite-dimensional affine-invariant Riemannian distance [18, 20], and the infinite-dimensional Alpha Log-Det divergences [17].

We first prove the following result.

**Theorem 23.** *Let $\mathcal{H}_1, \mathcal{H}_2$ be separable Hilbert spaces. Let $A, B : \mathcal{H}_1 \to \mathcal{H}_2$ be compact linear operators such that both $AA^* : \mathcal{H}_2 \to \mathcal{H}_2$ and $BB^* : \mathcal{H}_2 \to \mathcal{H}_2$ are trace class operators. Assume that $\dim(\mathcal{H}_2) = \infty$. Let $\alpha, \beta > 0$ be fixed. For any $r \in \mathbb{R}$, $r \neq 0$, for any $\gamma > 0$, $\mu > 0$,*

$$D_r^{(\alpha,\beta)}[(AA^* + \gamma I_{\mathcal{H}_2}), (BB^* + \mu I_{\mathcal{H}_2})] \tag{81}$$
$$= \frac{r(\delta - \frac{\alpha}{\alpha+\beta})}{\alpha\beta}\left(\log\frac{\gamma}{\mu}\right) + \frac{1}{\alpha\beta}\log\left(\frac{\alpha(\frac{\gamma}{\mu})^p + \beta(\frac{\gamma}{\mu})^{-q}}{\alpha + \beta}\right)$$
$$+ \frac{1}{\alpha\beta}\log\det\left[\frac{\alpha(\frac{\gamma}{\mu})^p(C + I_{\mathcal{H}_1} \otimes I_3)^p + \beta(\frac{\gamma}{\mu})^{-q}(C + I_{\mathcal{H}_1} \otimes I_3)^{-q}}{\alpha(\frac{\gamma}{\mu})^p + \beta(\frac{\gamma}{\mu})^{-q}}\right],$$

*where $\delta = \frac{\alpha\gamma^r}{\alpha\gamma^r + \beta\mu^r}$, $p = r(1-\delta)$, $q = r\delta$, and*

$$C = \begin{pmatrix} \frac{A^*A}{\gamma} & -\frac{A^*B}{\sqrt{\gamma\mu}}(I_{\mathcal{H}_1} + \frac{B^*B}{\mu})^{-1} & -\frac{A^*AA^*B}{\gamma\sqrt{\gamma\mu}}(I_{\mathcal{H}_1} + \frac{B^*B}{\mu})^{-1} \\ \frac{B^*A}{\sqrt{\gamma\mu}} & -\frac{B^*B}{\mu}(I_{\mathcal{H}_1} + \frac{B^*B}{\mu})^{-1} & -\frac{B^*AA^*B}{\gamma\mu}(I_{\mathcal{H}_1} + \frac{B^*B}{\mu})^{-1} \\ \frac{B^*A}{\sqrt{\gamma\mu}} & -\frac{B^*B}{\mu}(I_{\mathcal{H}_1} + \frac{B^*B}{\mu})^{-1} & -\frac{B^*AA^*B}{\gamma\mu}(I_{\mathcal{H}_1} + \frac{B^*B}{\mu})^{-1} \end{pmatrix}. \tag{82}$$

For comparison, the following is the corresponding version of $D_r^{(\alpha,\beta)}[(AA^* + \gamma I_{\mathcal{H}_2}), (BB^* + \mu I_{\mathcal{H}_2})]$, using the finite-dimensional formula given in Eq. (19), when $\dim(\mathcal{H}_2) < \infty$.

**Theorem 24.** *Let $\mathcal{H}_1, \mathcal{H}_2$ be separable Hilbert spaces. Let $A, B : \mathcal{H}_1 \to \mathcal{H}_2$ be compact linear operators such that both $AA^* : \mathcal{H}_2 \to \mathcal{H}_2$ and $BB^* : \mathcal{H}_2 \to \mathcal{H}_2$ are trace class operators. Assume that $\dim(\mathcal{H}_2) < \infty$. Let $\alpha, \beta > 0$ be fixed. For any*

$r \in \mathbb{R}$, $r \neq 0$, *for any $\gamma > 0$, $\mu > 0$,*

$$D_r^{(\alpha,\beta)}[(AA^* + \gamma I_{\mathcal{H}_2}), (BB^* + \mu I_{\mathcal{H}_2})] \tag{83}$$
$$= \frac{1}{\alpha\beta}\left[\log\left(\frac{\alpha(\frac{\gamma}{\mu})^p + \beta(\frac{\gamma}{\mu})^{-q}}{\alpha+\beta}\right)\right]\dim(\mathcal{H}_2)$$
$$+ \frac{1}{\alpha\beta}\log\det\left[\frac{\alpha(\frac{\gamma}{\mu})^p(C+I_{\mathcal{H}_1}\otimes I_3)^p + \beta(\frac{\gamma}{\mu})^{-q}(C+I_{\mathcal{H}_1}\otimes I_3)^{-q}}{\alpha(\frac{\gamma}{\mu})^p + \beta(\frac{\gamma}{\mu})^{-q}}\right],$$

*where $p = r\frac{\beta}{\alpha+\beta}$, $q = r\frac{\alpha}{\alpha+\beta}$, and $C$ is as given in Theorem 23.*

Let us briefly recall the RKHS covariance operators discussed in [17]. Let $\mathbf{x} = [x_1, \ldots, x_m]$ be a data matrix randomly sampled from $\mathcal{X}$ according to a Borel probability distribution $\rho$, where $m \in \mathbb{N}$ is the number of observations. Let $K$ be a positive definite kernel on $\mathcal{X} \times \mathcal{X}$ and $\mathcal{H}_K$ its induced reproducing kernel Hilbert space (RKHS). Let $\Phi : \mathcal{X} \to \mathcal{H}_K$ be the corresponding feature map, so that $K(x, y) = \langle \Phi(x), \Phi(y) \rangle_{\mathcal{H}_K}$ for all pairs $(x, y) \in \mathcal{X} \times \mathcal{X}$. The feature map $\Phi$ gives rise to the bounded linear operator

$$\Phi(\mathbf{x}) : \mathbb{R}^m \to \mathcal{H}_K, \quad \Phi(\mathbf{x})\mathbf{b} = \sum_{j=1}^{m} b_j \Phi(x_j), \quad \mathbf{b} \in \mathbb{R}^m. \tag{84}$$

The operator $\Phi(\mathbf{x})$ can also be viewed as the (potentially infinite) mapped data matrix $\Phi(\mathbf{x}) = [\Phi(x_1), \ldots, \Phi(x_m)]$ of size $\dim(\mathcal{H}_K) \times m$ in the feature space $\mathcal{H}_K$, with the $j$th column being $\Phi(x_j)$. The corresponding empirical covariance operator for $\Phi(\mathbf{x})$ is defined to be

$$C_{\Phi(\mathbf{x})} = \frac{1}{m}\Phi(\mathbf{x})J_m\Phi(\mathbf{x})^* : \mathcal{H}_K \to \mathcal{H}_K, \tag{85}$$

where $\Phi(\mathbf{x})^* : \mathcal{H}_K \to \mathbb{R}^m$ is the adjoint operator of $\Phi(\mathbf{x})$ and $J_m$ is the centering matrix, defined by $J_m = I_m - \frac{1}{m}\mathbf{1}_m\mathbf{1}_m^T$ with $\mathbf{1}_m = (1, \ldots, 1)^T \in \mathbb{R}^m$.

Let $\mathbf{x} = [x_i]_{i=1}^m$, $\mathbf{y} = [y_i]_{i=1}^m$, $m \in \mathbb{N}$, be two random data matrices sampled from $\mathcal{X}$ according to two Borel probability distributions and $C_{\Phi(\mathbf{x})}$, $C_{\Phi(\mathbf{y})}$ be the corresponding covariance operators induced by the kernel $K$. Let $K[\mathbf{x}]$, $K[\mathbf{y}]$, and $K[\mathbf{x}, \mathbf{y}]$ be the $m \times m$ Gram matrices defined by

$$(K[\mathbf{x}])_{ij} = K(x_i, x_j), \quad (K[\mathbf{y}])_{ij} = K(y_i, y_j),$$
$$(K[\mathbf{x}, \mathbf{y}])_{ij} = K(x_i, y_j), \quad 1 \leq i, j \leq m. \tag{86}$$

Let $A = \frac{1}{\sqrt{m}}\Phi(\mathbf{x})J_m : \mathbb{R}^m \to \mathcal{H}_K$, $B = \frac{1}{\sqrt{m}}\Phi(\mathbf{y})J_m : \mathbb{R}^m \to \mathcal{H}_K$, so that

$$AA^* = C_{\Phi(\mathbf{x})}, \;\; BB^* = C_{\Phi(\mathbf{y})}, \;\; A^*A = \frac{1}{m}J_m K[\mathbf{x}]J_m, \;\; B^*B = \frac{1}{m}J_m K[\mathbf{y}]J_m,$$

$$A^*B = \frac{1}{m}J_m K[\mathbf{x},\mathbf{y}]J_m, \;\; B^*A = \frac{1}{m}J_m K[\mathbf{y},\mathbf{x}]J_m. \tag{87}$$

Theorems 23 and 24 can then be applied to give closed form formulas for the divergences between $(C_{\Phi(\mathbf{x})} + \gamma I)$ and $(C_{\Phi(\mathbf{y})} + \mu I)$, as follows.

**Theorem 25** (**Alpha-Beta Log-Det divergences between RKHS covariance operators - Infinite-dimensional version**). *Let $\alpha, \beta > 0$ be fixed. Let $r \in \mathbb{R}$, $r \neq 0$ be fixed. Assume that $\dim(\mathcal{H}_K) = \infty$. For any $\gamma > 0$, $\mu > 0$, the divergence $D_r^{(\alpha,\beta)}[(C_{\Phi(\mathbf{x})} + \gamma I), (C_{\Phi(\mathbf{y})} + \mu I)]$ is given by*

$$D_r^{(\alpha,\beta)}[(C_{\Phi(\mathbf{x})} + \gamma I), (C_{\Phi(\mathbf{y})} + \mu I)] \tag{88}$$
$$= \frac{r(\delta - \frac{\alpha}{\alpha+\beta})}{\alpha\beta}\left(\log\frac{\gamma}{\mu}\right) + \frac{1}{\alpha\beta}\log\left(\frac{\alpha(\frac{\gamma}{\mu})^p + \beta(\frac{\gamma}{\mu})^{-q}}{\alpha+\beta}\right)$$
$$+ \frac{1}{\alpha\beta}\log\det\left[\frac{\alpha(\frac{\gamma}{\mu})^p(C+I_{3m})^p + \beta(\frac{\gamma}{\mu})^{-q}(C+I_{3m})^{-q}}{\alpha(\frac{\gamma}{\mu})^p + \beta(\frac{\gamma}{\mu})^{-q}}\right],$$

*where $\delta = \frac{\alpha\gamma^r}{\alpha\gamma^r + \beta\mu^r}$, $p = r(1-\delta)$, $q = r\delta$, and*

$$C = \begin{pmatrix} C_{11} & C_{12} & C_{13} \\ C_{21} & C_{22} & C_{23} \\ C_{21} & C_{22} & C_{23} \end{pmatrix} \in \mathbb{R}^{3m \times 3m}. \tag{89}$$

*Here the sub-matrices $C_{ij}$, $i = 1, 2$, $j = 1, 2, 3$, each of size $m \times m$, are given by*

$$C_{11} = \frac{1}{\gamma m}J_m K[\mathbf{x}]J_m, \tag{90}$$

$$C_{12} = -\frac{1}{\sqrt{\gamma\mu}m}J_m K[\mathbf{x},\mathbf{y}]J_m\left(I_m + \frac{1}{\mu m}J_m K[\mathbf{y}]J_m\right)^{-1}, \tag{91}$$

$$C_{13} = -\frac{1}{\gamma\sqrt{\gamma\mu}m^2}J_m K[\mathbf{x}]J_m K[\mathbf{x},\mathbf{y}]J_m\left(I_m + \frac{1}{\mu m}J_m K[\mathbf{y}]J_m\right)^{-1}, \tag{92}$$

$$C_{21} = \frac{1}{\sqrt{\gamma\mu}m}J_m K[\mathbf{y},\mathbf{x}]J_m, \tag{93}$$

$$C_{22} = -\frac{1}{\mu m}J_m K[\mathbf{y}]J_m\left(I_m + \frac{1}{\mu m}J_m K[\mathbf{y}]J_m\right)^{-1}, \tag{94}$$

$$C_{23} = -\frac{1}{\gamma\mu m^2}J_m K[\mathbf{y},\mathbf{x}]J_m K[\mathbf{x},\mathbf{y}]J_m\left(I_m + \frac{1}{\mu m}J_m K[\mathbf{y}]J_m\right)^{-1}. \tag{95}$$

**Theorem 26** (**Alpha-Beta Log-Det divergences between RKHS covariance operators - Finite-dimensional version**). *Let $\alpha, \beta > 0$ be fixed. Let $r \in \mathbb{R}$, $r \neq 0$ be fixed. Assume that $\dim(\mathcal{H}_K) < \infty$. For any $\gamma > 0$, $\mu > 0$, the divergence $D_r^{(\alpha,\beta)}[(C_{\Phi(\mathbf{x})} + \gamma I), (C_{\Phi(\mathbf{y})} + \mu I)]$ is given by*

$$D_r^{(\alpha,\beta)}[(C_{\Phi(\mathbf{x})} + \gamma I), (C_{\Phi(\mathbf{y})} + \mu I)] \tag{96}$$
$$= \frac{1}{\alpha\beta}\left[\log\left(\frac{\alpha(\frac{\gamma}{\mu})^p + \beta(\frac{\gamma}{\mu})^{-q}}{\alpha + \beta}\right)\right]\dim(\mathcal{H}_K)$$
$$+ \frac{1}{\alpha\beta}\log\det\left[\frac{\alpha(\frac{\gamma}{\mu})^p(C + I_{3m})^p + \beta(\frac{\gamma}{\mu})^{-q}(C + I_{3m})^{-q}}{\alpha(\frac{\gamma}{\mu})^p + \beta(\frac{\gamma}{\mu})^{-q}}\right],$$

*where $p = r\frac{\beta}{\alpha+\beta}$, $q = r\frac{\alpha}{\alpha+\beta}$, and $C$ is as given in Theorem 25.*

*Remark* 10. The closed form formulas for $D_r^{(\alpha,\beta)}[(C_{\Phi(\mathbf{x})} + \gamma I), (C_{\Phi(\mathbf{y})} + \mu I)]$ given in Eqs. (88) and (96) in Theorems 25 and 26, respectively, coincide if and only if $\gamma = \mu$. If $\gamma \neq \mu$, then the right hand side of Eq. (96) approaches infinity when $\dim(\mathcal{H}_K) \to \infty$. Thus in general, the infinite-dimensional version is *not* obtainable as the limit of the finite-dimensional version as the dimension goes to infinity.

*Remark* 11. The closed form formulas given by Eqs. (88) and (96) in Theorems 25 and 26, respectively, are derived under more general conditions than those in [17] and are consequently more general but more complicated than the corresponding closed form formulas for the Alpha Log-Det divergences in [17] (see Theorems 12,13, 15, 16 in [17]). Thus for practical applications involving the Alpha Log-Det divergences, the corresponding closed form formulas in [17] should be employed.

### Appendix A. Proofs of main results

*Appendix A.1. Proofs for the general Alpha-Beta Log-Determinant divergences*

In this section, we prove Lemma 1, Proposition 1, and Theorems 5, 6, 7, and 8.

***Proof of Lemma 1.*** Since any bounded operator $A$ commutes with the identity operator $I$, we have

$$\exp(A + \gamma I) = e^\gamma \exp(A) = e^\gamma\left(I + \sum_{j=1}^\infty \frac{A^j}{j!}\right) = e^\gamma I + e^\gamma \sum_{j=1}^\infty \frac{A^j}{j!},$$

where $\sum_{j=1}^{\infty} \frac{A^j}{j!}$ is trace class, since

$$\left\|\sum_{j=1}^{\infty} \frac{A^j}{j!}\right\|_{\text{tr}} \leq \sum_{j=1}^{\infty} \frac{||A||_{\text{tr}}^j}{j!} = \exp(||A||_{\text{tr}}) - 1 < \infty.$$

Thus $\exp(A + \gamma I) \in \text{Tr}_X(\mathcal{H})$. This completes the proof. $\square$

***Proof of Proposition 1***. For $(A+\gamma I), (B+\mu I) \in \text{PTr}(\mathcal{H})$, we have $(B+\mu I)^{-1/2}(A+\gamma I)(B+\mu I)^{-1/2} \in \text{PTr}(\mathcal{H})$ and the logarithm $\log[(B+\mu I)^{-1/2}(A+\gamma I)(B+\mu I)^{-1/2}] \in \text{Tr}_X(\mathcal{H})$ is well-defined. By the discussion preceding Proposition 1, we have

$$\begin{aligned}
&\log[(A + \gamma I)(B + \mu I)^{-1}] \\
&= \log[(B + \mu I)^{1/2}(B + \mu I)^{-1/2}(A + \gamma I)(B + \mu I)^{-1/2}(B + \mu I)^{-1/2}] \\
&= (B + \mu I)^{1/2} \log[(B + \mu I)^{-1/2}(A + \gamma I)(B + \mu I)^{-1/2}](B + \mu I)^{-1/2} \\
&= (B + \mu I)^{1/2} \log\left(\Lambda + \frac{\gamma}{\mu}I\right)(B + \mu I)^{-1/2} \in \text{Tr}_X(\mathcal{H}).
\end{aligned}$$

For the power function, we have

$$\begin{aligned}
[(A + \gamma I)(B + \mu I)^{-1}]^\alpha &= \exp(\alpha \log[(A + \gamma I)(B + \mu I)^{-1}]) \\
&= \exp\left[(B + \mu I)^{1/2}\alpha \log\left(\Lambda + \frac{\gamma}{\mu}I\right)(B + \mu I)^{-1/2}\right] \\
&= (B + \mu I)^{1/2} \exp\left[\alpha \log\left(\Lambda + \frac{\gamma}{\mu}I\right)\right](B + \mu I)^{-1/2} \\
&= (B + \mu I)^{1/2} \left(\Lambda + \frac{\gamma}{\mu}I\right)^\alpha (B + \mu I)^{-1/2}.
\end{aligned}$$

For the sum of two power functions, we then have

$$\begin{aligned}
&\frac{\alpha[(A + \gamma I)(B + \mu I)^{-1}]^p + \beta[(A + \gamma I)(B + \mu I)^{-1}]^q}{\alpha + \beta} \\
&= (B + \mu I)^{1/2} \left[\frac{\alpha(\Lambda + \frac{\gamma}{\mu}I)^p + \beta(\Lambda + \frac{\gamma}{\mu}I)^q}{\alpha + \beta}\right](B + \mu I)^{-1/2}.
\end{aligned}$$

By Lemma 5 in [17], $\det_X[C(A + \gamma I)C^{-1}] = \det_X(A + \gamma I)$ for any invertible

operator $C \in \mathcal{L}(\mathcal{H})$. It follows that

$$\text{det}_X \left[ \frac{\alpha[(A+\gamma I)(B+\mu I)^{-1}]^p + \beta[(A+\gamma I)(B+\mu I)^{-1}]^q}{\alpha + \beta} \right]$$
$$=\text{det}_X \left[ \frac{\alpha(\Lambda + \frac{\gamma}{\mu}I)^p + \beta(\Lambda + \frac{\gamma}{\mu}I)^q}{\alpha + \beta} \right].$$

This completes the proof. $\square$

*Proof of Theorem 5.* By definition of the power function, we have

$$\alpha(A+\gamma I)^p + (1-\alpha)(B+\mu I)^q = \alpha \exp[p \log(A+\gamma I)] + (1-\alpha)\exp[q \log(B+\mu I)]$$
$$= \alpha \exp\left[p \log\left(\frac{A}{\gamma}+I\right) + p(\log \gamma)I\right] + (1-\alpha)\exp\left[q\log\left(\frac{B}{\mu}+I\right) + q(\log \mu)I\right]$$
$$= \alpha \gamma^p \left(\frac{A}{\gamma}+I\right)^p + (1-\alpha)\mu^q\left(\frac{B}{\mu}+I\right)^q.$$

It follows that for $\delta = \frac{\alpha\gamma^p}{\alpha\gamma^p + (1-\alpha)\mu^q}$, $1-\delta = \frac{(1-\alpha)\mu^q}{\alpha\gamma^p + (1-\alpha)\mu^q}$, we have

$$\text{det}_X[\alpha(A+\gamma I)^p + (1-\alpha)(B+\mu I)^q]$$
$$= [\alpha\gamma^p + (1-\alpha)\mu^q]\det\left[\frac{\alpha\gamma^p}{\alpha\gamma^p + (1-\alpha)\mu^q}\left(\frac{A}{\gamma}+I\right)^p + \frac{(1-\alpha)\mu^q}{\alpha\gamma^p + (1-\alpha)\mu^q}\left(\frac{B}{\mu}+I\right)^q\right]$$
$$\geq [\alpha\gamma^p + (1-\alpha)\mu^q]\det\left(\frac{A}{\gamma}+I\right)^{p\delta}\det\left(\frac{B}{\mu}+I\right)^{q(1-\delta)}$$

by Proposition 7 in [17]

$$\geq \gamma^{p\alpha}\mu^{(1-\alpha)q}\det\left(\frac{A}{\gamma}+I\right)^{p\delta}\det\left(\frac{B}{\mu}+I\right)^{q(1-\delta)}$$

by Ky Fan's Inequality applied to $\alpha\gamma^p + (1-\alpha)\mu^q$

$$= \gamma^{p(\alpha-\delta)}\mu^{-q(\alpha-\delta)}\text{det}_X(A+\gamma I)^{p\delta}\text{det}_X(B+\mu I)^{q(1-\delta)}$$
$$= \left(\frac{\gamma^p}{\mu^q}\right)^{\alpha-\delta}\text{det}_X(A+\gamma I)^{p\delta}\text{det}_X(B+\mu I)^{q(1-\delta)}.$$

For $0 < \alpha < 1$, equality happens if and only if simultaneously, we have

$$\left(\frac{A}{\gamma}+I\right)^p = \left(\frac{B}{\mu}+I\right)^q \text{ and } \gamma^p = \mu^q \iff (A+\gamma I)^p = (B+\mu I)^q.$$

In particular, for $\gamma = \mu$, the condition $\gamma^p = \mu^q$ becomes

$$\gamma^p = \gamma^q \iff \gamma^{p-q} = 1 \iff p = q \text{ if } \gamma \neq 1.$$

With the conditions $\gamma = \mu \neq 1$ and $p = q$, we then have

$$\left(\frac{A}{\gamma} + I\right)^p = \left(\frac{B}{\gamma} + I\right)^p \iff A = B.$$

This completes the proof of the theorem. $\square$

***Proof of Theorem 6.*** Recall that we write the operator $(B + \mu I)^{-1/2}(A + \gamma I)(B + \mu I)^{-1/2}$ in the form

$$(B + \mu I)^{-1/2}(A + \gamma I)(B + \mu I)^{-1/2} = \Lambda + (\gamma/\mu)I \in \mathrm{PTr}(\mathcal{H}).$$

Its inverse has the form

$$(B + \mu I)^{1/2}(A + \gamma I)^{-1}(B + \mu I)^{1/2} = [\Lambda + (\gamma/\mu)I]^{-1}$$

$$= \frac{\mu}{\gamma}I - \left(\frac{\mu}{\gamma}\right)^2 \Lambda \left(I + \frac{\mu}{\gamma}\Lambda\right)^{-1} \in \mathrm{PTr}(\mathcal{H}).$$

It follows from Corollary 1 that

$$\det\nolimits_{\mathrm{X}} \left[\frac{\alpha[(\Lambda + (\gamma/\mu)I]^p + \beta[(\Lambda + (\gamma/\mu)I)^{-1}]^q}{\alpha + \beta}\right]$$

$$\geq \left(\frac{(\gamma/\mu)^p}{(\mu/\gamma)^q}\right)^{\frac{\alpha}{\alpha+\beta}-\delta} \det\nolimits_{\mathrm{X}}(\Lambda + (\gamma/\mu)I)^{p\delta} \det\nolimits_{\mathrm{X}}[(\Lambda + (\gamma/\mu)I)]^{-q(1-\delta)}$$

$$= \left(\frac{\gamma}{\mu}\right)^{(p+q)(\frac{\alpha}{\alpha+\beta}-\delta)} \det\nolimits_{\mathrm{X}}(\Lambda + (\gamma/\mu)I)^{p\delta} \det\nolimits_{\mathrm{X}}[(\Lambda + (\gamma/\mu)I)]^{-q(1-\delta)}, \quad (\text{A.1})$$

where $\delta = \frac{\alpha(\frac{\gamma}{\mu})^p}{\alpha(\frac{\gamma}{\mu})^p + \beta(\frac{\mu}{\gamma})^q} = \frac{\alpha(\frac{\gamma}{\mu})^{p+q}}{\alpha(\frac{\gamma}{\mu})^{p+q}+\beta}$, $1 - \delta = \frac{\beta(\frac{\mu}{\gamma})^q}{\alpha(\frac{\gamma}{\mu})^p + \beta(\frac{\mu}{\gamma})^q} = \frac{\beta}{\alpha(\frac{\gamma}{\mu})^{p+q}+\beta}$.

For the two determinants on the right hand side of (A.1) to cancel each other out, we need

$$p\delta = q(1-\delta) \iff \alpha p \left(\frac{\gamma}{\mu}\right)^p = \beta q \left(\frac{\mu}{\gamma}\right)^q \iff \alpha p \left(\frac{\gamma}{\mu}\right)^{p+q} = \beta q.$$

Assuming that this condition holds, then along with the definition of $D_{(p,q)}^{(\alpha,\beta)}$, (A.1) gives

$$\left[\left(\frac{\gamma}{\mu}\right)^{(p+q)(\delta-\frac{\alpha}{\alpha+\beta})} \det\nolimits_{\mathrm{X}}\left(\frac{\alpha(\Lambda + \frac{\gamma}{\mu}I)^p + \beta(\Lambda + \frac{\gamma}{\mu}I)^{-q}}{\alpha + \beta}\right)\right] \geq 1$$

$$\iff D_{(p,q)}^{(\alpha,\beta)}[(A + \gamma I), (B + \mu I)] \geq 0.$$

In the inequality in (A.1), the equality sign happens if and only if

$$[(\Lambda + (\gamma/\mu)I]^p = [(\Lambda + (\gamma/\mu)I]^{-q} \iff [(\Lambda + (\gamma/\mu)I]^{p+q} = I.$$

If $p + q = 0$, then this is always true, so that $D_{(p,q)}^{(\alpha,\beta)}[(A + \gamma I), (B + \mu I)] = 0$ for all pairs $(A + \gamma I), (B + \mu I) \in \mathrm{PTr}(\mathcal{H})$, which is not what we want. In fact, with $p + q = 0$, the condition $\alpha p \left(\frac{\gamma}{\mu}\right)^{p+q} = \beta q$ gives $(\alpha + \beta)p = 0 \Rightarrow p = 0 \Rightarrow q = 0$.

If $p + q \neq 0$, since $\Lambda + (\gamma/\mu)I > 0$, this happens if and only if

$$\Lambda + (\gamma/\mu)I = I \iff (B + \mu I)^{-1/2}(A + \gamma I)(B + \mu I)^{-1/2} = I$$
$$\iff A + \gamma I = B + \mu I \iff A = B \text{ and } \gamma = \mu.$$

This completes the proof. □

*Proof of Theorem 7.* Under the condition $p + q = r$, by Theorem 6, we have

$$\alpha p \left(\frac{\gamma}{\mu}\right)^r = \beta(r - p) \Rightarrow p = \frac{\beta r}{\alpha \left(\frac{\gamma}{\mu}\right)^r + \beta}$$

It follows then that $q = r - p = \frac{r\alpha\left(\frac{\gamma}{\mu}\right)^r}{\alpha\left(\frac{\gamma}{\mu}\right)^r + \beta}$. The equivalence of Eqs. (8) and (9) follows from Proposition 1. □

*Proof of Theorem 8.* We have

$$\frac{\alpha(\Lambda + \frac{\gamma}{\mu}I)^p + \beta(\Lambda + \frac{\gamma}{\mu}I)^{-q}}{\alpha + \beta} = \frac{\alpha(\frac{\gamma}{\mu})^p(\frac{\mu}{\gamma}\Lambda + I)^p + \beta(\frac{\gamma}{\mu})^{-q}(\frac{\mu}{\gamma}\Lambda + I)^{-q}}{\alpha + \beta}$$
$$= \frac{\alpha(\frac{\gamma}{\mu})^p(I + C_1) + \beta(\frac{\gamma}{\mu})^{-q}(I + C_2)}{\alpha + \beta}$$
$$= \frac{\left[\alpha(\frac{\gamma}{\mu})^p + \beta(\frac{\gamma}{\mu})^{-q}\right]I + \left[\alpha(\frac{\gamma}{\mu})^p C_1 + \beta(\frac{\gamma}{\mu})^{-q}C_2\right]}{\alpha + \beta}$$
$$= \frac{\alpha(\frac{\gamma}{\mu})^p + \beta(\frac{\gamma}{\mu})^{-q}}{\alpha + \beta}\left[I + \frac{\alpha(\frac{\gamma}{\mu})^p C_1 + \beta(\frac{\gamma}{\mu})^{-q}C_2}{\alpha(\frac{\gamma}{\mu})^p + \beta(\frac{\gamma}{\mu})^{-q}}\right],$$

where $C_1 = \sum_{k=1}^{\infty} \frac{p^k}{k!}\left[\log\left(\frac{\mu}{\gamma}\Lambda + I\right)\right]^k \in \mathrm{Tr}(\mathcal{H})$, $C_2 = \sum_{k=1}^{\infty} \frac{(-1)^k q^k}{k!}\left[\log\left(\frac{\mu}{\gamma}\Lambda + I\right)\right]^k \in$

Tr($\mathcal{H}$). By definition of the $\det_X$ function, we then have

$$\log \det_X \left[ \frac{\alpha(\Lambda + \frac{\gamma}{\mu}I)^p + \beta(\Lambda + \frac{\gamma}{\mu}I)^{-q}}{\alpha + \beta} \right]$$

$$= \log \left( \frac{\alpha(\frac{\gamma}{\mu})^p + \beta(\frac{\gamma}{\mu})^{-q}}{\alpha + \beta} \right) + \log \det \left[ I + \frac{\alpha(\frac{\gamma}{\mu})^p C_1 + \beta(\frac{\gamma}{\mu})^{-q} C_2}{\alpha(\frac{\gamma}{\mu})^p + \beta(\frac{\gamma}{\mu})^{-q}} \right]$$

$$= \log \left( \frac{\alpha(\frac{\gamma}{\mu})^p + \beta(\frac{\gamma}{\mu})^{-q}}{\alpha + \beta} \right) + \log \det \left[ \frac{\alpha(\Lambda + \frac{\gamma}{\mu}I)^p + \beta(\Lambda + \frac{\gamma}{\mu}I)^{-q}}{\alpha(\frac{\gamma}{\mu})^p + \beta(\frac{\gamma}{\mu})^{-q}} \right].$$

This, together with the definition of $D_{(p,q)}^{(\alpha,\beta)}$, gives us the desired expression. $\square$

*Appendix A.2. Proofs for the Affine-invariant Riemannian distance*

In this section, we prove Theorem 9. We first need the following preliminary results.

**Lemma 2.** *Let $\gamma > 0$. Assume that $r = r(\alpha)$ is smooth, with $r(0) = 0$. Let $\delta = \frac{\gamma^r}{\gamma^r + 1}$. Then*

$$\lim_{\alpha \to 0} \frac{r(\delta - \frac{1}{2})}{\alpha^2} = \frac{[r'(0)]^2}{4} \log \gamma. \qquad (A.2)$$

*In particular, for $r = 2\alpha$, we have*

$$\lim_{\alpha \to 0} \frac{r(\delta - \frac{1}{2})}{\alpha^2} = \log \gamma. \qquad (A.3)$$

*Proof of Lemma 2.* By L'Hopital's rule applied twice, we obtain

$$\lim_{\alpha \to 0} \frac{r(\delta - \frac{1}{2})}{\alpha^2} = \lim_{\alpha \to 0} \frac{r(\gamma^r - 1)}{2\alpha^2(\gamma^r + 1)} = \lim_{\alpha \to 0} \frac{r(\gamma^r - 1)}{4\alpha^2}$$

$$= \lim_{\alpha \to 0} \frac{r'(\alpha)(\gamma^r - 1) + r\gamma^r r'(\alpha) \log \gamma}{8\alpha}$$

$$= \lim_{\alpha \to 0} \frac{r''(\alpha)(\gamma^r - 1) + \gamma^r (r'(\alpha))^2 \log \gamma + \gamma^r (r'(\alpha))^2 \log \gamma}{8}$$

$$+ \lim_{\alpha \to 0} \frac{r\gamma^r (r'(\alpha) \log \gamma)^2 + r\gamma^r r''(\alpha) \log \gamma}{8}$$

$$= \frac{[r'(0)]^2 \log \gamma}{4}.$$

This completes the proof. $\square$

**Lemma 3.** *Let $\gamma > 0$ be fixed. Let $\lambda > 0$ be fixed. Assume that $r = r(\alpha)$ is smooth, with $r(0) = 0$. Define $\delta = \frac{\gamma^r}{\gamma^r + 1}$, $p = r(1 - \delta)$, $q = r\delta$. Then*

$$\lim_{\alpha \to 0} \frac{1}{\alpha^2} \log\left(\frac{\lambda^p + \lambda^{-q}}{2}\right) = \frac{[r'(0)]^2}{4}\left[-(\log \gamma)(\log \lambda) + \frac{1}{2}(\log \lambda)^2\right]. \quad (A.4)$$

*In particular, if $\gamma = \lambda$, then*

$$\lim_{\alpha \to 0} \frac{1}{\alpha^2} \log\left(\frac{\gamma^p + \gamma^{-q}}{2}\right) = -\frac{[r'(0)]^2}{8}(\log \gamma)^2. \quad (A.5)$$

*Proof of Lemma 3.* For $p$, $q$ sufficiently small,

$$\lambda^p = e^{p \log \lambda} = 1 + p \log \lambda + \frac{p^2}{2}(\log \lambda)^2 + o(p^3),$$

$$\lambda^{-q} = e^{-q \log \lambda} = 1 - q \log \lambda + \frac{q^2}{2}(\log \lambda)^2 + o(q^3).$$

Thus for $\alpha$ sufficiently small, so that $p = o(\alpha)$, $q = o(\alpha)$, we have

$$\frac{\lambda^p + \lambda^{-q}}{2} = 1 + \frac{p - q}{2} \log \lambda + \frac{p^2 + q^2}{4}(\log \lambda)^2 + o(p^3, q^3)$$

$$= 1 + r\left(\frac{1}{2} - \delta\right)(\log \lambda) + \frac{r^2}{4}\left[(1 - \delta)^2 + \delta^2\right](\log \lambda)^2 + o(\alpha^3).$$

By Lemma 2, we have

$$\lim_{\alpha \to 0} \frac{r\left(\frac{1}{2} - \delta\right)}{\alpha^2} = -\frac{[r'(0)]^2}{4} \log \gamma.$$

We have by L'Hopital's rule

$$\lim_{\alpha \to 0} \frac{r^2}{\alpha^2} = \lim_{\alpha \to 0} \frac{2rr'(\alpha)}{2\alpha} = \lim_{\alpha \to 0} [r'(\alpha)]^2 + rr''(\alpha) = [r'(0)]^2.$$

Since $\lim_{\alpha \to 0} \delta = \frac{1}{2}$, it follows then that

$$\lim_{\alpha \to 0} \frac{r^2}{4\alpha^2}[(1 - \delta)^2 + \delta^2] = \frac{[r'(0)]^2}{8}.$$

Combining these limits with $\lim_{x \to 0} \frac{\log(1+ax)}{x} = a$, we obtain

$$\lim_{\alpha \to 0} \frac{1}{\alpha^2} \log\left(\frac{\lambda^p + \lambda^{-q}}{2}\right) = \frac{[r'(0)]^2}{4}\left[-(\log \gamma)(\log \lambda) + \frac{1}{2}(\log \lambda)^2\right].$$

This completes the proof of the lemma. □

**Lemma 4.** *Let $\gamma > 0$ be fixed. Let $\lambda \in \mathbb{R}$ be fixed such that $\lambda + \gamma > 0$. Assume that $r = r(\alpha)$ is smooth, with $r(0) = 0$. Define $\delta = \frac{\gamma^r}{\gamma^r + 1}$, $p = r(1 - \delta)$, $q = r\delta$. Then*

$$\lim_{\alpha \to 0} \frac{1}{\alpha^2} \log\left(\frac{(\lambda + \gamma)^p + (\lambda + \gamma)^{-q}}{\gamma^p + \gamma^{-q}}\right) = \frac{[r'(0)]^2}{8} [\log(\lambda + \gamma) - \log \gamma]^2$$

$$= \frac{[r'(0)]^2}{8} \left[\log\left(\frac{\lambda}{\gamma} + 1\right)\right]^2. \quad (A.6)$$

*In particular, for $r = r(\alpha) = 2\alpha$, we have*

$$\lim_{\alpha \to 0} \frac{1}{\alpha^2} \log\left(\frac{(\lambda + \gamma)^p + (\lambda + \gamma)^{-q}}{\gamma^p + \gamma^{-q}}\right) = \frac{1}{2} [\log(\lambda + \gamma) - \log \gamma]^2$$

$$= \frac{1}{2} \left[\log\left(\frac{\lambda}{\gamma} + 1\right)\right]^2. \quad (A.7)$$

*Proof of Lemma 4.* We have by Lemma 3

$$\lim_{\alpha \to 0} \frac{1}{\alpha^2} \log\left(\frac{(\lambda + \gamma)^p + (\lambda + \gamma)^{-q}}{\gamma^p + \gamma^{-q}}\right)$$

$$= \lim_{\alpha \to 0} \frac{1}{\alpha^2} \log\left(\frac{(\lambda + \gamma)^p + (\lambda + \gamma)^{-q}}{2}\right) - \lim_{\alpha \to 0} \frac{1}{\alpha^2} \log\left(\frac{\gamma^p + \gamma^{-q}}{2}\right)$$

$$= \frac{[r'(0)]^2}{4} \left(-(\log \gamma)[\log(\lambda + \gamma)] + \frac{1}{2}[\log(\lambda + \gamma)]^2 - [-\frac{1}{2}(\log \gamma)^2]\right)$$

$$= \frac{[r'(0)]^2}{8} [\log(\lambda + \gamma) - \log \gamma]^2$$

$$= \frac{[r'(0)]^2}{8} \left[\log\left(\frac{\lambda}{\gamma} + 1\right)\right]^2.$$

This completes the proof. □

**Lemma 5.** *Let $\gamma > 0$ be fixed. Let $\lambda \in \mathbb{R}$ be fixed such that $\lambda + \gamma > 0$. Assume that $r = r(\alpha)$ is smooth, with $r(0) = 0$. Define $\delta = \frac{\gamma^r}{\gamma^r + 1}$, $p = r(1 - \delta)$, $q = r\delta$. Then*

$$\frac{(\lambda + \gamma)^p + (\lambda + \gamma)^{-q}}{\gamma^p + \gamma^{-q}} \geq 1, \quad (A.8)$$

$$\log\left(\frac{(\lambda + \gamma)^p + (\lambda + \gamma)^{-q}}{\gamma^p + \gamma^{-q}}\right) \geq 0. \quad (A.9)$$

*Proof of Lemma 5.* By Theorem 5, we have

$$\frac{(\lambda+\gamma)^p + (\lambda+\gamma)^{-q}}{\gamma^p + \gamma^{-q}} = \frac{\gamma^p}{\gamma^p + \gamma^{-q}} \left(\frac{\lambda}{\gamma}+1\right)^p + \frac{\gamma^{-q}}{\gamma^p + \gamma^{-q}} \left(\frac{\lambda}{\gamma}+1\right)^{-q}$$

$$= \alpha \left(\frac{\lambda}{\gamma}+1\right)^p + (1-\alpha)\left(\frac{\lambda}{\gamma}+1\right)^{-q}$$

$$\text{where } \alpha = \frac{\gamma^p}{\gamma^p + \gamma^{-q}} = \frac{\gamma^{p+q}}{\gamma^{p+q}+1} = \frac{\gamma^r}{\gamma^r+1} = \delta$$

$$\geq \left(\frac{\lambda}{\gamma}+1\right)^{p\delta} \left(\frac{\lambda}{\gamma}+1\right)^{-q(1-\delta)} = \left(\frac{\lambda}{\gamma}+1\right)^{(p+q)\delta - q} = \left(\frac{\lambda}{\gamma}+1\right)^{r\delta - q} = 1,$$

since $q = r\delta$. This completes the proof. $\square$

*Proof of Theorem 9.* For $\alpha = \beta$, we have

$$\delta = \frac{(\frac{\gamma}{\mu})^r}{(\frac{\gamma}{\mu})^r + 1}, \quad p = r(1-\delta), \quad q = r\delta.$$

Let $\{\lambda_j\}_{j\in\mathbb{N}}$ be the eigenvalues of $\Lambda$. By Theorem 8, we have

$$D_r^{(\alpha,\alpha)}[(A+\gamma I),(B+\mu I)] = \frac{r(\delta-\frac{1}{2})}{\alpha^2} \log\left(\frac{\gamma}{\mu}\right) + \frac{1}{\alpha^2} \log\left(\frac{(\frac{\gamma}{\mu})^p + (\frac{\gamma}{\mu})^{-q}}{2}\right)$$

$$+ \frac{1}{\alpha^2} \log \det \left[\frac{(\Lambda + \frac{\gamma}{\mu}I)^p + (\Lambda + \frac{\gamma}{\mu}I)^{-q}}{(\frac{\gamma}{\mu})^p + (\frac{\gamma}{\mu})^{-q}}\right]$$

$$= \frac{r(\delta-\frac{1}{2})}{\alpha^2} \log\left(\frac{\gamma}{\mu}\right) + \frac{1}{\alpha^2} \log\left(\frac{(\frac{\gamma}{\mu})^p + (\frac{\gamma}{\mu})^{-q}}{2}\right)$$

$$+ \frac{1}{\alpha^2} \sum_{j=1}^{\infty} \log \left(\frac{(\lambda_j + \frac{\gamma}{\mu})^p + (\lambda_j + \frac{\gamma}{\mu})^{-q}}{(\frac{\gamma}{\mu})^p + (\frac{\gamma}{\mu})^{-q}}\right).$$

By Lemma 2, we have

$$\lim_{\alpha \to 0} \frac{r(\delta - \frac{1}{2})}{\alpha^2} \log\left(\frac{\gamma}{\mu}\right) = \frac{[r'(0)]^2}{4} \left[\log\frac{\gamma}{\mu}\right]^2.$$

By Lemma 3, we have

$$\lim_{\alpha \to 0} \frac{1}{\alpha^2} \log\left(\frac{(\frac{\gamma}{\mu})^p + (\frac{\gamma}{\mu})^{-q}}{2}\right) = -\frac{[r'(0)]^2}{8} \left[\log\frac{\gamma}{\mu}\right]^2.$$

By Lemma 4, we have

$$\lim_{\alpha \to 0} \frac{1}{\alpha^2} \log \det \left[ \frac{(\Lambda + \frac{\gamma}{\mu}I)^p + (\Lambda + \frac{\gamma}{\mu}I)^{-q}}{(\frac{\gamma}{\mu})^p + (\frac{\gamma}{\mu})^{-q}} \right]$$

$$= \lim_{\alpha \to 0} \frac{1}{\alpha^2} \sum_{j=1}^{\infty} \log \left[ \frac{(\lambda_j + \frac{\gamma}{\mu})^p + (\lambda_j + \frac{\gamma}{\mu})^{-q}}{(\frac{\gamma}{\mu})^p + (\frac{\gamma}{\mu})^{-q}} \right]$$

$$= \sum_{j=1}^{\infty} \lim_{\alpha \to 0} \frac{1}{\alpha^2} \log \left[ \frac{(\lambda_j + \frac{\gamma}{\mu})^p + (\lambda_j + \frac{\gamma}{\mu})^{-q}}{(\frac{\gamma}{\mu})^p + (\frac{\gamma}{\mu})^{-q}} \right]$$

by Lebesgue's Monotone Convergence Theorem, since

$$\log \left[ \frac{(\lambda_j + \frac{\gamma}{\mu})^p + (\lambda_j + \frac{\gamma}{\mu})^{-q}}{(\frac{\gamma}{\mu})^p + (\frac{\gamma}{\mu})^{-q}} \right] \geq 0 \; \forall j \in \mathbb{N} \; \text{ by Lemma 5}$$

$$= \frac{[r'(0)]^2}{8} \sum_{j=1}^{\infty} \left[ \log\left(\lambda_j + \frac{\gamma}{\mu}\right) - \log\left(\frac{\gamma}{\mu}\right) \right]^2 = \frac{[r'(0)]^2}{8} \sum_{j=1}^{\infty} \left[ \log\left(\lambda_j \frac{\mu}{\gamma} + 1\right) \right]^2.$$

Summing up these three expressions, we obtain

$$\lim_{\alpha \to 0} D_r^{(\alpha,\alpha)}[(A + \gamma I), (B + \mu I)] = \frac{[r'(0)]^2}{8} \left( \left[\log \frac{\gamma}{\mu}\right]^2 + \sum_{j=1}^{\infty} \left[ \log\left(\lambda_j \frac{\mu}{\gamma} + 1\right) \right]^2 \right)$$

$$= \frac{[r'(0)]^2}{8} \left( \left[\log \frac{\gamma}{\mu}\right]^2 + \left\|\log\left(\Lambda \frac{\mu}{\gamma} + I\right)\right\|_{\text{HS}}^2 \right) = \frac{[r'(0)]^2}{8} \left\|\log\left(\Lambda + \frac{\gamma}{\mu}I\right)\right\|_{\text{eHS}}^2$$

$$= \frac{[r'(0)]^2}{8} ||\log[(B + \mu I)^{-1/2}(A + \gamma I)(B + \mu I)^{-1/2}]||_{\text{eHS}}^2$$

$$= \frac{[r'(0)]^2}{8} d_{\text{aiHS}}^2[(A + \gamma I), (B + \mu I)].$$

This completes the proof. □

*Appendix A.3. Proofs for the Alpha Log-Determinant divergences*

In this section, we prove Theorem 10.

**Proof of Theorem 10.** The proof for the cases $\alpha = 0$ and $\alpha = 1$ is a special case of the results discussed at the end of Section 5.3.

Consider now the case $0 < \alpha < 1$. We first note that

$$d_{\text{logdet}}^{1-2\alpha}[(A + \gamma I), (B + \mu I)] = \frac{1}{\alpha(1-\alpha)} \log \left[ \frac{\det_X(\alpha(A + \gamma I) + (1-\alpha)(B + \mu I))}{\det_X(A + \gamma I)^q \det_X(B + \mu I)^{1-q}} \right]$$

$$+ \frac{q - \alpha}{\alpha(1-\alpha)} \log \frac{\gamma}{\mu},$$

where $q = \frac{\alpha\gamma}{\alpha\gamma+(1-\alpha)\mu}$.

By Definition 6, we have

$$D_r^{(\alpha,1-\alpha)}[(A+\gamma I),(B+\mu I)]$$
$$= \frac{1}{\alpha(1-\alpha)}\log\left[\left(\frac{\gamma}{\mu}\right)^{r(\delta-\alpha)}\det_X\left(\alpha\left(\Lambda+\frac{\gamma}{\mu}I\right)^{r(1-\delta)}+(1-\alpha)\left(\Lambda+\frac{\gamma}{\mu}I\right)^{-r\delta}\right)\right]$$
$$= \frac{r(\delta-\alpha)}{\alpha(1-\alpha)}\log\left(\frac{\gamma}{\mu}\right)$$
$$+ \frac{1}{\alpha(1-\alpha)}\log\det_X\left(\alpha\left(\Lambda+\frac{\gamma}{\mu}I\right)^{r(1-\delta)}+(1-\alpha)\left(\Lambda+\frac{\gamma}{\mu}I\right)^{-r\delta}\right).$$

By Proposition 1, we have

$$\det_X\left(\alpha\left(\Lambda+\frac{\gamma}{\mu}I\right)^{r(1-\delta)}+(1-\alpha)\left(\Lambda+\frac{\gamma}{\mu}I\right)^{-r\delta}\right)$$
$$= \det_X\left[\alpha[(A+\gamma I)(B+\mu I)^{-1}]^{r(1-\delta)}+(1-\alpha)[(A+\gamma I)(B+\mu I)^{-1}]^{-r\delta}\right]$$
$$= \det_X[(A+\gamma I)(B+\mu I)^{-1}]^{-r\delta}\det_X[\alpha[(A+\gamma I)(B+\mu I)^{-1}]^r+(1-\alpha)I].$$

In particular, for $r=1$, we have

$$\det_X[\alpha[(A+\gamma I)(B+\mu I)^{-1}]+(1-\alpha)] = \frac{\det_X[\alpha(A+\gamma I)+(1-\alpha)(B+\mu I)]}{\det_X(B+\mu I)}.$$

Thus it follows that

$$\det_X\left(\alpha\left(\Lambda+\frac{\gamma}{\mu}I\right)^{(1-\delta)}+(1-\alpha)\left(\Lambda+\frac{\gamma}{\mu}I\right)^{-\delta}\right)$$
$$= \frac{\det_X[\alpha(A+\gamma I)+(1-\alpha)(B+\mu I)]}{\det_X(A+\gamma I)^\delta\det_X(B+\mu I)^{1-\delta}}.$$

Also for $r=1$, in Definition 6, we have $\delta=\delta(r=1)=\frac{\alpha\gamma}{\alpha\gamma+(1-\alpha)\mu}$. Combining all of these expressions and comparing with the expressions for $d_{\text{logdet}}^{1-2\alpha}$, we obtain the first desired statement.

For $r=-1$, we have

$$D_{-1}^{(\alpha,1-\alpha)}[(A+\gamma I),(B+\mu I)]$$
$$= \frac{-(\delta_{-1}-\alpha)}{\alpha(1-\alpha)}\log\left(\frac{\gamma}{\mu}\right)$$
$$+ \frac{1}{\alpha(1-\alpha)}\log\det_X\left(\alpha\left(\Lambda+\frac{\gamma}{\mu}I\right)^{-(1-\delta_{-1})}+(1-\alpha)\left(\Lambda+\frac{\gamma}{\mu}I\right)^{\delta_{-1}}\right),$$

where $\delta_{-1} = \delta(r = -1) = \frac{\alpha\frac{1}{\gamma}}{\alpha\frac{1}{\gamma}+(1-\alpha)\frac{1}{\mu}} = \frac{\alpha\mu}{\alpha\mu+(1-\alpha)\gamma}$.

Similar to the case $r = 1$, we have

$$\det_X[\alpha[(A+\gamma I)(B+\mu I)^{-1}]^{-1} + (1-\alpha)I] = \frac{\det_X[(1-\alpha)(A+\gamma I) + \alpha(B+\mu I)]}{\det_X(A+\gamma I)}.$$

Thus it follows that

$$\det_X\left(\alpha\left(\Lambda + \frac{\gamma}{\mu}I\right)^{-(1-\delta_{-1})} + (1-\alpha)\left(\Lambda + \frac{\gamma}{\mu}I\right)^{\delta_{-1}}\right)$$
$$= \frac{\det_X[(1-\alpha)(A+\gamma I) + \alpha(B+\mu I)]}{\det_X(A+\gamma I)^{1-\delta_{-1}}\det_X(B+\mu I)^{\delta_{-1}}}.$$

On the other hand, we have

$$d_{\text{logdet}}^{2\alpha-1}[(A+\gamma I), (B+\mu I)] = \frac{1}{\alpha(1-\alpha)}\log\left[\frac{\det_X((1-\alpha)(A+\gamma I) + \alpha(B+\mu I))}{\det_X(A+\gamma I)^p\det_X(B+\mu I)^{1-p}}\right]$$
$$+ \frac{p-(1-\alpha)}{\alpha(1-\alpha)}\log\frac{\gamma}{\mu},$$

where $p = \frac{(1-\alpha)\gamma}{(1-\alpha)\gamma+\alpha\mu} = 1 - \delta_{-1}$. Combining all of these expressions, we obtain the second desired statement, namely

$$D_{-1}^{(\alpha,1-\alpha)}[(A+\gamma I), (B+\mu I)] = d_{\text{logdet}}^{2\alpha-1}[(A+\gamma I), (B+\mu I)].$$

This completes the proof. □

*Appendix A.4. Proofs for the other limiting cases*

In this section, we prove Theorems 11 and 12. We need the following preliminary results.

**Lemma 6.** *Let $\mathcal{H}$ be a separable Hilbert space. Let $A \in \text{Sym}(\mathcal{H}) \cap \text{Tr}(\mathcal{H})$ be such that $A+I > 0$. Then $\forall \alpha \in \mathbb{R}$, the operator $(A+I)^\alpha$ is well defined and $(A+I)^\alpha - I \in \text{Sym}(\mathcal{H}) \cap \text{Tr}(\mathcal{H})$. Equivalently, let $\{\lambda_k\}_{k\in\mathbb{N}}$ be the eigenvalues of $A$, then*

$$\text{tr}[(A+I)^\alpha - I] = \sum_{k=1}^{\infty}[(\lambda_k+1)^\alpha - 1] \tag{A.10}$$

*has a finite value.*

*Proof of Lemma 6.* By Lemma 3 in [17], if $A \in \text{Sym}(\mathcal{H}) \cap \text{Tr}(\mathcal{H})$ and $A + I > 0$, then $\log(A + I) \in \text{Sym}(\mathcal{H}) \cap \text{Tr}(\mathcal{H})$. By definition of the power function, we have

$$(A + I)^\alpha = \exp[\alpha \log(A + I)] = I + \sum_{j=1}^\infty \frac{\alpha^j}{j!} [\log(A + I)]^j.$$

Since $\text{Tr}(\mathcal{H})$ is a Banach algebra under the trace norm, we have

$$||(A + I)^\alpha - I||_{\text{tr}} = \left\| \sum_{j=1}^\infty \frac{\alpha^j}{j!} [\log(A + I)]^j \right\|_{\text{tr}} \leq \sum_{j=1}^\infty \frac{|\alpha|^j}{j!} \|\log(A + I)\|_{\text{tr}}^j$$

$$= \exp(|\alpha| \, \|\log(A + I)\|_{\text{tr}}) - 1 < \infty.$$

Thus $(A + I)^\alpha - I \in \text{Tr}(\mathcal{H})$. The equivalent statement is then obvious. This completes the proof. $\square$

**Lemma 7.** *Let $\mathcal{H}$ be a separable Hilbert space. Assume that $(A + \gamma I) \in \text{PTr}(\mathcal{H})$. Then for any $\alpha \in \mathbb{R}$, we have $(A + \gamma I)^\alpha - \gamma^\alpha I \in \text{Sym}(\mathcal{H}) \cap \text{Tr}(\mathcal{H})$ and*

$$\text{tr}[(A + \gamma I)^\alpha - \gamma^\alpha I] = \gamma^\alpha \text{tr}\left[\left(\frac{A}{\gamma} + I\right)^\alpha - I\right], \tag{A.11}$$

$$\text{tr}_X[(A + \gamma I)^\alpha] = \gamma^\alpha \left(1 + \text{tr}\left[\left(\frac{A}{\gamma} + I\right)^\alpha - I\right]\right). \tag{A.12}$$

*Proof of Lemma 7.* By definition of the power function, we have

$$(A + \gamma I)^\alpha = \exp[\alpha \log(A + \gamma I)] = \exp\left[(\alpha \log \gamma)I + \alpha \log\left(\frac{A}{\gamma} + I\right)\right]$$

$$= \gamma^\alpha \left(\frac{A}{\gamma} + I\right)^\alpha = \gamma^\alpha \left[\left(\frac{A}{\gamma} + I\right)^\alpha - I\right] + \gamma^\alpha I,$$

where $\left[\left(\frac{A}{\gamma} + I\right)^\alpha - I\right] \in \text{Tr}(\mathcal{H})$ by Lemma 6. Thus it follows that $(A+\gamma I)^\alpha - \gamma^\alpha I \in \text{Sym}(\mathcal{H}) \cap \text{Tr}(\mathcal{H})$ and

$$\text{tr}[(A + \gamma I)^\alpha - \gamma^\alpha I] = \gamma^\alpha \text{tr}\left[\left(\frac{A}{\gamma} + I\right)^\alpha - I\right],$$

which is the first identity. By definition of the extended trace

$$\text{tr}_X[(A + \gamma I)^\alpha] = \text{tr}_X([(A + \gamma I)^\alpha - \gamma^\alpha I] + \gamma^\alpha I) = \gamma^\alpha \text{tr}\left[\left(\frac{A}{\gamma} + I\right)^\alpha - I\right] + \gamma^\alpha,$$

which is the second identity. This completes the proof. $\square$

**Lemma 8.** *Let $(A + \gamma I), (B + \mu I) \in \mathrm{PTr}(\mathcal{H})$. Let $\Lambda + \frac{\gamma}{\mu} I = (B + \mu I)^{-1/2}(A + \gamma I)(B + \mu I)^{-1/2}$. Then for any $\alpha \in \mathbb{R}$,*

$$\mathrm{tr}_X[(A+\gamma I)(B+\mu I)^{-1}]^\alpha = \mathrm{tr}_X\left[\left(\Lambda + \frac{\gamma}{\mu}\right)^\alpha\right]$$

$$= \mathrm{tr}_X[(B+\mu I)^{-1}(A+\gamma I)]^\alpha. \quad (\text{A.13})$$

$$\mathrm{det}_X[(A+\gamma I)(B+\mu I)^{-1}]^\alpha = \mathrm{det}_X\left[\left(\Lambda + \frac{\gamma}{\mu}\right)^\alpha\right]$$

$$= \mathrm{det}_X[(B+\mu I)^{-1}(A+\gamma I)]^\alpha. \quad (\text{A.14})$$

*Proof of Lemma 8.* By Proposition 1, we have

$$[(A+\gamma I)(B+\mu I)^{-1}]^\alpha = (B+\mu I)^{1/2}\left(\Lambda + \frac{\gamma}{\mu}\right)^\alpha (B+\mu I)^{-1/2}.$$

Similarly,

$$[(B+\mu I)^{-1}(A+\gamma I)]^\alpha = (B+\mu I)^{-1/2}\left(\Lambda + \frac{\gamma}{\mu}\right)^\alpha (B+\mu I)^{1/2}.$$

By the commutativity of the $\mathrm{tr}_X$ operation (Lemma 4 in [17]), we then have

$$\mathrm{tr}_X[(A+\gamma I)(B+\mu I)^{-1}]^\alpha = \mathrm{tr}_X\left[\left(\Lambda + \frac{\gamma}{\mu}\right)^\alpha\right] = \mathrm{tr}_X[(B+\mu I)^{-1}(A+\gamma I)]^\alpha.$$

Similarly, by the product property of the $\mathrm{det}_X$ operation (Proposition 4 in [17]),

$$\mathrm{det}_X[(A+\gamma I)(B+\mu I)^{-1}]^\alpha = \mathrm{det}_X\left[\left(\Lambda + \frac{\gamma}{\mu}\right)^\alpha\right] = \mathrm{det}_X[(B+\mu I)^{-1}(A+\gamma I)]^\alpha.$$

This completes the proof. □

**Lemma 9.** *Assume that $\lambda > 0, \gamma > 0, \alpha > 0$ are fixed. Assume that $r = r(\beta)$ is smooth. Then for $\delta = \frac{\alpha \gamma^r}{\alpha \gamma^r + \beta}$, $p = r(1-\delta)$, $q = r\delta$, we have*

$$\lim_{\beta \to 0} \frac{1}{\alpha\beta} \log\left(\frac{\alpha\lambda^p + \beta\lambda^{-q}}{\alpha + \beta}\right) = \frac{1}{\alpha^2}\left((\log \lambda)\frac{r(0)}{\gamma^{r(0)}} + \lambda^{-r(0)} - 1\right). \quad (\text{A.15})$$

*In particular, for $\lambda = \gamma$, we have*

$$\lim_{\beta \to 0} \frac{1}{\alpha\beta} \log\left(\frac{\alpha\gamma^p + \beta\gamma^{-q}}{\alpha + \beta}\right) = \frac{1}{\alpha^2}\left([(\log \gamma)r(0) + 1]\gamma^{-r(0)} - 1\right). \quad (\text{A.16})$$

*Proof of Lemma 9.* We have for $\alpha > 0$, $\lim_{\beta \to 0} \delta = 1$, $\lim_{\beta \to 0} p = 0$, $\lim_{\beta \to 0} q = r(0)$, so that $\lim_{\beta \to 0}(\alpha\lambda^p + \beta\lambda^{-q}) = \alpha$. With $p = r(1-\delta) = \frac{r\beta}{\alpha\gamma^r + \beta}$, we have

$$\frac{\partial p}{\partial \beta} = \frac{(\frac{\partial r}{\partial \beta}\beta + r)(\alpha\gamma^r + \beta) - r\beta(\alpha\gamma^r \log\gamma \frac{\partial r}{\partial \beta} + 1)}{(\alpha\gamma^r + \beta)^2},$$

$$\lim_{\beta \to 0} \frac{\partial p}{\partial \beta} = \frac{r(0)}{\alpha\gamma^{r(0)}}.$$

With $q = r\delta = \frac{r\alpha\gamma^r}{\alpha\gamma^r + \beta}$, we have

$$\frac{\partial q}{\partial \beta} = \frac{(\frac{\partial r}{\partial \beta}\alpha\gamma^r + r\alpha\gamma^r \log\gamma \frac{\partial r}{\partial \beta})(\alpha\gamma^r + \beta) - r\alpha\gamma^r(\alpha\gamma^r \log\gamma \frac{\partial r}{\partial \beta} + 1)}{(\alpha\gamma^r + \beta)^2},$$

$$\lim_{\beta \to 0} \frac{\partial q}{\partial \beta} = \frac{\partial r}{\partial \beta}(0) - \frac{r(0)}{\alpha\gamma^{r(0)}}.$$

The required limit is of the form $\frac{0}{0}$ and L'Hopital's rule can be applied to give

$$\lim_{\beta \to 0} \frac{1}{\alpha\beta} \log\left(\frac{\alpha\lambda^p + \beta\lambda^{-q}}{\alpha + \beta}\right)$$

$$= \frac{1}{\alpha} \lim_{\beta \to 0} \frac{\alpha + \beta}{\alpha\lambda^p + \beta\lambda^{-q}} \frac{[\alpha\lambda^p(\log\lambda)\frac{\partial p}{\partial \beta} + \lambda^{-q} - \beta\lambda^{-q}(\log\lambda)\frac{\partial q}{\partial \beta}](\alpha + \beta) - (\alpha\lambda^p + \beta\lambda^{-q})}{(\alpha + \beta)^2}$$

$$= \frac{\alpha(\log\lambda)\frac{\partial p}{\partial \beta}(0) + \lambda^{-r(0)} - 1}{\alpha^2} = \frac{1}{\alpha^2}\left((\log\lambda)\frac{r(0)}{\gamma^{r(0)}} + \lambda^{-r(0)} - 1\right).$$

This completes the proof. $\square$

**Lemma 10.** *Assume that $\gamma > 0, \alpha > 0$ are fixed. Assume that $\lambda \in \mathbb{R}$ is also fixed, such that $\lambda + \gamma > 0$. Assume that $r = r(\beta)$ is smooth. Then for $\delta = \frac{\alpha\gamma^r}{\alpha\gamma^r + \beta}$, $p = r(1-\delta)$, $q = r\delta$, we have*

$$\lim_{\beta \to 0} \frac{1}{\alpha\beta} \log\left(\frac{\alpha(\lambda+\gamma)^p + \beta(\lambda+\gamma)^{-q}}{\alpha\gamma^p + \beta\gamma^{-q}}\right)$$

$$= \frac{1}{\alpha^2}\left[\log\left(\frac{\lambda}{\gamma} + 1\right)\frac{r(0)}{\gamma^{r(0)}} + (\lambda+\gamma)^{-r(0)} - \gamma^{-r(0)}\right]. \quad (A.17)$$

*Proof of Lemma 10.* By Lemma 9, we have

$$\lim_{\beta \to 0} \frac{1}{\alpha\beta} \log\left(\frac{\alpha(\lambda+\gamma)^p + \beta(\lambda+\gamma)^{-q}}{\alpha\gamma^p + \beta\gamma^{-q}}\right)$$

$$= \lim_{\beta \to 0} \frac{1}{\alpha\beta} \log\left(\frac{\alpha(\lambda+\gamma)^p + \beta(\lambda+\gamma)^{-q}}{\alpha + \beta}\right) - \lim_{\beta \to 0} \frac{1}{\alpha\beta} \log\left(\frac{\alpha\gamma^p + \beta\gamma^{-q}}{\alpha + \beta}\right)$$

$$= \frac{1}{\alpha^2}\left((\log(\lambda+\gamma))\frac{r(0)}{\gamma^{r(0)}} + (\lambda+\gamma)^{-r(0)} - 1\right) - \frac{1}{\alpha^2}\left((\log\gamma)\frac{r(0)}{\gamma^{r(0)}} + \gamma^{-r(0)} - 1\right)$$

$$= \frac{1}{\alpha^2}\left[\log\left(\frac{\lambda}{\gamma} + 1\right)\frac{r(0)}{\gamma^{r(0)}} + (\lambda+\gamma)^{-r(0)} - \gamma^{-r(0)}\right].$$

This completes the proof. □

**Lemma 11.** *Assume that $\gamma > 0, \alpha > 0$ are fixed. Assume that $\lambda \in \mathbb{R}$ is also fixed, such that $\lambda + \gamma > 0$. Assume that $r = r(\beta)$ is smooth. Then for $\delta = \frac{\alpha\gamma^r}{\alpha\gamma^r + \beta}$, $p = r(1-\delta)$, $q = r\delta$, we have*

$$\frac{\alpha(\lambda+\gamma)^p + \beta(\lambda+\gamma)^{-q}}{\alpha\gamma^p + \beta\gamma^{-q}} \geq 1, \tag{A.18}$$

$$\log\left(\frac{\alpha(\lambda+\gamma)^p + \beta(\lambda+\gamma)^{-q}}{\alpha\gamma^p + \beta\gamma^{-q}}\right) \geq 0. \tag{A.19}$$

*Proof of Lemma 11.* We proceed as in the proof of Lemma 5, by applying Theorem 5 as follows

$$\frac{\alpha(\lambda+\gamma)^p + \beta(\lambda+\gamma)^{-q}}{\alpha\gamma^p + \beta\gamma^{-q}} = \frac{\alpha\gamma^p}{\alpha\gamma^p + \beta\gamma^{-q}}\left(\frac{\lambda}{\gamma}+1\right)^p + \frac{\beta\gamma^{-q}}{\alpha\gamma^p + \beta\gamma^{-q}}\left(\frac{\lambda}{\gamma}+1\right)^{-q}$$

$$= s\left(\frac{\lambda}{\gamma}+1\right)^p + (1-s)\left(\frac{\lambda}{\gamma}+1\right)^{-q},$$

where $s = \dfrac{\alpha\gamma^p}{\alpha\gamma^p + \beta\gamma^{-q}} = \dfrac{\alpha\gamma^{p+q}}{\alpha\gamma^{p+q}+\beta} = \dfrac{\alpha\gamma^r}{\alpha\gamma^r+\beta} = \delta,$

$$\geq \left(\frac{\lambda}{\gamma}+1\right)^{p\delta}\left(\frac{\lambda}{\gamma}+1\right)^{-q(1-\delta)} = \left(\frac{\lambda}{\gamma}+1\right)^{(p+q)\delta-q} = \left(\frac{\lambda}{\gamma}+1\right)^{r\delta-q} = 1,$$

since $r\delta = q$. This completes the proof. □

**Lemma 12.** *Assume that $\gamma > 0, \alpha > 0$ are fixed. Assume that $r = r(\beta)$ is smooth. Then for $\delta = \frac{\alpha\gamma^r}{\alpha\gamma^r+\beta}$,*

$$\lim_{\beta \to 0} \frac{r(\delta - \frac{\alpha}{\alpha+\beta})}{\alpha\beta} = \frac{1}{\alpha^2}r(0)[-\gamma^{-r(0)}+1]. \tag{A.20}$$

*Proof of Lemma 12.* We first have

$$\frac{\partial \delta}{\partial \beta} = \frac{\alpha\gamma^r \log\gamma \frac{\partial r}{\partial \beta}(\alpha\gamma^r+\beta) - \alpha\gamma^r(\alpha\gamma^r \log\gamma \frac{\partial r}{\partial \beta}+1)}{(\alpha\gamma^r+\beta)^2}$$

$$\lim_{\beta \to 0}\frac{\partial \delta}{\partial \beta} = -\frac{1}{\alpha\gamma^{r(0)}}.$$

Since the required limit has the form $\frac{0}{0}$, we apply L'Hopital's rule to get

$$\lim_{\beta \to 0}\frac{r(\delta - \frac{\alpha}{\alpha+\beta})}{\alpha\beta} = \lim_{\beta \to 0}\frac{1}{\alpha}\left[\frac{\partial r}{\partial \beta}\left(\delta - \frac{\alpha}{\alpha+\beta}\right) + r\left(\frac{\partial \delta}{\partial \beta} + \frac{\alpha}{(\alpha+\beta)^2}\right)\right]$$

$$= \frac{1}{\alpha}\left[r(0)\left(-\frac{1}{\alpha\gamma^{r(0)}}+\frac{1}{\alpha}\right)\right] = \frac{1}{\alpha^2}r(0)[-\gamma^{-r(0)}+1].$$

This completes the proof. □

**Proof of Theorem 11.** Let $\{\lambda_j\}_{j=1}^\infty$ be the eigenvalues of $\Lambda$. By Theorem 8, we have

$$D_r^{(\alpha,\beta)}[(A+\gamma I),(B+\mu I)] = \frac{r(\delta - \frac{\alpha}{\alpha+\beta})}{\alpha\beta}\log\left(\frac{\gamma}{\mu}\right) + \frac{1}{\alpha\beta}\log\left(\frac{\alpha(\frac{\gamma}{\mu})^p + \beta(\frac{\gamma}{\mu})^{-q}}{\alpha+\beta}\right)$$

$$+ \frac{1}{\alpha\beta}\log\det\left(\frac{\alpha(\Lambda + \frac{\gamma}{\mu}I)^p + \beta(\Lambda + \frac{\gamma}{\mu}I)^{-q}}{\alpha(\frac{\gamma}{\mu})^p + \beta(\frac{\gamma}{\mu})^{-q}}\right)$$

$$= \frac{r(\delta - \frac{\alpha}{\alpha+\beta})}{\alpha\beta}\log\left(\frac{\gamma}{\mu}\right) + \frac{1}{\alpha\beta}\log\left(\frac{\alpha(\frac{\gamma}{\mu})^p + \beta(\frac{\gamma}{\mu})^{-q}}{\alpha+\beta}\right)$$

$$+ \frac{1}{\alpha\beta}\sum_{j=1}^\infty \log\left(\frac{\alpha(\lambda_j + \frac{\gamma}{\mu})^p + \beta(\lambda_j + \frac{\gamma}{\mu})^{-q}}{\alpha(\frac{\gamma}{\mu})^p + \beta(\frac{\gamma}{\mu})^{-q}}\right),$$

where $p = p(\beta) = r(1-\delta) = \frac{r\beta}{\alpha(\frac{\gamma}{\mu})^r+\beta}$, $q = q(\beta) = r\delta = \frac{r\alpha(\frac{\gamma}{\mu})^r}{\alpha(\frac{\gamma}{\mu})^r+\beta}$.

For $\alpha > 0$ fixed, as functions of $\beta$, we have

$$\lim_{\beta\to 0} p(\beta) = 0, \quad \lim_{\beta\to 0} q(\beta) = r(0).$$

For simplicity, in the following, we replace $\frac{\gamma}{\mu}$ by $\gamma$. By Lemma 9,

$$\lim_{\beta\to 0}\frac{1}{\alpha\beta}\log\left(\frac{\alpha\gamma^p + \beta\gamma^{-q}}{\alpha+\beta}\right) = \frac{1}{\alpha^2}\left([(\log\gamma)r(0)+1]\gamma^{-r(0)} - 1\right).$$

By Lemma 10,

$$\lim_{\beta\to 0}\frac{1}{\alpha\beta}\log\left(\frac{\alpha(\lambda_j+\gamma)^p + \beta(\lambda_j+\gamma)^{-q}}{\alpha\gamma^p + \beta\gamma^{-q}}\right)$$
$$= \frac{1}{\alpha^2}\left[\log\left(\frac{\lambda_j}{\gamma}+1\right)\frac{r(0)}{\gamma^{r(0)}} + (\lambda_j+\gamma)^{-r(0)} - \gamma^{-r(0)}\right].$$

By Lemma 11, we have $\log\left(\frac{\alpha(\lambda_j+\gamma)^p+\beta(\lambda_j+\gamma)^{-q}}{\alpha\gamma^p+\beta\gamma^{-q}}\right) \geq 0 \; \forall j \in \mathbb{N}$, so that by Lebesgue's Monotone Convergence Theorem, we obtain

$$\lim_{\beta\to 0}\frac{1}{\alpha\beta}\sum_{j=1}^\infty \log\left(\frac{\alpha(\lambda_j+\gamma)^p + \beta(\lambda_j+\gamma)^{-q}}{\alpha\gamma^p + \beta\gamma^{-q}}\right)$$
$$= \sum_{j=1}^\infty \lim_{\beta\to 0}\frac{1}{\alpha\beta}\log\left(\frac{\alpha(\lambda_j+\gamma)^p + \beta(\lambda_j+\frac{\gamma}{\mu})^{-q}}{\alpha\gamma^p + \beta\gamma^{-q}}\right)$$
$$= \frac{1}{\alpha^2}\sum_{j=1}^\infty\left[\log\left(\frac{\lambda_j}{\gamma}+1\right)\frac{r(0)}{\gamma^{r(0)}} + (\lambda_j+\gamma)^{-r(0)} - \gamma^{-r(0)}\right].$$

By Lemma 12

$$\log(\gamma)\lim_{\beta\to 0}\frac{r(\delta - \frac{\alpha}{\alpha+\beta})}{\alpha\beta} = \frac{1}{\alpha^2}r(0)[-\gamma^{-r(0)} + 1]\log(\gamma).$$

Combining all three expressions, we obtain the desired limit as the sum

$$\frac{1}{\alpha^2}[\gamma^{-r(0)} + r(0)\log(\gamma) - 1]$$
$$+ \frac{1}{\alpha^2}\left\{\frac{r(0)}{\gamma^{r(0)}}\sum_{j=1}^{\infty}\log\left(\frac{\lambda_j}{\gamma}+1\right) + \sum_{j=1}^{\infty}\left[\frac{1}{(\lambda_j+\gamma)^{r(0)}} - \frac{1}{\gamma^{r(0)}}\right]\right\}. \quad (A.21)$$

By Lemmas 6 and 7, we have

$$\sum_{j=1}^{\infty}\left[\frac{1}{(\lambda_j+\gamma)^{r(0)}} - \frac{1}{\gamma^{r(0)}}\right] = \gamma^{-r(0)}\sum_{j=1}^{\infty}\left[\left(\frac{\lambda_j}{\gamma}+1\right)^{-r(0)} - 1\right]$$
$$= \gamma^{-r(0)}\mathrm{tr}\left[\left(\frac{\Lambda}{\gamma}+I\right)^{-r(0)} - I\right] = \mathrm{tr}[(\Lambda+\gamma I)^{-r(0)} - \gamma^{-r(0)}I].$$

Thus it follows that

$$\gamma^{-r(0)} - 1 + \sum_{j=1}^{\infty}\left[\frac{1}{(\lambda_j+\gamma)^{r(0)}} - \frac{1}{\gamma^{r(0)}}\right]$$
$$= \gamma^{-r(0)} - 1 + \mathrm{tr}[(\Lambda+\gamma I)^{-r(0)} - \gamma^{-r(0)}I] = \mathrm{tr}_X[(\Lambda+\gamma I)^{-r(0)} - I].$$

Furthermore,

$$\frac{r(0)}{\gamma^{r(0)}}\sum_{j=1}^{\infty}\log\left(\frac{\lambda_j}{\gamma}+1\right) = r(0)\gamma^{-r(0)}\log\det\left(\frac{\Lambda}{\gamma}+I\right)$$
$$= r(0)\gamma^{-r(0)}\log\det_X(\Lambda+\gamma I) - r(0)\gamma^{-r(0)}\log\gamma$$
$$= -\gamma^{-r(0)}\log\det_X(\Lambda+\gamma I)^{-r(0)} - r(0)\gamma^{-r(0)}\log\gamma.$$

Plugging the last two expressions into (A.21), we obtain the desired limit as

$$\frac{1}{\alpha^2}\left\{r(0)(1-\gamma^{-r(0)})\log\gamma\right\}$$
$$+ \frac{1}{\alpha^2}\left\{\mathrm{tr}_X[(\Lambda+\gamma I)^{-r(0)} - I] - \gamma^{-r(0)}\log\det_X(\Lambda+\gamma I)^{-r(0)}\right\}. \quad (A.22)$$

We now replace $\gamma$ by $\frac{\gamma}{\mu}$. We have by Lemma 8,

$$\mathrm{tr}_X\left[\left(\Lambda + \frac{\gamma}{\mu}I\right)^{-r(0)}\right] = \mathrm{tr}_X[(B+\mu I)^{-1}(A+\gamma I)]^{-r(0)}$$
$$= \mathrm{tr}_X[(A+\gamma I)^{-1}(B+\mu I)]^{r(0)},$$
$$\det_X\left(\Lambda + \frac{\gamma}{\mu}I\right)^{-r(0)} = \det_X[(B+\mu I)^{-1}(A+\gamma I)]^{-r(0)}$$
$$= \det_X\left[(A+\gamma I)^{-1}(B+\mu I)\right]^{r(0)}.$$

Then (A.22) becomes

$$\frac{r(0)}{\alpha^2}\left[\left(\frac{\mu}{\gamma}\right)^{r(0)}-1\right]\log\frac{\mu}{\gamma}+\frac{1}{\alpha^2}\mathrm{tr}_X([(A+\gamma I)^{-1}(B+\mu I)]^{r(0)}-I)$$

$$-\frac{1}{\alpha^2}\left(\frac{\mu}{\gamma}\right)^{r(0)}\log\det{}_X[(A+\gamma I)^{-1}(B+\mu I)]^{r(0)}.$$

This completes the proof of the theorem. $\square$

*Proof of Theorem 12.* The dual symmetry in Theorem 13 gives

$$\lim_{\alpha\to 0}D_r^{(\alpha,\beta)}[(A+\gamma I),(B+\mu I)]=\lim_{\alpha\to 0}D_r^{(\beta,\alpha)}[(B+\mu I),(A+\gamma I)].$$

The limit on the right hand side then follows from Theorem 11. $\square$

*Appendix A.5. Proofs of the properties of the Alpha-Beta Log-Determinant divergences*

In this section, we prove Theorems 13, 14, 15, 16, 17, and 18. For the case $\alpha=\beta=0$, we have $D_0^{(0,0)}[(A+\gamma I),(B+\mu I)]=\frac{1}{2}d_{\mathrm{aiHS}}^2[(A+\gamma I),(B+\mu I)]$, with $d_{\mathrm{aiHS}}$ being the affine-invariant Riemannian distance on $\mathrm{PTr}(\mathcal{H})$. Thus these properties are either automatic or straightforward to verify. We thus focus on the three cases ($\alpha>0,\beta>0$), ($\alpha>0,\beta=0$), and ($\alpha=0,\beta>0$).

*Proof of Theorem 13 (Dual symmetry).* For the case $\alpha>0,\beta=0$ and $\alpha=0,\beta>0$, from Eqs. (10) and (11), we immediately have

$$D_r^{(\alpha,0)}[(A+\gamma I),(B+\mu I)]=\frac{r}{\alpha^2}\left[\left(\frac{\mu}{\gamma}\right)^r-1\right]\log\frac{\mu}{\gamma}$$

$$+\frac{1}{\alpha^2}\mathrm{tr}_X([(A+\gamma I)^{-1}(B+\mu I)]^r-I)$$

$$-\frac{1}{\alpha^2}\left(\frac{\mu}{\gamma}\right)^r\log\det{}_X[(A+\gamma I)^{-1}(B+\mu I)]^r$$

$$=D_r^{(0,\alpha)}[(B+\mu I),(A+\gamma I)].$$

Consider now the case $\alpha>0,\beta>0$. Write $\delta=\delta(\alpha,\beta)$ to emphasize its dependence on $\alpha$ and $\beta$, we have $\delta(\alpha,\beta)=\frac{\alpha\gamma^r}{\alpha\gamma^r+\beta\mu^r}$ in $D_r^{(\alpha,\beta)}[(A+\gamma I),(B+\mu I)]$. Then for $D_r^{(\beta,\alpha)}[(B+\mu I),(A+\gamma I)]$, we have

$$\delta(\beta,\alpha)=\frac{\beta\mu^r}{\alpha\gamma^r+\beta\mu^r}=1-\delta(\alpha,\beta),\quad 1-\delta(\beta,\alpha)=\delta(\alpha,\beta),$$

$$\delta(\beta,\alpha)-\frac{\beta}{\alpha+\beta}=1-\delta(\alpha,\beta)-\frac{\beta}{\alpha+\beta}=-\left(\delta(\alpha,\beta)-\frac{\alpha}{\alpha+\beta}\right).$$

By Definition 1, we have

$$D_r^{(\beta,\alpha)}[(B+\mu I),(A+\gamma I)]$$
$$=\frac{1}{\alpha\beta}\log\left(\frac{\mu}{\gamma}\right)^{r(\delta(\beta,\alpha)-\frac{\beta}{\alpha+\beta})}$$
$$+\frac{1}{\alpha\beta}\log\det{}_X\left(\frac{\beta[(B+\mu I)(A+\gamma I)^{-1}]^{r(1-\delta(\beta,\alpha))}+\alpha[(B+\mu I)(A+\gamma I)^{-1}]^{-r\delta(\beta,\alpha)}}{\alpha+\beta}\right)$$
$$=\frac{1}{\alpha\beta}\log\left(\frac{\gamma}{\mu}\right)^{r(\delta(\alpha,\beta)-\frac{\alpha}{\alpha+\beta})}$$
$$+\frac{1}{\alpha\beta}\log\det{}_X\left(\frac{\beta[(A+\gamma I)(B+\mu I)^{-1}]^{-r\delta(\alpha,\beta)}+\alpha[(A+\gamma I)(B+\mu I)^{-1}]^{r(1-\delta(\alpha,\beta))}}{\alpha+\beta}\right)$$
$$=D_r^{(\alpha,\beta)}[(A+\gamma I),(B+\mu I)].$$

This completes the proof of the theorem. $\square$

*Proof of Theorem 14 (Dual invariance under inversion).* We have

$$(A+\gamma I)^{-1}=\frac{1}{\gamma}I-\frac{A}{\gamma}(A+\gamma I)^{-1},\quad (B+\mu I)^{-1}=\frac{1}{\mu}I-\frac{B}{\mu}(B+\mu I)^{-1},$$
$$(B+\mu I)^{1/2}(A+\gamma I)^{-1}(B+\mu I)^{1/2}=[(B+\mu I)^{-1/2}(A+\gamma I)(B+\mu I)^{-1/2}]^{-1}.$$

Consider the case $\alpha>0, \beta>0$. By Definition 1, we have

$$D_r^{(\alpha,\beta)}[(A+\gamma I)^{-1},(B+\mu I)^{-1}]$$
$$=\frac{1}{\alpha\beta}\log\left(\frac{1/\gamma}{1/\mu}\right)^{r(\delta_2-\frac{\alpha}{\alpha+\beta})}+\frac{1}{\alpha\beta}\log\det{}_X\left(\frac{\alpha(\Lambda+\frac{\gamma}{\mu}I)^{-r(1-\delta_2)}+\beta(\Lambda+\frac{\gamma}{\mu})^{r\delta_2}}{\alpha+\beta}\right)$$

where $\delta_2=\frac{\alpha(1/\gamma)^r}{\alpha(1/\gamma)^r+\beta(1/\mu)^r}=\frac{\alpha\mu^r}{\alpha\mu^r+\beta\gamma^r}=\delta(-r)$. Thus

$$D_r^{(\alpha,\beta)}[(A+\gamma I)^{-1},(B+\mu I)^{-1}]=D_{-r}^{(\alpha,\beta)}[(A+\gamma I),(B+\mu I)].$$

Consider the case $\alpha=0, \beta>0$ (the case $\alpha>0, \beta=0$ then follows by dual symme-

try). We have

$$]D_r^{(0,\beta)}[(A+\gamma I)^{-1},(B+\mu I)^{-1}] = \frac{r}{\beta^2}\left[\left(\frac{1/\gamma}{1/\mu}\right)^r - 1\right]\log\frac{1/\gamma}{1/\mu}$$

$$+ \frac{1}{\beta^2}\mathrm{tr}_X([(B+\mu I)(A+\gamma I)^{-1}]^r - I) - \frac{1}{\beta^2}\left(\frac{1/\gamma}{1/\mu}\right)^r \log\det{}_X[(B+\mu I)(A+\gamma I)^{-1}]^r$$

$$= -\frac{r}{\beta^2}\left[\left(\frac{\gamma}{\mu}\right)^{-r} - 1\right]\log\frac{\gamma}{\mu} + \frac{1}{\beta^2}\mathrm{tr}_X([(A+\gamma I)(B+\mu I)^{-1}]^{-r} - I)$$

$$- \frac{1}{\beta^2}\left(\frac{\gamma}{\mu}\right)^{-r}\log\det{}_X[(A+\gamma I)(B+\mu I)^{-1}]^{-r}.$$

By Lemma 8, we have

$$\mathrm{tr}_X[(A+\gamma I)(B+\mu I)^{-1}]^{-r} = \mathrm{tr}_X\left[\left(\Lambda+\frac{\gamma}{\mu}\right)^{-r}\right] = \mathrm{tr}_X[(B+\mu I)^{-1}(A+\gamma I)]^{-r},$$

$$\det{}_X[(A+\gamma I)(B+\mu I)^{-1}]^{-r} = \det{}_X\left[\left(\Lambda+\frac{\gamma}{\mu}\right)^{-r}\right] = \det{}_X[(B+\mu I)^{-1}(A+\gamma I)]^{-r}.$$

Thus it follows that

$$D_r^{(0,\beta)}[(A+\gamma I)^{-1},(B+\mu I)^{-1}]$$

$$= -\frac{r}{\beta^2}\left[\left(\frac{\gamma}{\mu}\right)^{-r} - 1\right]\log\frac{\gamma}{\mu} + \frac{1}{\beta^2}\mathrm{tr}_X([(B+\mu I)^{-1}(A+\gamma I)]^{-r} - I)$$

$$- \frac{1}{\beta^2}\left(\frac{\gamma}{\mu}\right)^{-r}\log\det{}_X[(B+\mu I)^{-1}(A+\gamma I)]^{-r}$$

$$= D_{-r}^{(0,\beta)}[(A+\gamma I),(B+\mu I)].$$

This completes the proof. $\square$

***Proof of Theorem 15 (Affine-invariance).*** We have for $(A+\gamma I) \in \mathrm{PTr}(\mathcal{H})$ and $(C+\nu I) \in \mathrm{Tr}_X(\mathcal{H})$, $\nu \neq 0$,

$$(C+\nu I)(A+\gamma I)(C+\nu I)^*$$
$$= CAC^* + \nu(CA + AC^*) + \nu^2 A + \gamma CC^* + \gamma\nu(C+C^*) + \gamma\nu^2 I \in \mathrm{Tr}_X(\mathcal{H}).$$

Since $(C+\nu I)$ is assumed to be invertible, the operator $(C+\nu I)(A+\gamma I)(C+\nu I)^*$ is also invertible, with inverse $[(C+\nu I)^*]^{-1}(A+\gamma I)^{-1}(C+\nu I)^{-1}$. Furthermore,

$\forall x \in \mathcal{H}$,

$$\langle x, (C+\nu I)(A+\gamma I)(C+\nu I)^* x \rangle = \langle (C+\nu I)^* x, (A+\gamma I)(C+\nu I)^* x \rangle$$
$$\geq M_A \|(C+\nu I)^* x\| \geq 0,$$

with equality if and only if $(C+\nu I)^* x = 0 \iff x = 0$. Thus $(C+\nu I)(A+\gamma I)(C+\nu I)^*$ is strictly positive. Together with its invertibility, this shows that this is a positive definite operator. Hence $(C+\nu I)(A+\gamma I)(C+\nu I)^* \in \mathrm{PTr}(\mathcal{H})$.

For two operators $(A+\gamma I), (B+\mu I) \in \mathrm{PTr}(\mathcal{H})$, we then have

$$[(C+\nu I)(A+\gamma I)(C+\nu I)^*][(C+\nu I)(B+\mu I)(C+\nu I)^*]^{-1}$$
$$= (C+\nu I)[(A+\gamma I)(B+\mu I)^{-1}](C+\nu I)^{-1}.$$

Then for any $p \in \mathbb{R}$, we have

$$([(C+\nu I)(A+\gamma I)(C+\nu I)^*][(C+\nu I)(B+\mu I)(C+\nu I)^*]^{-1})^p$$
$$= (C+\nu I)[(A+\gamma I)(B+\mu I)^{-1}]^p (C+\nu I)^{-1}.$$

Thus for any $a, b > 0$ and any $p, q \in \mathbb{R}$,

$$a([(C+\nu I)(A+\gamma I)(C+\nu I)^*][(C+\nu I)(B+\mu I)(C+\nu I)^*]^{-1})^p$$
$$+ b([(C+\nu I)(A+\gamma I)(C+\nu I)^*][(C+\nu I)(B+\mu I)(C+\nu I)^*]^{-1})^q$$
$$= (C+\nu I)(a[(A+\gamma I)(B+\mu I)^{-1}]^p + b[(A+\gamma I)(B+\mu I)^{-1}]^q)(C+\nu I)^{-1}.$$

By the definition of $D_r^{(\alpha,\beta)}$ and the following invariances of the extended Fredholm determinant $\det_X$ as well as of the extended trace operation $\mathrm{tr}_X$, namely,

$$\det\nolimits_X[C(A+\gamma I)C^{-1}] = \det\nolimits_X[(A+\gamma I)],$$
$$\mathrm{tr}_X[C(A+\gamma I)C^{-1}] = \mathrm{tr}_X[(A+\gamma I)],$$

for $A+\gamma I \in \mathrm{Tr}_X(\mathcal{H})$, $\gamma \neq 0$, and $C \in \mathcal{L}(\mathcal{H})$ invertible (Lemma 5 in [17]), we then obtain the desired affine invariance for $D_r^{(\alpha,\beta)}$, namely

$$D_r^{(\alpha,\beta)}[(C+\nu I)(A+\gamma I)(C+\nu I)^*, (C+\nu I)(B+\mu I)(C+\nu I)^*]$$
$$= D_r^{(\alpha,\beta)}[(A+\gamma I), (B+\mu I)].$$

This completes the proof. $\square$

***Proof of Theorem 16 (Invariance under unitary transformations)***. The proof of this theorem is similar to that of the proof for Theorem 15, using the fact that $C^* = C^{-1}$ and the properties

$$\det_X[C(A+\gamma I)C^{-1}] = \det_X[(A+\gamma I)],$$
$$\text{tr}_X[C(A+\gamma I)C^{-1}] = \text{tr}_X[(A+\gamma I)],$$

of the operations $\det_X$ and $\text{tr}_X$. $\square$

***Proof of Theorem 17***. For the case $\alpha > 0$, $\beta > 0$, this follows immediately from Definition 1. For the case $\alpha > 0, \beta = 0$, by Definition 2 and Lemma 8, we have

$$D_r^{(\alpha,0)}[(A+\gamma I),(B+\mu I)] = \frac{r}{\alpha^2}\left[\left(\frac{\mu}{\gamma}\right)^r - 1\right]\log\left(\frac{\mu}{\gamma}\right)$$
$$+ \frac{1}{\alpha^2}\text{tr}_X(\Lambda + \frac{\gamma}{\mu})^{-r} - I) - \frac{1}{\alpha^2}\left(\frac{\mu}{\gamma}\right)^r \log\det_X(\Lambda + \frac{\gamma}{\mu})^{-r}$$
$$= D_r^{(\alpha,0)}[(\Lambda + \frac{\gamma}{\mu}), I].$$

The case $\alpha = 0, \beta > 0$ is entirely similar. $\square$

***Proof of Theorem 18***. We first note that $(\Lambda + \frac{\gamma}{\mu}I)^\omega = (\frac{\gamma}{\mu})^\omega(\frac{\mu}{\gamma}\Lambda + I)^\omega$. Then for $\alpha > 0, \beta > 0$, the statement of the theorem follows immediately from Definition 1. For the case $\alpha > 0, \beta = 0$, by Definition 2 and Lemma 8, we have

$$D_{\omega r}^{(\omega\alpha,0)}[(A+\gamma I),(B+\mu I)] = \frac{r}{\omega^2\alpha^2}\left[\left(\frac{\mu}{\gamma}\right)^{\omega r} - 1\right]\log\left(\frac{\mu}{\gamma}\right)^\omega$$
$$+ \frac{1}{\omega^2\alpha^2}\text{tr}_X(\Lambda + \frac{\gamma}{\mu})^{-\omega r} - I) - \frac{1}{\omega^2\alpha^2}\left(\frac{\mu}{\gamma}\right)^{\omega r}\log\det_X(\Lambda + \frac{\gamma}{\mu})^{-\omega r}$$
$$= \frac{1}{\omega^2}D_r^{(\alpha,0)}[(\Lambda + \frac{\gamma}{\mu})^\omega, I].$$

The case $\alpha = 0, \beta > 0$ is entirely similar. $\square$

*Appendix A.6. Proofs of Theorems 1, 2, and 3*

We are now ready to provide the proofs for Theorems 1, 2, and 3.

For the proof of positivity, we first need the following technical result.

**Lemma 13.** *(i) Let $r \neq 0$ be fixed. The function $f(x) = x^r - 1 - r\log(x)$ for $x > 0$ has a unique global minimum $f_{\min} = f(1) = 0$. In other words, $f(x) \geq 0\ \forall x > 0$, with equality if and only if $x = 1$.*

*(ii) Let $\nu > 0, r \neq 0$ be fixed. For $r \neq 0$, the function $g(x) = (\frac{x}{\nu} + 1)^r - 1 - r\log(\frac{x}{\nu} + 1)$ for $x > -\nu$ has a unique global minimum $g_{\min} = g(0) = 0$. In other words, $g(x) \geq 0\ \forall x > -\nu$, with equality if and only if $x = 0$.*

*Proof of Lemma 13.* (i) We have $f'(x) = \frac{r(x^r - 1)}{x}$. When $r > 0$, we have $x^r < 1$ for $0 < x < 1$ and $x^r > 1$ for $x > 1$. When $r < 0$, we have $x^r > 1$ for $0 < x < 1$ and $x^r < 1$ for $x > 1$. Thus, for all $r \neq 0$, we have $f'(x) < 0$ when $0 < x < 1$ and $f'(x) > 0$ when $x > 1$. Hence $f$ has a unique global minimum $f_{\min} = f(1) = 0$.

(ii) The proof for $g$ follows that for $f$ by the change of variable $y = \frac{x}{\nu} + 1$. □

*Proof of Theorem 1 (Positivity).* For the case $\alpha > 0, \beta > 0$, this is a special case of Theorem 6, with $p + q = r$. Consider now the case $\alpha = 0, \beta > 0$ (the case $\alpha > 0, \beta = 0$ then follows by dual symmetry). For the proof of positivity, we can ignore the positive factor $\beta^2$ and thus it suffices to consider $D_r^{(0,1)}$. We recall that we define $\Lambda + \nu I = (B + \mu I)^{-1/2}(A + \gamma I)(B + \mu I)^{-1/2}$, where $\nu = \frac{\gamma}{\mu}$. Then, since $\det_X[(B + \mu I)^{-1/2}(A + \gamma I)(B + \mu I)^{-1/2}] = \det_X[(B + \mu I)^{-1}(A + \gamma I)]$ and $\text{tr}_X[(B + \mu I)^{-1/2}(A + \gamma I)(B + \mu I)^{-1/2}] = \text{tr}_X[(B + \mu I)^{-1}(A + \gamma I)]$, we have

$$D_r^{(0,1)}[(A + \gamma I), (B + \mu I)]$$
$$= r(\nu^r - 1)\log\nu + \text{tr}_X[(\Lambda + \nu I)^r - I] - \nu^r \log\det_X(\Lambda + \nu I)^r$$

By Lemma 7,

$$\text{tr}_X[(\Lambda + \nu I)^r - I] = \nu^r - 1 + \nu^r \text{tr}\left[\left(\frac{\Lambda}{\nu} + I\right)^r - I\right].$$

Also

$$\log\det_X(\Lambda + \nu I)^r = \log\left[\nu^r \det\left(\frac{\Lambda}{\nu} + I\right)^r\right] = r\log\det\left(\frac{\Lambda}{\nu} + I\right) + r\log\nu.$$

Thus we have

$$D_r^{(0,1)}[(A+\gamma I),(B+\mu I)]$$
$$= \nu^r - 1 - r\log\nu + \nu^r \left(\mathrm{tr}\left[\left(\frac{\Lambda}{\nu}+I\right)^r - I\right] - r\log\det\left(\frac{\Lambda}{\nu}+I\right)\right)$$
$$= \nu^r - 1 - r\log\nu + \nu^r \left[\sum_{k=1}^{\infty}\left(\frac{\lambda_k}{\nu}+1\right)^r - 1 - r\log\left(\frac{\lambda_k}{\nu}+1\right)\right].$$

By the first part of Lemma 13, we have for all $\nu > 0$

$$\nu^r - 1 - r\log\nu \geq 0,$$

with equality if and only if $\nu = 1$. By the second part of the Lemma 13, we have for all $k \in \mathbb{N}$

$$\left(\frac{\lambda_k}{\nu}+1\right)^r - 1 - r\log\left(\frac{\lambda_k}{\nu}+1\right) \geq 0,$$

with equality if and only $\lambda_k = 0$. Combining these two inequalities, we obtain

$$D_r^{(0,1)}[(A+\gamma I),(B+\mu I)] \geq 0,$$

with equality if and only if $\nu = \frac{\gamma}{\mu} = 1$ and $\lambda_k = 0 \forall k \in \mathbb{N} \iff \Lambda = I$, that is if and only $(B+\mu I)^{-1/2}(A+\gamma I)(B+\mu I)^{-1/2} = I \iff A+\gamma I = B+\mu I \iff A = B$ and $\gamma = \mu$. This completes the proof. $\square$

***Proof of Theorem 2 (Special cases - I)***. The first statement of the theorem is the content of Theorem 9. The second statement is the content of Theorem 10. $\square$

***Proof of Theorem 3 (Special cases - II)***. This theorem follows from Theorems 9 and 10 as well as the symmetry of $D_r^{(\alpha,\alpha)}[(A+\gamma I),(B+\mu I)]$ as proved in Theorem 13. $\square$

*Appendix A.7. Proofs for the divergences between RKHS covariance operators*

In this section, we prove Theorems 23, 24, 25, and 26. We first need the following preliminary results.

**Lemma 14.** *Let $\mathcal{H}_1, \mathcal{H}_2$ be separable Hilbert spaces. Let $A : \mathcal{H}_1 \to \mathcal{H}_2$ and $B : \mathcal{H}_2 \to \mathcal{H}_1$ be compact linear operators such that both $AB : \mathcal{H}_2 \to \mathcal{H}_2$ and $BA : \mathcal{H}_1 \to \mathcal{H}_1$ are trace class operators. Let $\alpha, \beta > 0$ be fixed. For any $p, q \in \mathbb{R}$,*

$$\det\left[\frac{\alpha(AB + I_{\mathcal{H}_2})^p + \beta(AB + I_{\mathcal{H}_2})^q}{\alpha + \beta}\right]$$
$$= \det\left[\frac{\alpha(BA + I_{\mathcal{H}_1})^p + \beta(BA + I_{\mathcal{H}_1})^q}{\alpha + \beta}\right]. \tag{A.23}$$

*Proof of Lemma 14.* Since the nonzero eigenvalues of $AB : \mathcal{H}_2 \to \mathcal{H}_2$ and $BA : \mathcal{H}_1 \to \mathcal{H}_1$ are the same, we have for any $p \in \mathbb{R}$

$$\det[(AB + I_{\mathcal{H}_2})^p] = \det[(BA + I_{\mathcal{H}_1})^p].$$

For any $p, q \in \mathbb{R}$,

$$\det\left[\frac{\alpha(AB + I_{\mathcal{H}_2})^p + \beta(AB + I_{\mathcal{H}_2})^q}{\alpha + \beta}\right]$$
$$= \det\left[\frac{\alpha(BA + I_{\mathcal{H}_1})^p + \beta(BA + I_{\mathcal{H}_1})^q}{\alpha + \beta}\right].$$

In the above equality, we have used the fact that a zero eigenvalue of $AB$ and $BA$ corresponds to an eigenvalue equal to 1 for $\frac{\alpha(AB+I_{\mathcal{H}_2})^p+\beta(AB+I_{\mathcal{H}_2})^q}{\alpha+\beta} : \mathcal{H}_2 \to \mathcal{H}_2$ and $\frac{\alpha(BA+I_{\mathcal{H}_1})^p+\beta(BA+I_{\mathcal{H}_1})^q}{\alpha+\beta} : \mathcal{H}_1 \to \mathcal{H}_1$, respectively, which does not change the determinant. This completes the proof. $\square$

**Lemma 15.** *Let $\mathcal{H}_1, \mathcal{H}_2$ be separable Hilbert spaces. Let $A, B : \mathcal{H}_1 \to \mathcal{H}_2$ be compact linear operators such that both $AA^* : \mathcal{H}_2 \to \mathcal{H}_2$ and $BB^* : \mathcal{H}_2 \to \mathcal{H}_2$ are trace class operators. Let $\alpha, \beta > 0$ be fixed. For any $p, q \in \mathbb{R}$,*

$$\det\left[\frac{\alpha[(AA^* + I_{\mathcal{H}_2})(BB^* + I_{\mathcal{H}_2})^{-1}]^p + \beta[(AA^* + I_{\mathcal{H}_2})(BB^* + I_{\mathcal{H}_2})^{-1}]^q}{\alpha + \beta}\right]$$
$$= \det\left[\frac{\alpha(C + I_{\mathcal{H}_1} \otimes I_3)^p + \beta(C + I_{\mathcal{H}_1} \otimes I_3)^q}{\alpha + \beta}\right], \tag{A.24}$$

*where*

$$C = \begin{pmatrix} A^*A & -A^*B(I_{\mathcal{H}_1} + B^*B)^{-1} & -A^*AA^*B(I_{\mathcal{H}_1} + B^*B)^{-1} \\ B^*A & -B^*B(I_{\mathcal{H}_1} + B^*B)^{-1} & -B^*AA^*B(I_{\mathcal{H}_1} + B^*B)^{-1} \\ B^*A & -B^*B(I_{\mathcal{H}_1} + B^*B)^{-1} & -B^*AA^*B(I_{\mathcal{H}_1} + B^*B)^{-1} \end{pmatrix}. \tag{A.25}$$

***Proof of Lemma 15.*** We make use of the following notation. Let $A, B, C : \mathcal{H}_1 \to \mathcal{H}_2$ be three bounded linear operators. Consider the operator $(A\ B\ C) : \mathcal{H}_1^3 \to \mathcal{H}_2$, with

$$(A\ B\ C)^* = \begin{pmatrix} A^* \\ B^* \\ C^* \end{pmatrix} : \mathcal{H}_2 \to \mathcal{H}_1^3.$$

Here $\mathcal{H}_1^3 = \mathcal{H}_1 \oplus \mathcal{H}_1 \oplus \mathcal{H}_1$ denotes the direct sum of $\mathcal{H}_1$ with itself, that is

$$\mathcal{H}_1^3 = \mathcal{H}_1 \oplus \mathcal{H}_1 \oplus \mathcal{H}_1 = \{(v_1, v_2, v_3) \ : \ v_1, v_2, v_3 \in \mathcal{H}_1\},$$

equipped with the inner product

$$\langle (v_1, v_2, v_3), (w_1, w_2, w_3) \rangle_{\mathcal{H}_1^3} = \langle v_1, w_1 \rangle_{\mathcal{H}_1} + \langle v_2, w_2 \rangle_{\mathcal{H}_1} + \langle v_3, w_3 \rangle_{\mathcal{H}_1}.$$

If $\{e_i\}_{i=1}^\infty$ is an orthonormal basis for $\mathcal{H}_1$, then $\{(e_i, 0, 0)\}_{i=1}^\infty \cup \{(0, e_i, 0)\}_{i=1}^\infty \cup \{(0, 0, e_i)\}_{i=1}^\infty$ is an orthonormal basis for $\mathcal{H}_1^3$.

We now utilize this notation in our setting. By the Sherman-Morrison-Woodbury formula, we have

$$(BB^* + I_{\mathcal{H}_2})^{-1} = I_{\mathcal{H}_2} - B(I_{\mathcal{H}_1} + B^*B)^{-1} B^*.$$

Thus it follows that

$$(AA^* + I_{\mathcal{H}_2})(BB^* + I_{\mathcal{H}_2})^{-1} = I_{\mathcal{H}_2} + AA^* - B(I_{\mathcal{H}_1} + B^*B)^{-1} B^*$$
$$- AA^* B(I_{\mathcal{H}_1} + B^*B)^{-1} B^*$$
$$= I_{\mathcal{H}_2} + C_1 C_2.$$

Here the operators $C_1, C_2$ are defined as follows.

$$C_1 = [A\ \ -B(I_{\mathcal{H}_1} + B^*B)^{-1}\ \ -AA^*B(I_{\mathcal{H}_1} + B^*B)^{-1}] : \mathcal{H}_1^3 \to \mathcal{H}_2,$$

$$C_2 = \begin{pmatrix} A^* \\ B^* \\ B^* \end{pmatrix} : \mathcal{H}_2 \to \mathcal{H}_1^3.$$

The operator $C_2 C_1 : \mathcal{H}_1^3 \to \mathcal{H}_1^3$ is given by

$$C_2 C_1 = \begin{pmatrix} A^*A & -A^*B(I_{\mathcal{H}_1} + B^*B)^{-1} & -A^*AA^*B(I_{\mathcal{H}_1} + B^*B)^{-1} \\ B^*A & -B^*B(I_{\mathcal{H}_1} + B^*B)^{-1} & -B^*AA^*B(I_{\mathcal{H}_1} + B^*B)^{-1} \\ B^*A & -B^*B(I_{\mathcal{H}_1} + B^*B)^{-1} & -B^*AA^*B(I_{\mathcal{H}_1} + B^*B)^{-1} \end{pmatrix}.$$

It follows from Lemma 14 that

$$\det\left[\frac{\alpha[(AA^* + I_{\mathcal{H}_2})(BB^* + I_{\mathcal{H}_2})^{-1}]^p + \beta[(AA^* + I_{\mathcal{H}_2})(BB^* + I_{\mathcal{H}_2})^{-1}]^q}{\alpha + \beta}\right]$$

$$= \det\left[\frac{\alpha(I_{\mathcal{H}_2} + C_1 C_2)^p + \beta(I_{\mathcal{H}_2} + C_1 C_2)^q}{\alpha + \beta}\right]$$

$$= \det\left[\frac{\alpha(C_2 C_1 + I_{\mathcal{H}_1} \otimes I_3)^p + \beta(C_2 C_1 + I_{\mathcal{H}_1} \otimes I_3)^q}{\alpha + \beta}\right].$$

This completes the proof. □

***Proof of Theorem 23.*** Let $\Lambda + \frac{\gamma}{\mu}I = (BB^* + \mu I_{\mathcal{H}_2})^{-1/2}(AA^* + \gamma I)(BB^* + \mu I)^{-1/2}$ and $Z + \frac{\gamma}{\mu}I = (AA^* + \gamma I)(BB^* + \mu I)^{-1}$, with $\frac{\mu}{\gamma}Z + I = (\frac{AA^*}{\gamma} + I)(\frac{BB^*}{\mu} + I)^{-1}$.

By Theorem 8, we have

$$D_r^{(\alpha,\beta)}[(AA^* + \gamma I_{\mathcal{H}_2}), (BB^* + \mu I_{\mathcal{H}_2})]$$

$$= \frac{r(\delta - \frac{\alpha}{\alpha+\beta})}{\alpha\beta}\left(\log\frac{\gamma}{\mu}\right) + \frac{1}{\alpha\beta}\log\left(\frac{\alpha(\frac{\gamma}{\mu})^p + \beta(\frac{\gamma}{\mu})^{-q}}{\alpha + \beta}\right)$$

$$+ \frac{1}{\alpha\beta}\log\det\left[\frac{\alpha(\Lambda + \frac{\gamma}{\mu}I)^p + \beta(\Lambda + \frac{\gamma}{\mu}I)^{-q}}{\alpha(\frac{\gamma}{\mu})^p + \beta(\frac{\gamma}{\mu})^{-q}}\right],$$

with $p = r(1 - \delta)$ and $q = r\delta$. The determinant in the last term is

$$\det\left[\frac{\alpha(\Lambda + \frac{\gamma}{\mu}I)^p + \beta(\Lambda + \frac{\gamma}{\mu}I)^{-q}}{\alpha(\frac{\gamma}{\mu})^p + \beta(\frac{\gamma}{\mu})^{-q}}\right] = \det\left[\frac{\alpha(\frac{\gamma}{\mu})^p(\frac{\mu}{\gamma}\Lambda + I)^p + \beta(\frac{\gamma}{\mu})^{-q}(\frac{\mu}{\gamma}\Lambda + I)^{-q}}{\alpha(\frac{\gamma}{\mu})^p + \beta(\frac{\gamma}{\mu})^{-q}}\right]$$

$$= \det\left[\frac{\alpha(\frac{\gamma}{\mu})^p(\frac{\mu}{\gamma}Z + I)^p + \beta(\frac{\gamma}{\mu})^{-q}(\frac{\mu}{\gamma}Z + I)^{-q}}{\alpha(\frac{\gamma}{\mu})^p + \beta(\frac{\gamma}{\mu})^{-q}}\right]$$

$$= \det\left[\frac{\alpha(\frac{\gamma}{\mu})^p(C + I_{\mathcal{H}_1} \otimes I_3)^p + \beta(\frac{\gamma}{\mu})^{-q}(C + I_{\mathcal{H}_1} \otimes I_3)^{-q}}{\alpha(\frac{\gamma}{\mu})^p + \beta(\frac{\gamma}{\mu})^{-q}}\right]$$

by Lemma 15, where

$$C = \begin{pmatrix} \frac{A^*A}{\gamma} & -\frac{A^*B}{\sqrt{\gamma\mu}}(I_{\mathcal{H}_1} + \frac{B^*B}{\mu})^{-1} & -\frac{A^*AA^*B}{\gamma\sqrt{\gamma\mu}}(I_{\mathcal{H}_1} + \frac{B^*B}{\mu})^{-1} \\ \frac{B^*A}{\sqrt{\gamma\mu}} & -\frac{B^*B}{\mu}(I_{\mathcal{H}_1} + \frac{B^*B}{\mu})^{-1} & -\frac{B^*AA^*B}{\gamma\mu}(I_{\mathcal{H}_1} + \frac{B^*B}{\mu})^{-1} \\ \frac{B^*A}{\sqrt{\gamma\mu}} & -\frac{B^*B}{\mu}(I_{\mathcal{H}_1} + \frac{B^*B}{\mu})^{-1} & -\frac{B^*AA^*B}{\gamma\mu}(I_{\mathcal{H}_1} + \frac{B^*B}{\mu})^{-1} \end{pmatrix},$$

which is obtained by replacing $AA^*$ and $BB^*$ in Lemma 15 with $\frac{AA^*}{\gamma}$ and $\frac{BB^*}{\mu}$, respectively. This completes the proof of the theorem. □

***Proof of Theorem 24.*** Let $Z + \frac{\gamma}{\mu}I = (AA^* + \gamma I)(BB^* + \mu I)^{-1}$. By the finite-dimensional formula given in Eq. (19), we have

$$D_r^{(\alpha,\beta)}[(AA^* + \gamma I_{\mathcal{H}_2}), (BB^* + \mu I_{\mathcal{H}_2})] = \frac{1}{\alpha\beta}\log\det\left[\frac{\alpha(Z + \frac{\gamma}{\mu}I)^p + \beta(Z + \frac{\gamma}{\mu}I)^{-q}}{\alpha + \beta}\right]$$

$$= \frac{1}{\alpha\beta}\left[\log\left(\frac{\alpha(\frac{\gamma}{\mu})^p + \beta(\frac{\gamma}{\mu})^{-q}}{\alpha + \beta}\right)\right]\dim(\mathcal{H}_2)$$

$$+ \frac{1}{\alpha\beta}\log\det\left[\frac{\alpha(Z + \frac{\gamma}{\mu}I)^p + \beta(Z + \frac{\gamma}{\mu}I)^{-q}}{\alpha(\frac{\gamma}{\mu})^p + \beta(\frac{\gamma}{\mu})^{-q}}\right].$$

As in the proof of Theorem 23, the determinant in last term in the above expression is

$$\det\left[\frac{\alpha(\Lambda + \frac{\gamma}{\mu}I)^p + \beta(\Lambda + \frac{\gamma}{\mu}I)^{-q}}{\alpha(\frac{\gamma}{\mu})^p + \beta(\frac{\gamma}{\mu})^{-q}}\right]$$

$$= \det\left[\frac{\alpha(\frac{\gamma}{\mu})^p(C + I_{\mathcal{H}_1} \otimes I_3)^p + \beta(\frac{\gamma}{\mu})^{-q}(C + I_{\mathcal{H}_1} \otimes I_3)^{-q}}{\alpha(\frac{\gamma}{\mu})^p + \beta(\frac{\gamma}{\mu})^{-q}}\right].$$

This gives us the final expression. □

***Proof of Theorem 25.*** We consider the linear operators

$$A = \frac{1}{\sqrt{m}}\Phi(\mathbf{x})J_m : \mathbb{R}^m \to \mathcal{H}_K, \quad B = \frac{1}{\sqrt{m}}\Phi(\mathbf{y})J_m : \mathbb{R}^m \to \mathcal{H}_K.$$

The desired expression then follows from Theorem 23. □

***Proof of Theorem 26.*** This is proved in the same way as Theorem 25, except that we invoke Theorem 24. □

*Appendix A.8. Proofs for the metric properties*

In this section, we prove Theorems 19, 20, 21, which lead to the proofs of Theorems 4 and 22. We present two sets of separate proofs for Theorems 4 and 22, one simpler proof for the particular case $\alpha = 1/2$, which corresponds to the infinite-dimensional symmetric Stein divergence, and one general proof for any $\alpha > 0$. The former case utilizes Theorem 28 and the latter case utilizes Theorem 30, both of which should be of interest in their own right.

*Appendix A.8.1. The case of the infinite-dimensional symmetric Stein divergence*

Consider the first case $\alpha = 1/2$, which corresponds to the infinite-dimensional symmetric Stein divergence.

**Lemma 16.** *Let $\mathcal{H}$ be a separable Hilbert space. Let $A, B, C : \mathcal{H} \to \mathcal{H}$ be self-adjoint finite-rank operators, such that $A + I > 0$, $B + I > 0$, $C + I > 0$. Then*

$$\sqrt{\log \frac{\det(\frac{A+B}{2} + I)}{\sqrt{\det(A+I)\det(B+I)}}} \leq \sqrt{\log \frac{\det(\frac{A+C}{2} + I)}{\sqrt{\det(A+I)\det(C+I)}}}$$
$$+ \sqrt{\log \frac{\det(\frac{C+B}{2} + I)}{\sqrt{\det(C+I)\det(B+I)}}}. \quad (A.26)$$

***Proof of Lemma 16.*** Since $A, B, C$ are all finite-rank operators, there exists a finite-dimensional subspace $\mathcal{H}_n \subset \mathcal{H}$, with $\dim(\mathcal{H}_n) = n$ for some $n \in \mathbb{N}$, such that $\text{range}(A) \subset \mathcal{H}_n$, $\text{range}(B) \subset \mathcal{H}_n$, and $\text{range}(C) \subset \mathcal{H}_n$. Let

$$A_n = A\big|_{\mathcal{H}_n} : \mathcal{H}_n \to \mathcal{H}_n, \quad B_n = B\big|_{\mathcal{H}_n} : \mathcal{H}_n \to \mathcal{H}_n, \quad C_n = C\big|_{\mathcal{H}_n} : \mathcal{H}_n \to \mathcal{H}_n.$$

Then $A_n, B_n, C_n$ are linear operators on the finite-dimensional space $\mathcal{H}_n$ and thus are represented by $n \times n$ matrices, which we denote by the same symbols. We also have

$$(A+B)_n = (A+B)\big|_{\mathcal{H}_n} = A\big|_{\mathcal{H}_n} + B\big|_{\mathcal{H}_n} = A_n + B_n,$$
$$(A+C)_n = A_n + C_n, \quad (C+B)_n = B_n + C_n.$$

Applying the finite-dimensional result in [16], we then obtain

$$\sqrt{\log \frac{\det(\frac{A_n+B_n}{2} + I_n)}{\sqrt{\det(A_n+I_n)\det(B_n+I_n)}}} \leq \sqrt{\log \frac{\det(\frac{A_n+C_n}{2} + I_n)}{\sqrt{\det(A_n+I_n)\det(C_n+I_n)}}}$$
$$+ \sqrt{\log \frac{\det(\frac{C_n+B_n}{2} + I_n)}{\sqrt{\det(C_n+I_n)\det(B_n+I_n)}}}.$$

It is clear that the non-zero eigenvalues of $A$ and $A_n$ are the same, so that $\det(A+I) = \det(A_n + I_n)$ and the same holds true for the other operators. This gives us the final result. $\square$

***Proof of Theorem 21 (Triangle inequality- square root of symmetric Stein divergence).*** Let $\{A_n\}_{n \in \mathbb{N}}, \{B_n\}_{n \in \mathbb{N}}, \{C_n\}_{n \in \mathbb{N}}$ be sequences of finite-rank operators with

$$||A_n - A||_{\text{tr}} \to 0, \quad ||B_n - B||_{\text{tr}} \to 0, \quad ||C_n - C||_{\text{tr}} \to 0, \quad \text{as } n \to \infty.$$

By Lemma 16, we have

$$\sqrt{\log \frac{\det(\frac{A_n+B_n}{2}+I)}{\sqrt{\det(A_n+I)\det(B_n+I)}}} \leq \sqrt{\log \frac{\det(\frac{A_n+C_n}{2}+I)}{\sqrt{\det(A_n+I)\det(C_n+I)}}}$$
$$+ \sqrt{\log \frac{\det(\frac{C_n+B_n}{2}+I)}{\sqrt{\det(C_n+I)\det(B_n+I)}}}.$$

By Theorem 3.5 in [22], as $n \to \infty$, we have

$$\det(A_n+I) \to \det(A+I), \quad \det(B_n+I) \to \det(B+I),$$
$$\det(\frac{A_n+B_n}{2}+I) \to \det(\frac{A+B}{2}+I),$$

and the same holds true for the other operators. Thus by taking the limit as $n \to \infty$ in the above triangle inequality for $(A_n+I), (B_n+I)$ and $(C_n+I)$, we obtain the final triangle inequality for $(A+I), (B+I)$, and $(C+I)$. $\square$

The following is the specialization of Theorem 4 when $\alpha = 1/2$.

**Theorem 27** (**Metric property - square root of symmetric Stein divergence**). *Let $\gamma > 0, \gamma \in \mathbb{R}$ be fixed. The square root of the infinite-dimensional symmetric Stein divergence $\sqrt{D_1^{(1/2,1/2)}[(A+\gamma I),(B+\gamma I)]}$ is a metric on* $\mathrm{PTr}(\mathcal{H})(\gamma)$.

***Proof of Theorem 27***. We have already shown the positivity and symmetry of $D_1^{(1/2,1/2)}[(A+\gamma I),(B+\gamma I)]$. It remains for us to show the triangle inequality, namely

$$\sqrt{D_1^{(1/2,1/2)}[(A+\gamma I),(B+\gamma I)]} \leq \sqrt{D_1^{(1/2,1/2)}[(A+\gamma I),(C+\gamma I)]}$$
$$+ \sqrt{D_1^{(1/2,1/2)}[(C+\gamma I),(B+\gamma I)]},$$

for any three operators $(A+\gamma I), (B+\gamma I), (C+\gamma I) \in \mathrm{PTr}(\mathcal{H})$. We have

$$D_1^{(1/2,1/2)}[(A+\gamma I),(B+\gamma I)] = 4\log\left[\frac{\det_{\mathrm{X}}(\frac{(A+\gamma I)+(B+\gamma I)}{2})}{\det_{\mathrm{X}}(A+\gamma I)^{1/2}\det_{\mathrm{X}}(B+\gamma I)^{1/2}}\right]$$
$$= 4\log\left[\frac{\det(\frac{A+B}{2\gamma}+I)}{\det(\frac{A}{\gamma}+I)^{1/2}\det(\frac{B}{\gamma}+I)^{1/2}}\right].$$

Thus the triangle inequality for $D_1^{(1/2,1/2)}[(A+\gamma I),(B+\gamma I)]$ follows that stated in Theorem 21. $\square$

**Lemma 17.** *Let $\mathcal{H}$ be a separable Hilbert space. Let $A, B : \mathcal{H} \to \mathcal{H}$ be self-adjoint finite-rank operators, with maximum rank $n$, $n \in \mathbb{N}$, such that $A + I > 0$, $B + I > 0$. Then*

$$\prod_{j=1}^{n} \left[ \frac{\lambda_j(A) + \lambda_j(B)}{2} + 1 \right] \leq \det\left( \frac{A+B}{2} + I \right). \tag{A.27}$$

*Proof of Lemma 17.* Since $A, B$ are both finite-rank operators, there exists a finite-dimensional subspace $\mathcal{H}_n \subset \mathcal{H}$, with $\dim(\mathcal{H}_n) = n$, such that $\text{range}(A) \subset \mathcal{H}_n, \text{range}(B) \subset \mathcal{H}_n$. Let

$$A_n = A\big|_{\mathcal{H}_n} : \mathcal{H}_n \to \mathcal{H}_n, \quad B_n = B\big|_{\mathcal{H}_n} : \mathcal{H}_n \to \mathcal{H}_n.$$

Then $A_n, B_n$ are linear operators on the finite-dimensional space $\mathcal{H}_n$ and thus are represented by $n \times n$ matrices, which we denote by the same symbols. We also have

$$(A+B)_n = (A+B)\big|_{\mathcal{H}_n} = A\big|_{\mathcal{H}_n} + B\big|_{\mathcal{H}_n} = A_n + B_n.$$

Thus we can apply the following inequality for finite-dimensional SPD matrices ([23])

$$\prod_{j=1}^{n} \left[ \frac{\lambda_j(A_n) + \lambda_j(B_n)}{2} + 1 \right] = \prod_{j=1}^{n} \left[ \frac{\lambda_j(A_n + I_n) + \lambda_j(B_n + I_n)}{2} \right]$$

$$\leq \det\left( \frac{A_n + B_n}{2} + I_n \right).$$

We note that the non-zero eigenvalues of $A_n, B_n$ are the same as those of $A, B$, respectively, with the maximum number being $n$, and $\det(\frac{A+B}{2} + I) = \det(\frac{A_n+B_n}{2} + I_n)$. Together with the previous inequality, this gives us the final result. $\square$

**Theorem 28.** *Let $\mathcal{H}$ be a separable Hilbert space. Let $A, B : \mathcal{H} \to \mathcal{H}$ be self-adjoint trace class operators, such that $A + I > 0$, $B + I > 0$. Then*

$$\prod_{j=1}^{\infty} \left[ \frac{\lambda_j(A) + \lambda_j(B)}{2} + 1 \right] \leq \det\left( \frac{A+B}{2} + I \right). \tag{A.28}$$

*Proof of Theorem 28.* Let $A = \sum_{j=1}^{\infty} \lambda_j(A) \phi_j \otimes \phi_j$ denote the spectral decomposition for $A$. For each $n \in \mathbb{N}$, define

$$A_n = \sum_{j=1}^{n} \lambda_j(A) \phi_j \otimes \phi_j.$$

Then $A_n$ is a finite-rank operator with the eigenvalues being the first $n$ eigenvalues of $A$ and $\lim_{n\to\infty} ||A_n - A||_{\text{tr}} = 0$. In the same way, we construct a sequence of finite-rank operators $B_n$ with $\lim_{n\to\infty} ||B_n - B||_{\text{tr}} = 0$, so that

$$\lim_{n\to\infty} ||(A_n + B_n) - (A + B)||_{\text{tr}} = 0.$$

By Theorem 3.5 in [22], as $n \to \infty$, we then have

$$\lim_{n\to\infty} \det\left(\frac{A_n + B_n}{2} + I\right) = \det\left(\frac{A + B}{2} + I\right).$$

Applying Lemma 17 to $A_n$ and $B_n$, we have

$$\prod_{j=1}^{n}\left[\frac{\lambda_j(A_n) + \lambda_j(B_n)}{2} + 1\right] \leq \det\left(\frac{A_n + B_n}{2} + I\right). \quad (A.29)$$

The final result is then obtained by taking the limit as $n \to \infty$, noting that the eigenvalues of $A_n, B_n$, are precisely the first $n$ eigenvalues of $A, B$, respectively. $\square$

The following is the specialization of Theorem 22 when $\alpha = 1/2$.

**Theorem 29.** *Let $\mathcal{H}$ be a separable Hilbert space. Let $A, B : \mathcal{H} \to \mathcal{H}$ be self-adjoint trace class operators, such that $A + I > 0$, $B + I > 0$. Let $\text{Eig}(A), \text{Eig}(B) : \ell^2 \to \ell^2$ be diagonal operators with the diagonals consisting of the eigenvalues of $A$ and $B$, respectively, in decreasing order. Then*

$$D_1^{(1/2,1/2)}[(\text{Eig}(A) + I), (\text{Eig}(B) + I)] \leq D_1^{(1/2,1/2)}[(A + I), (B + I)]. \quad (A.30)$$

*Proof of Theorem 29.* By definition, we have

$$D_1^{(1/2,1/2)}[(\text{Eig}(A) + I), (\text{Eig}(B) + I)]$$
$$= 4\log\left[\frac{\det(\frac{\text{Eig}(A)+\text{Eig}(B)}{2} + I)}{\sqrt{\det(\text{Eig}(A) + I)\det(\text{Eig}(B) + I)}}\right]$$
$$= 4\log\left[\frac{\prod_{j=1}^{\infty}\left[\frac{\lambda_j(A)+\lambda_j(B)}{2} + 1\right]}{\sqrt{\det(A + I)\det(B + I)}}\right] \leq 4\log\left[\frac{\det(\frac{A+B}{2} + I)}{\sqrt{\det(A + I)\det(B + I)}}\right]$$

by Theorem 28
$$= D_1^{(1/2,1/2)}[(A + I), (B + I)].$$

This completes the proof. $\square$

*Appendix A.8.2. The general case*

We now consider the general case $\alpha > 0$. We need the following results.

In the following, let $\mathscr{C}_p(\mathcal{H})$ denote the class of $p$th Schatten class operators on $\mathcal{H}$, under the norm $|| \ ||_p$, $1 \leq p \leq \infty$, which is defined by

$$||A||_p = [\sum_{k=1}^{\infty} \lambda_k^p(A^*A)^{1/2})]^{1/p}, \tag{A.31}$$

with $\mathscr{C}_1(\mathcal{H})$ being the space of trace class operators $\text{Tr}(\mathcal{H})$, $\mathscr{C}_2(\mathcal{H})$ being the space of Hilbert-Schmidt operators $\text{HS}(\mathcal{H})$, and $\mathscr{C}_\infty(\mathcal{H})$ being the set of compact operators under the operator norm $|| \ ||$.

**Theorem 30.** *Let $r \in \mathbb{R}$ be fixed but arbitrary. Assume that $1 \leq p \leq \infty$. Let $\{A_n\}_{n \in \mathbb{N}} \in \text{Sym}(\mathcal{H}) \cap \mathscr{C}_p(\mathcal{H})$, $A \in \text{Sym}(\mathcal{H}) \cap \mathscr{C}_p(\mathcal{H})$ be such that $I + A > 0$, $I + A_n > 0 \ \forall n \in \mathbb{N}$. Assume that $\lim_{n \to \infty} ||A_n - A||_p = 0$. Then*

$$\lim_{n \to \infty} ||(I + A_n)^r - (I + A)^r||_p = 0. \tag{A.32}$$

***Proof of Theorem 30.*** (i) We first prove that

$$\lim_{n \to \infty} ||(I + A_n)^r - (I + A)^r||_p = 0, \ \ 0 \leq r \leq 1. \tag{A.33}$$

The case $r = 0$ is trivial. Let us prove this for $0 < r \leq 1$. For this limit, we make use of the following result from [24] (Corollary 3.2), which states that for any two positive operators $A, B$ on $\mathcal{H}$ such that $A \geq c > 0$, $B \geq c > 0$, and any operator $X$ on $\mathcal{H}$,

$$||A^r X - X B^r||_p \leq rc^{r-1}||AX - XB||_p, \tag{A.34}$$

where $0 < r \leq 1$ and $|| \ ||_p$, $1 \leq p \leq \infty$, denotes the Schatten $p$-norm.

By the assumption $I + A > 0$, there exists $M_A > 0$ such that

$$\langle x, (I + A)x \rangle \geq M_A ||x||^2 \ \ \forall x \in \mathcal{H}.$$

By the assumption $\lim_{n \to \infty} ||A_n - A||_p = 0$, for any $\epsilon$ satisfying $0 < \epsilon < M_A$, there exists $N = N(\epsilon) \in \mathbb{N}$ such that $||A_n - A||_p < \epsilon \ \forall n \geq N$. Then $\forall x \in \mathcal{H}$,

$$|\langle x, (A_n - A)x \rangle| \leq ||A_n - A|| \ ||x||^2 \leq ||A_n - A||_p ||x||^2 \leq \epsilon ||x||^2.$$

It thus follows that $\forall x \in \mathcal{H}$,

$$\langle x, (I + A_n)x \rangle = \langle x, (I + A)x \rangle + \langle x, (A_n - A)x \rangle \geq (M_A - \epsilon)||x||^2.$$

Thus we have $I + A \geq M_A > 0$, $I + A_n \geq M_A - \epsilon > 0$ $\forall n \geq N = N(\epsilon)$. Then, applying Eq. (A.34), we have for all $n \geq N$,

$$||(I + A_n)^r - (I + A)^r||_p \leq r(M_A - \epsilon)^{r-1}||(I + A_n) - (I + A)||_p$$
$$= r\left(\frac{1}{M_A - \epsilon}\right)^{1-r} ||A_n - A||_p,$$

which implies

$$\lim_{n \to \infty} ||(I + A_n)^r - (I + A)^r||_p = 0.$$

This completes the proof of the first limit.

(ii) For $r > 1$, we proceed by induction as follows. We have

$$||(I + A_n)^r - (I + A)^r||_p$$
$$\leq ||(I + A_n)^r - (I + A_n)(I + A)^{r-1}||_p + ||(I + A_n)(I + A)^{r-1} - (I + A)^r||_p$$
$$\leq ||I + A_n|| \, ||(I + A_n)^{r-1} - (I + A)^{r-1}||_p + ||A_n - A||_p ||(I + A)^{r-1}||.$$

Thus this case follows from the case $0 \leq r \leq 1$ by induction.

(iii) We now prove that

$$\lim_{n \to \infty} ||(I + A_n)^{-1} - (I + A)^{-1}||_p = 0. \tag{A.35}$$

We have $\forall n \geq N = N(\epsilon)$,

$$||(I + A_n)^{-1} - (I + A)^{-1}||_p = ||(I + A_n)^{-1}[(I + A_n) - (I + A)](I + A)^{-1}||_p$$
$$\leq ||(I + A_n)^{-1}|| \, ||A_n - A||_p ||(I + A)^{-1}|| \leq \frac{1}{M_A(M_A - \epsilon)}||A_n - A||_p,$$

which implies that

$$\lim_{n \to \infty} ||(I + A_n)^{-1} - (I + A)^{-1}||_p = 0.$$

(iii) We next prove that

$$\lim_{n \to \infty} ||(I + A_n)^{-r} - (I + A)^{-r}||_p = 0, \quad 0 < r \leq 1. \tag{A.36}$$

We have
$$(I+A)^{-1} \geq \frac{1}{\max\{(1+\lambda_k(A)) : k \in \mathbb{N}\}} = \frac{1}{||I+A||} > 0.$$
From the limit $\lim_{n\to\infty} ||A_n - A|| = 0$, it follows that for any $\epsilon$ satisfying $0 < \epsilon < ||I+A||$, there exists $M = M(\epsilon) \in \mathbb{N}$ such that $\forall n \geq M$,
$$||I+A|| - \epsilon \leq ||I + A_n|| \leq ||I+A|| + \epsilon.$$
It follows that $\forall n \geq M$,
$$(I+A_n)^{-1} \geq \frac{1}{\max\{(1+\lambda_k(A_n)) : k \in \mathbb{N}\}} = \frac{1}{||I+A_n||} \geq \frac{1}{||I+A|| + \epsilon}.$$
Hence invoking Eq. (A.34) again, we obtain $\forall n \geq M$
$$||(I+A_n)^{-r} - (I+A)^{-r}||_p \leq r(||I+A|| + \epsilon)^{1-r} ||(I+A_n)^{-1} - (I+A)^{-1}||_p,$$
which implies that
$$\lim_{n\to\infty} ||(I+A_n)^{-r} - (I+A)^{-r}||_p = 0$$
by the previous limit, when $r = 1$.

(iv) By an induction argument as in step (ii), we then obtain that
$$\lim_{n\to\infty} ||(I+A_n)^{-r} - (I+A)^{-r}||_p = 0, \quad \forall r > 1. \tag{A.37}$$
This completes the proof. $\square$

**Lemma 18.** *Let $\mathcal{H}$ be a separable Hilbert space. Assume that $\{A_n\}_{n\in\mathbb{N}}$, $A$ are trace class operators on $\mathcal{H}$ such that $(I+A) > 0$, $(I+A_n) > 0$ $\forall n \in \mathbb{N}$. Assume that $||A_n - A||_{\text{tr}} = 0$ as $n \to \infty$. Then $A_n(I+A_n)^{-1}$ and $A(I+A)^{-1}$ are trace class operators and*
$$\lim_{n\to\infty} ||A_n(I+A_n)^{-1} - A(I+A)^{-1}||_{\text{tr}} = 0. \tag{A.38}$$

*Proof of Lemma 18.* It is obvious that, given that $A_n$ and $A$ are trace class operators, both $A_n(I+A_n)^{-1}$ and $A(I+A)^{-1}$ are trace class operators. We have
$$||A_n(I+A_n)^{-1} - A(I+A)^{-1}||_{\text{tr}} = ||(I+A_n)^{-1}A_n - A(I+A)^{-1}||_{\text{tr}}$$
$$= ||(I+A_n)^{-1}[A_n(I+A) - (I+A_n)A](I+A)^{-1}||_{\text{tr}}$$
$$= ||(I+A_n)^{-1}[A_n - A](I+A)^{-1}||_{\text{tr}} \leq ||(I+A_n)^{-1}||\, ||A_n - A||_{\text{tr}}\, ||(I+A)^{-1}||.$$

By the assumption $I + A > 0$, there exists $M_A > 0$ such that

$$\langle x, (I + A)x \rangle \geq M_A ||x||^2 \quad \forall x \in \mathcal{H}.$$

By the assumption $\lim_{n\to\infty} ||A_n - A||_{\text{tr}} = 0$, for any $\epsilon$ satisfying $0 < \epsilon < M_A$, there exists $N = N(\epsilon) \in \mathbb{N}$ such that $||A_n - A||_{\text{tr}} < \epsilon \ \forall n \geq N$. Then $\forall x \in \mathcal{H}$,

$$|\langle x, (A_n - A)x \rangle| \leq ||A_n - A|| \ ||x||^2 \leq ||A_n - A||_{\text{tr}} ||x||^2 \leq \epsilon ||x||^2.$$

It thus follows that $\forall x \in \mathcal{H}$,

$$\langle x, (I + A_n)x \rangle = \langle x, (I + A)x \rangle + \langle x, (A_n - A)x \rangle \geq (M_A - \epsilon)||x||^2.$$

Thus we have $I + A \geq M_A > 0$, $I + A_n \geq M_A - \epsilon > 0 \ \forall n \geq N = N(\epsilon)$, from which it follows that

$$||(I + A_n)^{-1}|| \leq \frac{1}{M_A - \epsilon} \forall N \geq N(\epsilon), \quad ||(I + A)^{-1}|| \leq \frac{1}{M_A}.$$

Combining this with the first inequality, we have

$$||A_n(I + A_n)^{-1} - A(I + A)^{-1}||_{\text{tr}} \leq \frac{1}{M_A(M_A - \epsilon)} ||A_n - A||_{\text{tr}} \forall n \geq N,$$

which implies that

$$\lim_{n\to\infty} ||A_n(I + A_n)^{-1} - A(I + A)^{-1}||_{\text{tr}} = 0.$$

This completes the proof. $\square$

**Lemma 19.** *Let $\mathcal{H}$ be a separable Hilbert space. Let $\{A_n\}_{n\in\mathbb{N}}$, $A$, $\{B_n\}_{n\in\mathbb{N}}$, $B$, be self-adjoint, trace class operators on $\mathcal{H}$, with $\lim_{n\to\infty} ||A_n - A||_{\text{tr}} = 0$, $\lim_{n\to\infty} ||B_n - B||_{\text{tr}} = 0$. Assume that $I + A > 0, I + B > 0, I + A_n > 0, I + B_n > 0 \ \forall n \in \mathbb{N}$. Then $(I + B_n)^{-1/2}(I + A_n)(I + B_n)^{-1/2} - I$ and $(I + B)^{-1/2}(I + A)(I + B)^{-1/2} - I$ are self-adjoint, trace class operators on $\mathcal{H}$ and*

$$\lim_{n\to\infty} ||(I + B_n)^{-1/2}(I + A_n)(I + B_n)^{-1/2} - (I + B)^{-1/2}(I + A)(I + B)^{-1/2}||_{\text{tr}}$$
$$= 0. \tag{A.39}$$

*Proof of Lemma 19.* We write

$$(I + B_n)^{-1/2}(I + A_n)(I + B_n)^{-1/2} = I - B_n(I + B_n)^{-1} - (I + B_n)^{-1/2}A_n(I + B_n)^{-1/2},$$

$$(I + B)^{-1/2}(I + A)(I + B)^{-1/2} = I - B(I + B)^{-1} - (I + B)^{-1/2}A(I + B)^{-1/2}.$$

It follows immediately that $[(I + B_n)^{-1/2}(I + A_n)(I + B_n)^{-1/2} - I]$ and $[(I + B)^{-1/2}(I + A)(I + B)^{-1/2} - I]$ are self-adjoint, trace class operators on $\mathcal{H}$.

By Lemma 18, we have

$$\lim_{n \to \infty} ||B_n(I + B_n)^{-1} - B(I + B)^{-1}||_{\text{tr}} = 0.$$

Consider next the difference between the third terms of the above two expressions

$$||(I + B_n)^{-1/2}A_n(I + B_n)^{-1/2} - (I + B)^{-1/2}A(I + B)^{-1/2}||_{\text{tr}}$$
$$\leq ||(I + B_n)^{-1/2}A_n(I + B_n)^{-1/2} - (I + B_n)^{-1/2}A(I + B_n)^{-1/2}||_{\text{tr}}$$
$$+ ||(I + B_n)^{-1/2}A(I + B_n)^{-1/2} - (I + B_n)^{-1/2}A(I + B)^{-1/2}||_{\text{tr}}$$
$$+ ||(I + B_n)^{-1/2}A(I + B)^{-1/2} - (I + B)^{-1/2}A(I + B)^{-1/2}||_{\text{tr}}. \quad (A.40)$$

By the assumption $I + A > 0$, $I + B > 0$, there exist constants $M_A > 0$, $M_B > 0$ such that $I + A \geq M_A$, $I + B \geq M_B$. As in the proof of Lemma 18, since $\lim_{n \to \infty} ||A_n - A|| = 0$, $\lim_{n \to \infty} ||B_n - B|| = 0$, for any $0 < \epsilon < \min\{M_A, M_B\}$, there exist $N_A = N_A(\epsilon) \in \mathbb{N}$, $N_B = N_B(\epsilon) \in \mathbb{N}$, such that

$$I + A_n \geq M_A - \epsilon, \ \forall n \geq N_A, \quad I + B_n \geq M_B - \epsilon \ \forall n \geq N_B.$$

The first term on the right hand side of the inequality in Eq. (A.40) is

$$||(I + B_n)^{-1/2}(A_n - A)(I + B_n)^{-1/2}||_{\text{tr}} \leq ||A_n - A||_{\text{tr}}||(I + B_n)^{-1/2}||^2$$
$$\leq \frac{1}{M_B - \epsilon}||A_n - A||_{\text{tr}} \ \forall n \geq N_B.$$

The second term is

$$||(I + B_n)^{-1/2}A[(I + B_n)^{-1/2} - (I + B)^{-1/2}]||_{\text{tr}}$$
$$\leq ||(I + B_n)^{-1/2}|| \ ||A||_{\text{tr}}||[(I + B_n)^{-1/2} - (I + B)^{-1/2}]||$$
$$\leq \frac{1}{\sqrt{M_B - \epsilon}}||A||_{\text{tr}}||[(I + B_n)^{-1/2} - (I + B)^{-1/2}]||.$$

Similarly, for the third term, we have

$$||(I + B_n)^{-1/2}A(I + B)^{-1/2} - (I + B)^{-1/2}A(I + B)^{-1/2}||_{\text{tr}}$$
$$\leq ||A(I + B)^{-1/2}||_{\text{tr}} ||[(I + B_n)^{-1/2} - (I + B)^{-1/2}]||.$$

By Theorem 30, we have

$$||(I + B_n)^{-1/2} - (I + B)^{-1/2}|| \leq ||(I + B_n)^{-1/2} - (I + B)^{-1/2}||_{\text{tr}} \to 0$$

as $n \to \infty$. The final result is obtained by combining all of the above inequalities. $\square$

**Lemma 20.** *Let $\mathcal{H}$ be a separable Hilbert space. Let $A, B, C : \mathcal{H} \to \mathcal{H}$ be self-adjoint, finite-rank operators such that $(I + A) > 0$, $(I + B) > 0$, $(I + C) > 0$. Then*

$$D_{2\alpha}^{(\alpha,\alpha)}[(I + A), (I + B)] \leq D_{2\alpha}^{(\alpha,\alpha)}[(I + A), (I + C)]$$
$$+ D_{2\alpha}^{(\alpha,\alpha)}[(I + C), (I + B)]. \quad \text{(A.41)}$$

*Proof of Lemma 20.* Since $A, B, C$ are all finite-rank operators, there exists a finite-dimensional subspace $\mathcal{H}_n \subset \mathcal{H}$, with $\dim(\mathcal{H}_n) = n$ for some $n \in \mathbb{N}$, such that $\text{range}(A) \subset \mathcal{H}_n, \text{range}(B) \subset \mathcal{H}_n$, and $\text{range}(C) \subset \mathcal{H}_n$. Let

$$A_n = A\big|_{\mathcal{H}_n} : \mathcal{H}_n \to \mathcal{H}_n, \quad B_n = B\big|_{\mathcal{H}_n} : \mathcal{H}_n \to \mathcal{H}_n, \quad C_n = C\big|_{\mathcal{H}_n} : \mathcal{H}_n \to \mathcal{H}_n.$$

Then $A_n, B_n, C_n$ are linear operators on the finite-dimensional space $\mathcal{H}_n$ and thus are represented by $n \times n$ matrices, which we denote by the same symbols. We have

$$(I + A_n)(I + B_n)^{-1} = (I + A_n)[I - B_n(I + B_n)^{-1}]$$
$$= I + A_n - B_n(I + B_n)^{-1} - A_n B_n(I + B_n)^{-1},$$
$$(I + A)(I + B)^{-1} = I + A - B(I + B)^{-1} - AB(I + B)^{-1},$$

where $A - B(I + B)^{-1} - AB(I + B)^{-1}$ is of finite rank, since both $A$ and $B$ are, with range in $\mathcal{H}_n$. It is clear that

$$[A - B(I + B)^{-1} - AB(I + B)^{-1}]\big|_{\mathcal{H}_n} = A_n - B_n(I + B_n)^{-1} - A_n B_n(I + B_n)^{-1}.$$

Thus the nonzero eigenvalues of $(I + A)(I + B)^{-1} - I = [A - B(I + B)^{-1} - AB(I + B)^{-1}]$ and $(I + A_n)(I + B_n)^{-1} - I = [A_n - B_n(I + B_n)^{-1} - A_n B_n(I + B_n)^{-1}]$ are the same. It follows that

$$
\begin{aligned}
&D_{2\alpha}^{(\alpha,\alpha)}[(I+A),(I+B)] \\
&= \frac{1}{\alpha^2} \log \det \left[ \frac{[(I+A)(I+B)^{-1}]^\alpha + [(I+A)(I+B)^{-1}]^{-\alpha}}{2} \right] \\
&= \frac{1}{\alpha^2} \log \det \left[ \frac{[(I+A_n)(I+B_n)^{-1}]^\alpha + [(I+A_n)(I+B_n)^{-1}]^{-\alpha}}{2} \right] \\
&= D_{2\alpha}^{(\alpha,\alpha)}[(I+A_n),(I+B_n)].
\end{aligned}
$$

Similarly, we have

$$
\begin{aligned}
D_{2\alpha}^{(\alpha,\alpha)}[(I+A),(I+C)] &= D_{2\alpha}^{(\alpha,\alpha)}[(I+A_n),(I+C_n)], \\
D_{2\alpha}^{(\alpha,\alpha)}[(I+C),(I+B)] &= D_{2\alpha}^{(\alpha,\alpha)}[(I+C_n),(I+B_n)].
\end{aligned}
$$

Applying the triangle inequality from the finite-dimensional setting [15], we get

$$
\begin{aligned}
D_{2\alpha}^{(\alpha,\alpha)}[(I+A_n),(I+B_n)] &\leq D_{2\alpha}^{(\alpha,\alpha)}[(I+A_n),(I+C_n)] \\
&\quad + D_{2\alpha}^{(\alpha,\alpha)}[(I+C_n),(I+B_n)].
\end{aligned}
$$

Together with the above expressions, this gives us the final result. $\square$

*Proof of Theorem 19 (Convergence in trace norm).* Let $I + \Lambda = (I+B)^{-1/2}(I+A)(I+B)^{-1/2}$ and $I + \Lambda_n = (I+B_n)^{-1/2}(I+A_n)(I+B_n)^{-1/2}$, with $\Lambda, \Lambda_n \in \mathrm{Sym}(\mathcal{H}) \cap \mathrm{Tr}(\mathcal{H})$.

By Lemma 19, we have $\lim_{n\to\infty} ||\Lambda_n - \Lambda||_{\mathrm{tr}} = 0$.

Thus by Theorem 30, we have

$$\lim_{n\to\infty} ||(I+\Lambda_n)^\alpha - (I+\Lambda)^\alpha||_{\mathrm{tr}} = 0 \ \ \forall \alpha \in \mathbb{R}.$$

By Definition 5, we have

$$D_{2\alpha}^{(\alpha,\alpha)}[(I+A_n),(I+B_n)] = \frac{1}{\alpha^2} \log \det \left[ \frac{(I+\Lambda_n)^\alpha + (I+\Lambda_n)^{-\alpha}}{2} \right].$$

Taking limit as $n \to \infty$ and applying the continuity of the Fredholm determinant in the trace norm (e.g. Theorem 3.5 in [22]), we obtain

$$
\begin{aligned}
&\lim_{n\to\infty} D_{2\alpha}^{(\alpha,\alpha)}[(I+A_n),(I+B_n)] \\
&= \frac{1}{\alpha^2} \log \det \left[ \frac{(I+\Lambda)^\alpha + (I+\Lambda)^{-\alpha}}{2} \right] = D_{2\alpha}^{(\alpha,\alpha)}[(I+A),(I+B)].
\end{aligned}
$$

This completes the proof. □

*Proof of Theorem 20 (Triangle inequality).* For a fixed $\gamma > 0$, we have

$$D_{2\alpha}^{(\alpha,\alpha)}[(A+\gamma I),(B+\gamma I)]$$
$$= \frac{1}{\alpha^2} \log \det{}_{\mathrm{X}} \left( \frac{[(A+\gamma I)(B+\gamma I)^{-1}]^\alpha + (A+\gamma I)(B+\gamma I)^{-1}]^{-\alpha}}{2} \right)$$
$$= \frac{1}{\alpha^2} \mathrm{logdet} \left( \frac{[(\frac{A}{\gamma}+I)(\frac{B}{\gamma}+I)^{-1}]^\alpha + (\frac{A}{\gamma}+I)(\frac{B}{\gamma}+I)^{-1}]^{-\alpha}}{2} \right),$$

which thus reduces to the case $\gamma = 1$. Thus it suffices for us to prove in triangle inequality for $\gamma = 1$.

Let $\{A_n\}_{n\in\mathbb{N}}$, $\{B_n\}_{n\in\mathbb{N}}$, and $\{C_n\}_{n\in\mathbb{N}}$ be sequences of finite-rank operators such that

$$\lim_{n\to\infty} ||A_n - A||_{\mathrm{tr}} = 0, \quad \lim_{n\to\infty} ||B_n - B||_{\mathrm{tr}} = 0, \quad \lim_{n\to\infty} ||C_n - C||_{\mathrm{tr}} = 0.$$

By Lemma 20, we have the triangle inequality

$$\sqrt{D_{2\alpha}^{(\alpha,\alpha)}[(I+A_n),(I+B_n)]} \le \sqrt{D_{2\alpha}^{(\alpha,\alpha)}[(I+A_n),(I+C_n)]}$$
$$+ \sqrt{D_{2\alpha}^{(\alpha,\alpha)}[(I+C_n),(I+B_n)]}.$$

Taking limits on both side as $n \to \infty$ and invoking Theorem 19, we then obtain

$$\sqrt{D_{2\alpha}^{(\alpha,\alpha)}[(I+A),(I+B)]} \le \sqrt{D_{2\alpha}^{(\alpha,\alpha)}[(I+A),(I+C)]}$$
$$+ \sqrt{D_{2\alpha}^{(\alpha,\alpha)}[(I+C),(I+B)]}.$$

This completes the proof of the theorem. □

*Proof of Theorem 4 (Metric property).* The case $\alpha = 0$ corresponds to the affine-invariant Riemannian distance on the Hilbert manifold $\Sigma(\mathcal{H})$ [18], which is still a metric when restricted to $\mathrm{PTr}(\mathcal{H})$.

Consider the case $\alpha > 0$. The positivity and symmetry of the divergence $D_{2\alpha}^{(\alpha,\alpha)}[(A+\gamma I),(B+\gamma I)]$ are from Theorems 1 and 13, respectively. The triangle inequality for $\sqrt{D_{2\alpha}^{(\alpha,\alpha)}[(A+\gamma I),(B+\gamma I)]}$ is from Theorem 20. Thus $\sqrt{D_{2\alpha}^{(\alpha,\alpha)}[(A+\gamma I),(B+\gamma I)]}$ is a metric on $\mathrm{PTr}(\mathcal{H})(\gamma)$. □

***Proof of Theorem 22 (Diagonalization).*** Consider first the case $\alpha > 0$. As in the proof of Theorem 20, it suffices for us to prove this theorem for the case $\gamma = 1$. Let $A = \sum_{j=1}^{\infty} \lambda_j(A)\phi_j \otimes \phi_j$ denote the spectral decomposition for $A$. For each $n \in \mathbb{N}$, define

$$A_n = \sum_{j=1}^{n} \lambda_j(A)\phi_j \otimes \phi_j.$$

Then $A_n$ is a finite-rank operator with the eigenvalues being the first $n$ eigenvalues of $A$ and $\lim_{n\to\infty} ||A_n - A||_{\mathrm{tr}} = 0$. In the same way, we construct a sequence of finite-rank operators $B_n$ with $\lim_{n\to\infty} ||B_n - B||_{\mathrm{tr}} = 0$. By construction, we also have

$$\lim_{n\to\infty} ||\mathrm{Eig}(A_n) - \mathrm{Eig}(A)||_{\mathrm{tr}} = 0, \quad \lim_{n\to\infty} ||\mathrm{Eig}(B_n) - \mathrm{Eig}(B)||_{\mathrm{tr}} = 0.$$

Thus by Theorem 19, we have

$$\lim_{n\to\infty} D_{2\alpha}^{(\alpha,\alpha)}[(\mathrm{Eig}(A_n) + I), (\mathrm{Eig}(B_n) + I)] = D_{2\alpha}^{(\alpha,\alpha)}[(\mathrm{Eig}(A) + I), (\mathrm{Eig}(B) + I)],$$

$$\lim_{n\to\infty} D_{2\alpha}^{(\alpha,\alpha)}[(A_n + I), (B_n + I)] = D_{2\alpha}^{(\alpha,\alpha)}[(A + I), (B + I)].$$

Since $A_n, B_n$ can be identified with finite-dimensional matrices, as in the proof of Lemma 16, we can apply the corresponding finite-dimensional result in [15] to obtain

$$D_{2\alpha}^{(\alpha,\alpha)}[(\mathrm{Eig}(A_n) + I), (\mathrm{Eig}(B_n) + I)] \leq D_{2\alpha}^{(\alpha,\alpha)}[(A_n + I), (B_n + I)].$$

Thus taking limits as $n \to \infty$ gives

$$D_{2\alpha}^{(\alpha,\alpha)}[(\mathrm{Eig}(A) + I), (\mathrm{Eig}(B) + I)] \leq D_{2\alpha}^{(\alpha,\alpha)}[(A + I), (B + I)].$$

Letting $\alpha \to 0$ on both sides of the above expression, we also obtain the result for the case $\alpha = 0$. This completes the proof of the theorem. □